\documentclass[12pt]{article}
\usepackage{amsmath}
\usepackage{amssymb}
\renewcommand\smallskip{\vskip\smallskipamount}
\renewcommand\medskip{\vskip\medskipamount}
\renewcommand\bigskip{\vskip\bigskipamount}

\begin{document}

\footnotetext{The first author acknowledges the support of NSF
Grant DMS-0654261. The second author acknowledges the support of
NSF Grant DMS-1007156 and a Sloan Research Fellowship.}

\begin{center}
\begin{large}
\textbf{On the Local Isometric Embedding in $\mathbb{R}^{3}$ of
Surfaces with Gaussian Curvature of Mixed Sign}
\end{large}

\bigskip\bigskip
QING HAN \text{ } \& \text{ } MARCUS KHURI
\bigskip\bigskip
\end{center}

\begin{abstract}
We study the old problem of isometrically embedding a
2-dimensional Riemannian manifold into Euclidean 3-space.  It is
shown that if the Gaussian curvature vanishes to finite order and
its zero set consists of two Lipschitz curves intersecting
transversely at a point, then local sufficiently smooth isometric
embeddings exist.
\end{abstract}

\bigskip

\begin{center}
\textbf{1.  Introduction}
\end{center}\setcounter{equation}{0}
\setcounter{section}{1}

  Does every smooth 2-dimensional Riemannian manifold admit a
smooth local isometric embedding into $\mathbb{R}^{3}$, or
heuristically, can every abstract surface be visualized at least
locally?  This natural question was first posed in 1873 by Schlaefli
[16], and remarkably has remained to a large extent unanswered.  It
is the purpose of this paper to provide a general sufficient
condition under which local embeddings exist.\par
  The local isometric embedding problem for surfaces is equivalent
to finding local solutions of a particular Monge-Amp\`{e}re
equation, usually referred to as the Darboux equation. The primary
difficulty in analyzing this equation arises from the fact that it
changes from elliptic to hyperbolic type, whenever the Gaussian
curvature of the given metric passes from positive to negative
curvature. Consequently the hypotheses of any result must take into
account the manner in which the Gaussian curvature, $K$, vanishes.
The classical results deal with the cases in which the curvature
does not vanish, or the metric is analytic.  It was not until
1985/86 that the first degenerate cases (when $K$ vanishes) were
treated, by Lin.  He showed the existence of sufficiently smooth
embeddings if the metric is sufficiently smooth and $K\geq 0$ [11],
or $K(0)=0$, $|\nabla K(0)|\neq 0$ [12]. Smooth embeddings of smooth
surfaces were obtained by Han, Hong, and Lin [5] when $K\leq 0$ and
$\nabla K$ possesses a certain nondegeneracy, and by Han [3] when
$K$ vanishes across a single smooth curve (see also [1], [2], and
[8] for related results). Lastly if $K=|\nabla K(0)|=0$,
$|\nabla^{2}K(0)|\neq 0$ then Khuri [9] has proven the existence of
sufficiently smooth embeddings for sufficiently smooth surfaces. For
more details on this problem and other related topics the reader is
referred to [4].  Here we will show\medskip

\textbf{Theorem 1.1.}  \textit{Let $g\in C^{m_{*}}$, $m_{*}\geq
36(N+10)$, be a Riemannian metric defined on a neighborhood of the
origin in the plane, with Gaussian curvature $K$ vanishing there to
finite order $N$.  If the zero set $K^{-1}(0)$ consists of two
$C^{m_{*}}$ curves intersecting transversely at the origin, then $g$
admits a $C^{m}$, $m\leq \frac{1}{12}m_{*}-N-24$, isometric
embedding into $\mathbb{R}^{3}$ on some neighborhood of the
origin.}\medskip

\textbf{Remark 1.1.} \textit{Our methods actually treat a slightly
more general situation in that the zero set $K^{-1}(0)$ may consist
of more than two curves intersecting transversely at the origin.
However, in this setting there should not be more than two regions
on which $K\geq 0$.}\medskip

  The embeddings produced by this theorem are referred to as
sufficiently smooth, since higher regularity of the metric implies
higher regularity for the embedding.  However, this theorem does not
guarantee that $C^{\infty}$ metrics give rise to $C^{\infty}$
embeddings, as the methods used here require the domain of existence
to shrink whenever higher regularity is demanded of the solution. On
the other hand, it is likely that techniques similar to those found
in [3] and [6] may lead to a $C^{\infty}$ version of Theorem 1.1. We
also point out that counterexamples to the existence of local
isometric embeddings have been found for metrics of low regularity,
by Pogorelov [15] when $g\in C^{2,1}$, and by Nadirashvili and Yuan
[13] who recently generalized this result. Moreover counterexamples
to the local solvability of smooth Monge-Amp\`{e}re equations have
been found in [10].  Yet it is still very much an open question,
whether or not there are any smooth (or sufficiently smooth)
counterexamples to the local isometric embedding problem.\par
  As mentioned above this problem is equivalent to finding local
solutions of a particular Monge-Amp\`{e}re equation.  To see this we
employ a standard method originally introduced by Weingarten [20].
That is, we search for a function $z(u,v)$ defined in a neighborhood
of the origin such that the new metric $g-dz^{2}$ is flat.  Note
that $g-dz^{2}$ will be Riemannian as long as $|\nabla_{g}z|<1$.
Since flat metrics are locally isometric to Euclidean space, there
exist two $C^{m}$ functions $x(u,v)$, $y(u,v)$ (if $g\in C^{m}$ and
$z\in C^{m+1}$, see [7]) such that $g-dz^{2}=dx^{2}+dy^{2}$. The map
$(u,v)\mapsto(x(u,v),y(u,v),z(u,v))$ then provides the desired
embedding.  Furthermore a straightforward calculation shows that
$g-dz^{2}$ is flat if and only if $z$ satisfies
\begin{equation}
\det\mathrm{Hess}_{g}z=K(\det g)(1-|\nabla_{g}z|^{2}).
\end{equation}
This Monge-Amp\`{e}re equation is the so called Darboux equation,
and the question of its local solvability is equivalent to the
local isometric embedding problem.\par
  It is trivial to construct approximate local solutions of (1.1),
and thus it is natural to employ an implicit function theorem to
prove existence.  Under the hypotheses of Theorem 1.1 the
linearization will be of mixed type, and thus we will necessarily
lose derivatives when it is inverted.  This suggests that we use a
version of the Nash-Moser implicit function theorem, which
essentially reduces the problem to a study of the linearized
equation.  We will show that the linearization has a particularly
nice canonical form, when an appropriate coordinate system is
chosen and certain perturbation terms which behave like quadratic
error in the Nash-Moser iteration are removed.  This was first
observed by Han in [3].  More precisely, the canonical form is
given by:
\begin{equation}
Lu=(aKu_{x})_{x}+bu_{yy}+cKu_{x}+du_{y}
\end{equation}
where $a,b>0$.  The significance of this particular structure is
that it explicitly illustrates how the Gaussian curvature affects
the type of the linearization.  Moreover it is also important that
the first order coefficient $cK$ vanishes whenever the principal
symbol changes type, as this leads to the so called Levi conditions
[14] in the hyperbolic regions, which facilitate the making of
estimates.  Under the assumptions of Theorem 1.1 on the Gauss
curvature, there are four separate regions of elliptic or hyperbolic
type for (1.2), each having a Lipschitz smooth boundary. We will
develop the appropriate existence and regularity theory for (1.2) in
each of these regions, and show that combined with a Nash-Moser
iteration this leads to a corresponding solution of the nonlinear
equation (1.1) in each region.  These separate solutions will then
be patched together to form a solution on a full neighborhood of the
origin.\par
  This paper is organized as follows.  In section $\S 2$ we obtain
the canonical form (1.2), and in sections $\S 3$ and $\S 4$ the
linear existence theory is established in the elliptic and
hyperbolic regions, respectively.  Lastly in section $\S 5$ we
employ a version of the Nash-Moser iteration to solve (1.1) in the
elliptic and hyperbolic regions separately, and also show how the
solutions obtained can be patched together to yield the desired
solution.

\begin{center}
\textbf{2. The Linearized Canonical Form}
\end{center}\setcounter{equation}{0}
\setcounter{section}{2}

  In this section we will bring the linearization of (1.1) into
the canonical form (1.2).  Before doing this, however, we must
specify at which function the linearization will be evaluated. For
this we need an appropriate approximate solution, $z_{0}$.  We
will then search for a solution of (1.1) in the form
\begin{equation*}
z=z_{0}+\varepsilon^{5}w,
\end{equation*}
where $\varepsilon>0$ is a small parameter.  Let $y=(y^{1},y^{2})$
be local coordinates in a neighborhood of the origin with
$g=g_{ij}dy^{i}dy^{j}$, then we are interested in solving
\begin{equation}
\Phi(w):=\det\nabla_{ij}z-K|g|(1-|\nabla_{g}z|^{2})=0
\end{equation}
where $|g|=\det g_{ij}$ and $\nabla_{ij}$ are covariant
derivatives with respect to these coordinates.  We choose
\begin{equation*}
z_{0}=\frac{1}{2}(y^{1})^{2}+\sum_{n=3}^{m_{*}}p_{n}(y)
\end{equation*}
where each $p_{n}$ is a homogeneous polynomial of degree $n$,
chosen so that
\begin{equation}
\partial^{\alpha}\Phi(0)=0,\text{ }\text{ }\text{ }\text{ }\text{
}|\alpha|\leq m_{*}-2.
\end{equation}
Here $\alpha=(\alpha_{1},\alpha_{2})$ is a multi-index, and
$m_{*}$ is as in Theorem 1.1.  Note that such a polynomial $z_{0}$
may be found in the usual way by following the proof of the
Cauchy-Kowalevski Theorem, since the line $y^{2}=0$ is
noncharacteristic for (2.1) as $\nabla_{11}z_{0}(0)=1$.\par
  Upon rescaling coordinates by $y^{i}=\varepsilon^{2} x^{i}$, the
linearization of $\Phi$ at a function $w$,
\begin{equation*}
\mathcal{L}(w)u=\frac{d}{dt}\Phi(w+tu)|_{t=0},
\end{equation*}
is given by
\begin{equation}
\varepsilon^{-1}\mathcal{L}(w)u
=a^{ij}u_{;ij}-2\varepsilon^{4}K|g|\langle\nabla_{g}z,\nabla_{g}u\rangle
\end{equation}
where $u_{;ij}$ denote covariant derivatives in $x^{i}$
coordinates, $\langle\cdot,\cdot\rangle$ is the inner product
associated with $g$, and
\begin{equation*}
(a^{ij})=\left(\begin{array}{cc}
\nabla_{22}z & -\nabla_{12}z \\
-\nabla_{12}z & \nabla_{11}z
\end{array}\right)
\end{equation*}
is the cofactor matrix of $\mathrm{Hess}_{g}z$.  Note that the
quantity $|g|^{-1}a^{ij}$ transforms like a contravariant
2-tensor. According to the assumptions of Theorem 1.1, $K^{-1}(0)$
divides a small neighborhood of the origin into domains
$\{\Omega_{\kappa}^{+}\}_{\kappa=1}^{\kappa_{0}}$ on which $K>0$,
and $\{\Omega_{\varrho}^{-}\}_{\varrho=1}^{\varrho_{0}}$ on which
$K<0$ (obviously $\kappa_{0}+\varrho_{0}=4$). The following lemma
gives the desired canonical form. We will denote the Sobolev space
of square integrable derivatives by $H^{m}$, with norm
$\parallel\cdot\parallel_{H^{m}}$.\medskip

\textbf{Lemma 2.1.}  \textit{Let $g\in C^{m_{*}}$ and $w\in
C^{\infty}$ with $|w|_{C^{3}}<1$.  Given a domain
$\Omega_{\kappa}^{+}$, or $\Omega_{\varrho}^{-}$, $\varrho=1,2$,
or $\Omega_{\varrho}^{-}$, $\varrho=3,4$, and given small
$\sigma,\delta>0$, there exists a local $C^{m_{*}-2}$ change of
coordinates $\xi^{i}=\xi^{i}(x)$ such that}
\begin{equation*}
\Omega_{\kappa}^{+}\cap B_{\sigma}(0)=\{(\xi^{1},\xi^{2})\mid
0<\xi^{2}<(\tan\delta)\xi^{1},\text{ }|\xi|<\sigma\},
\end{equation*}
\textit{or}
\begin{equation*}
\Omega_{\varrho}^{-}\cap B_{\sigma}(0)=\{(\xi^{1},\xi^{2})\mid
h(\xi^{1})<\xi^{2},\text{ }|\xi|<\sigma\},\text{ }\text{ }\text{
}\text{ }\varrho=1,2,
\end{equation*}
\textit{or}
\begin{equation*}
\Omega_{\varrho}^{-}\cap B_{\sigma}(0)=\{(\xi^{1},\xi^{2})\mid
h(\xi^{2})<\xi^{1},\text{ }|\xi|<\sigma\},\text{ }\text{ }\text{
}\text{ }\varrho=3,4,
\end{equation*}
\textit{for some Lipschitz function $h(\xi^{i})$ (not necessarily
the same for different regions) satisfying $h(0)=0$ and
$|h(\xi^{i})-|\xi^{i}||_{C^{1}}=O(\sigma)$. In this new coordinate
system the linearization takes the form}
\begin{equation*}
\varepsilon^{-1}\mathcal{L}(w)u=a^{22}L(w)u+(a^{22})^{-1}
\Phi(w)[\partial_{x^{1}}^{2}u-\partial_{x^{1}}\log
(a^{22}\sqrt{|g|})\partial_{x^{1}}u],
\end{equation*}
\textit{where}
\begin{equation*}
L(w)u=\partial_{\xi^{1}}(k\partial_{\xi^{1}}u)+\partial_{\xi^{2}}^{2}u+
c\partial_{\xi^{1}}u +d\partial_{\xi^{2}}u
\end{equation*}
\textit{with}
\begin{eqnarray*}
k&=&K\overline{k}(x,w,\nabla w,\nabla^{2}w,\nabla\xi),\\
c&=& K\overline{c}(x,w,\nabla
w,\nabla^{2}w,\nabla^{3}w,\nabla\xi,\nabla^{2}\xi)
+(a^{22})^{-2}\partial_{x^{1}}\Phi(w)\partial_{x^{1}}\xi^{1},\\
d&=&\varepsilon^{2}\overline{d}(x,w,\nabla w,\nabla^{2}w),
\end{eqnarray*}
\textit{for some $\overline{k},\overline{c},\overline{d}\in
C^{m_{*}-4}$ such that $\overline{k}>1/2$ if
$\varepsilon=\varepsilon(m)$ is chosen sufficiently small.
Moreover there exists a constant $C_{m}$ independent of
$\varepsilon$, $\delta$ such that}
\begin{equation}
\parallel\xi\parallel_{H^{m}}\leq\delta^{-1}C_{m}(1+\parallel
w\parallel_{H^{m+4}}),\text{ }\text{ }\text{ }\text{ }\text{
}m\leq m_{*}-2.
\end{equation}

\textbf{Remark 2.1.}  \textit{In the estimate (2.4), $\delta$ is
only relevant for the regions
$\{\Omega_{\kappa}^{+}\}_{\kappa=1}^{\kappa_{0}}$. Furthermore,
since the curvature $K$ vanishes at least to second order, it may be
possible to eliminate the role of $\delta$ in the arguments of the
next section.}\medskip

\textit{Proof.} We may choose an initial coordinate system
$x=(x^{1},x^{2})$ so that each of the elliptic and hyperbolic
regions $\Omega_{\kappa}^{+}$, $\Omega_{\varrho}^{-}$ are sector
domains, that is, each occupies the region between two lines
passing through the origin. Furthermore we may assume that
$\Omega_{1}^{-}$ ($\Omega_{2}^{-}$) contains the positive
(negative) $x^{2}$-axis. Note that according to the hypotheses of
Theorem 1.1, $\partial\Omega_{\kappa}^{+}-\{(0,0)\}$ and
$\partial\Omega_{\varrho}^{-}-\{(0,0)\}$ are both $C^{m_{*}}$
smooth, so that this initial transformation is also $C^{m_{*}}$.
The approximate solution $z_{0}$ is chosen with respect to this
initial coordinate system (recall that
$y^{i}=\varepsilon^{2}x^{i}$), and therefore
\begin{equation*}
a^{22}>0,\text{ }\text{ }\text{ }\text{
}a^{12}=O(\varepsilon^{2}).
\end{equation*}
It is now an easy exercise in linear algebra to show that for each
domain $\Omega_{\kappa}^{+}$, or $\Omega_{\varrho}^{-}$,
$\varrho=1,2$, or $\Omega_{\varrho}^{-}$, $\varrho=3,4$, there
exists a linear change of coordinates
$\overline{x}=(\overline{x}^{1},\overline{x}^{2})$ such that
\begin{equation*}
\Omega_{\kappa}^{+}=\{(\overline{x}^{1},\overline{x}^{2})\mid
0<\overline{x}^{2}<\overline{x}^{1},\text{ }|\overline{x}|<1\},
\end{equation*}
or
\begin{equation*}
\Omega_{\varrho}^{-}=\{(\overline{x}^{1},\overline{x}^{2})\mid
|\overline{x}^{1}|<\overline{x}^{2}<1\},\text{ }\text{ }\text{
}\text{ }\varrho=1,2,
\end{equation*}
or
\begin{equation*}
\Omega_{\varrho}^{-}=\{(\overline{x}^{1},\overline{x}^{2})\mid
|\overline{x}^{2}|<\overline{x}^{1}<1\},\text{ }\text{ }\text{
}\text{ }\varrho=3,4,
\end{equation*}
and such that
\begin{equation*}
\partial_{\overline{y}^{1}}^{2}(y^{1})^{2}>0,\text{ }\text{
}\text{ }\text{ }
\partial_{\overline{y}^{1}}\partial_{\overline{y}^{2}}(y^{1})^{2}=0.
\end{equation*}
Here $\overline{y}^{i}=\varepsilon^{2}\overline{x}^{i}$. It
follows that $a^{22}>0$ and $a^{12}=O(\varepsilon^{2})$ are
preserved under this linear change of coordinates. For convenience
we will still denote $\overline{y}^{i}$ by $y^{i}$ and
$\overline{x}^{i}$ by $x^{i}$.\par

We may write (2.3) as
\begin{equation*}
L_{1}(w)u=a_{1}^{ij}u_{x^{i}x^{j}}
+a_{1}^{i}u_{x^{i}}:=\varepsilon^{-1}\mathcal{L}(w)u,
\end{equation*}
where $a_{1}^{ij}=a^{ij}$,
\begin{equation}
a_{1}^{l}=-\varepsilon^{2}(a^{ij}\Gamma_{ij}^{l}+2K|g|z^{l})
\end{equation}
with $z^{l}=g^{li}z_{y^{i}}$, and $\Gamma_{ij}^{l}$ are
Christoffel symbols in $y^{i}$ coordinates.  According to (2.1)
\begin{equation*}
a_{1}^{11}=\nabla_{22}z=(a^{22})^{-1}[K|g|(1-|\nabla_{g}z|^{2})+(\nabla_{12}z)^{2}
+\Phi(w)].
\end{equation*}
We then set
\begin{equation*}
L_{2}(w)u=a_{2}^{ij}u_{x^{i}x^{j}}+a_{2}^{i}u_{x^{i}}:=
L_{1}(w)u-(a^{22})^{-1}\Phi(w)u_{x^{1}x^{1}},
\end{equation*}
that is $a_{1}^{ij}=a_{2}^{ij}$ and $a_{1}^{i}=a_{2}^{i}$ except
for
\begin{equation}
a_{2}^{11}=(a^{22})^{-1}[K|g|(1-|\nabla_{g}z|^{2})+(\nabla_{12}z)^{2}].
\end{equation}
Also let
\begin{equation*}
L_{3}(w)u=a_{3}^{ij}u_{x^{i}x^{j}}+a_{3}^{i}u_{x^{i}}:=(a^{22})^{-1}L_{2}(w)u.
\end{equation*}\par
  We now define the desired change of coordinates by
\begin{equation*}
\xi^{1}=\xi^{1}(x^{1},x^{2}),\text{ }\text{ }\text{ }\text{
}\xi^{2}=x^{2},
\end{equation*}
with
\begin{equation}
a^{12}\xi^{1}_{x^{1}}+a^{22}\xi^{1}_{x^{2}}=0,
\end{equation}
so that if
\begin{equation*}
L_{4}(w)u=a_{4}^{ij}u_{\xi^{i}\xi^{j}}+a_{4}^{i}u_{\xi^{i}}:=L_{3}(w)u
\end{equation*}
then
\begin{equation*}
a_{4}^{12}=a_{3}^{ij}\xi^{1}_{x^{i}}\xi^{2}_{x^{j}}=0.
\end{equation*}
In order to obtain the correct expression in these new coordinates
for the domains $\Omega_{\kappa}^{+}$, or $\Omega_{\varrho}^{-}$,
$\varrho=1,2$, or $\Omega_{\varrho}^{-}$, $\varrho=3,4$, we impose
the initial conditions
\begin{equation}
\xi^{1}(x^{1},x^{1})=(\tan\delta)^{-1}x^{1},\text{ }\text{ }\text{
or }\text{ }\text{ }\xi^{1}(x^{1},0)=x^{1},\text{ }\text{ }\text{
or }\text{ }\text{ }\xi^{1}(x^{1},0)=x^{1},
\end{equation}
respectively.  Note that since the curves
$x^{1}\mapsto(x^{1},x^{1})$ and $x^{1}\mapsto(x^{1},0)$ are
noncharacteristic for (2.7), equation (2.7) with initial condition
(2.8) has a unique $C^{m_{*}-2}$ solution on some neighborhood of
the origin.  Furthermore, standard methods for first order
equations combined with the Gagliardo-Nirenberg inequalities
(Lemma 5.2 below) yields (2.4).  Note also that
$(\xi^{1},\xi^{2})$ forms a new coordinate system near the origin
since
\begin{equation*}
(\tan\delta)^{-1}=\xi^{1}_{x^{1}}(0,0)+\xi^{1}_{x^{2}}(0,0)=
\xi^{1}_{x^{1}}(0,0)\left(1-\frac{a^{12}}{a^{22}}(0,0)\right),
\end{equation*}
or
\begin{equation*}
\xi^{1}_{x^{1}}(0,0)=1,
\end{equation*}
according to the respective initial conditions given by (2.8).\par
  We shall now calculate the coefficients of $L_{4}(w)$.
Observe that (2.6) yields
\begin{eqnarray*}
a_{4}^{11}=a_{3}^{ij}\xi^{1}_{x^{i}}\xi^{1}_{x^{j}}\!\!\!
&=&\!\!\!(\nabla_{11}z)^{-2}[K|g|(1-|\nabla_{g}z|^{2})+(\nabla_{12}z)^{2}](\xi^{1}_{x^{1}})^{2}\\
& &\!\!\!
-2(\nabla_{11}z)^{-1}(\nabla_{12}z)\xi^{1}_{x^{1}}\xi^{1}_{x^{2}}+(\xi_{x^{2}}^{1})^{2},
\end{eqnarray*}
but (2.7) gives
\begin{equation*}
\xi^{1}_{x^{2}}=(\nabla_{11}z)^{-1}(\nabla_{12}z)\xi^{1}_{x^{1}}
\end{equation*}
so
\begin{equation*}
a_{4}^{11}=(a^{22})^{-2}K|g|(1-|\nabla_{g}z|^{2})(\xi^{1}_{x^{1}})^{2}.
\end{equation*}\par
  Next we examine $a_{4}^{1}$.  By (2.7)
\begin{equation*}
\xi^{1}_{x^{1}x^{2}}=-(a^{12}_{3})_{x^{1}}\xi^{1}_{x^{1}}-a^{12}_{3}\xi^{1}_{x^{1}x^{1}},\text{
}\text{ }\text{ }\text{
}\xi^{1}_{x^{2}x^{2}}=-(a^{12}_{3})_{x^{2}}\xi^{1}_{x^{1}}-a^{12}_{3}\xi^{1}_{x^{1}x^{2}},
\end{equation*}
so (2.6) produces
\begin{eqnarray}
a^{1}_{4}\!\!&=&\!\!a^{ij}_{3}\xi^{1}_{x^{i}x^{j}}+
a^{i}_{3}\xi^{1}_{x^{i}}\\
&=&\!\!(a^{22})^{-2}K|g|(1-|\nabla_{g}z|^{2})\xi^{1}_{x^{1}x^{1}}
-\left[\left(\frac{a^{12}}{a^{22}}\right)\left(\frac{a^{12}}{a^{22}}\right)_{x^{1}}
+\left(\frac{a^{12}}{a^{22}}\right)_{x^{2}}\right]\xi^{1}_{x^{1}}+a^{i}_{3}
\xi^{1}_{x^{i}}.\nonumber
\end{eqnarray}
Calculating the second term on the right-hand side yields,
\begin{eqnarray*}
(a^{22})^{2}\left[\left(\frac{a^{12}}{a^{22}}\right)\left(\frac{a^{12}}{a^{22}}\right)_{x^{1}}
+\left(\frac{a^{12}}{a^{22}}\right)_{x^{2}}\right]\!\!&=&\!\!
a^{12}a^{12}_{x^{1}}-(a^{22})^{-1}(a^{12})^{2}a^{22}_{x^{1}}+a^{22}a^{12}_{x^{2}}
-a^{12}a^{22}_{x^{2}}\\
&=&\!\!
a^{12}a^{12}_{x^{1}}-a^{11}a^{22}_{x^{1}}+a^{22}a^{12}_{x^{2}}-a^{12}a^{22}_{x^{2}}\\
& &\!\!
+(a^{22})^{-1}a^{22}_{x^{1}}(\det a^{ij})\\
&=&\!\!-a^{12}_{x^{1}}a^{12}+a^{11}_{x^{1}}a^{22}+a^{22}a^{12}_{x^{2}}
-a^{12}a^{22}_{x^{2}}\\
& &\!\!+(a^{22})^{-1}a^{22}_{x^{1}}(\det a^{ij})-(\det
a^{ij})_{x^{1}}.
\end{eqnarray*}
Therefore (2.5), (2.7), and (2.9) imply that
\begin{eqnarray}
a^{22}a^{1}_{4}\!\!&=&\!\!-[a_{x^{j}}^{ij}+\varepsilon^{2}(a^{lj}\Gamma_{lj}^{i}
+2K|g|z^{i})]\xi^{1}_{x^{i}}+((a^{22})^{-1}\det
a^{ij})_{x^{1}}\xi^{1}_{x^{1}}\\
&
&+(a^{22})^{-1}K|g|(1-|\nabla_{g}z|^{2})\xi^{1}_{x^{1}x^{1}}.\nonumber
\end{eqnarray}
A computation shows
\begin{eqnarray*}
&
&\!\!a^{11}_{y^{1}}+a^{12}_{y^{2}}+a^{lj}\Gamma_{lj}^{1}\nonumber\\
&=&\!\!
-\Gamma_{j2}^{j}z_{y^{1}y^{2}}+\Gamma_{j1}^{j}z_{y^{2}y^{2}}\nonumber\\
& &\!\!+(\Gamma_{12,y^{2}}^{i}
-\Gamma_{22,y^{1}}^{i}-\Gamma_{11}^{1}\Gamma_{22}^{i}+2\Gamma_{12}^{1}\Gamma_{12}^{i}
-\Gamma_{22}^{1}\Gamma_{11}^{i})z_{y^{i}}\\
&=&\!\! \Gamma_{j2}^{j}a^{12}+\Gamma_{j1}^{j}a^{11}\nonumber\\
& &\!\!+(\Gamma_{12,y^{2}}^{i}
-\Gamma_{22,y^{1}}^{i}-\Gamma_{11}^{1}\Gamma_{22}^{i}+2\Gamma_{12}^{1}\Gamma_{12}^{i}
-\Gamma_{22}^{1}\Gamma_{11}^{i}-\Gamma_{j2}^{j}\Gamma_{12}^{i}
+\Gamma_{j1}^{j}\Gamma_{22}^{i})z_{y^{i}}.\nonumber
\end{eqnarray*}
However, we see that the coefficient of $z_{y^{i}}$ is in fact a
curvature term.  More precisely, if we denote it by $\chi^{i}$
then
\begin{equation*}
\chi^{i}=\Gamma_{12,y^{2}}^{i}-\Gamma_{22,y^{1}}^{i}
+\Gamma_{12}^{j}\Gamma_{j2}^{i}-\Gamma_{22}^{j}\Gamma_{j1}^{i}
=-R^{i}_{212}=-g^{i1}|g|K
\end{equation*}
where $R^{i}_{jkl}$ is the Riemann tensor for $g$ in $y^{i}$
coordinates (recall that $\Gamma_{lj}^{i}$ are Christoffel symbols
in $y^{i}$ coordinates).  A similar calculation shows that
\begin{equation*}
a^{12}_{y^{1}}+a^{22}_{y^{2}}+a^{lj}\Gamma_{lj}^{2}
=\Gamma_{j1}^{j}a^{12}+\Gamma_{j2}^{j}a^{22}-g^{i2}|g|Kz_{y^{i}}.
\end{equation*}
Therefore after solving for $\xi^{1}_{x^{2}}$ in (2.7), (2.10)
becomes
\begin{eqnarray*}
a^{22}a_{4}^{1}&=&(a^{22})^{-1}K|g|(1-|\nabla_{g}z|^{2})\xi^{1}_{x^{1}x^{1}}
-\varepsilon^{2}K|g|z^{i}\xi^{1}_{x^{i}}\\
& &-[\varepsilon^{2}\Gamma_{j1}^{j}(a^{22})^{-1}\det
a^{ij}-((a^{22})^{-1}\det a^{ij})_{x^{1}}]\xi^{1}_{x^{1}}.
\end{eqnarray*}
It now follows from
\begin{equation*}
\det a^{ij}=\Phi(w)+K|g|(1-|\nabla_{g}z|^{2})
\end{equation*}
that we have
\begin{eqnarray*}
a^{22}a_{4}^{1}&=&-\varepsilon^{2}K|g|z^{i}\xi^{1}_{x^{i}}
+\partial_{x^{1}}[(a^{22})^{-1}K|g|(1-|\nabla_{g}z|^{2})\xi^{1}_{x^{1}}]\\
&
&-[\varepsilon^{2}\Gamma_{j1}^{j}(a^{22})^{-1}K|g|(1-|\nabla_{g}z|^{2})
+\varepsilon^{2}\Gamma_{j1}^{j}(a^{22})^{-1}\Phi(w)
-((a^{22})^{-1}\Phi(w))_{x^{1}}]\xi^{1}_{x^{1}}.
\end{eqnarray*}\par
  Lastly, it is trivial to calculate the remaining coefficients:
\begin{equation*}
a_{4}^{22}=1,\text{ }\text{ }\text{ }\text{ }a_{4}^{2}
=-\varepsilon^{2}(a^{22})^{-1}(a^{ij}\Gamma_{ij}^{2}+2K|g|z^{2}).
\end{equation*}
Then by defining
\begin{equation*}
L(w)u:=L_{4}(w)u+(\varepsilon^{2}\Gamma_{j1}^{j}+\partial_{x^{1}}\log
a^{22})(a^{22})^{-2}\Phi(w)\xi^{1}_{x^{1}}u_{\xi^{1}}
\end{equation*}
and recalling that
$\Gamma_{j1}^{j}=\frac{1}{2}\partial_{y^{1}}\log|g|$, we obtain
the desired result.  Q.E.D.

\begin{center}
\textbf{3. Linear Theory in the Elliptic Regions}
\end{center}\setcounter{equation}{0}
\setcounter{section}{3}

  In light of Lemma 2.1, it will be sufficient for our purposes to
study the question of existence and regularity for the operator
$L(w)$, instead of the pure linearization $\mathcal{L}(w)$.  More
precisely, in this section we will study the Dirichlet problem for
a modified version of $L(w)$ in an elliptic region.  First note
that by using polar coordinates $\xi^{1}=r\cos\theta$,
$\xi^{2}=r\sin\theta$, we can transform the elliptic region
$\Omega_{\kappa}^{+}\cap B_{\sigma}(0)$ of Lemma 2.1 into a
rectangle
\begin{equation*}
\Omega=\{(r,\theta)\mid 0<r<\sigma,\text{ }0<\theta<\delta\}.
\end{equation*}
Under these coordinates we find that
\begin{equation*}
L(w)u=\mathcal{K}u_{rr}+\mathcal{A}u_{r\theta}+\mathcal{B}u_{\theta\theta}
+\mathcal{C}u_{r}+\mathcal{D}u_{\theta}
\end{equation*}
where
\begin{eqnarray*}
\mathcal{K}&=&k\cos^{2}\theta+\sin^{2}\theta,\\
\mathcal{A}&=&2(1-k)\frac{\sin\theta\cos\theta}{r},\\
\mathcal{B}&=&k\frac{\sin^{2}\theta}{r^{2}}+\frac{\cos^{2}\theta}{r^{2}},\\
\mathcal{C}&=&k\frac{\sin^{2}\theta}{r}+\frac{\cos^{2}\theta}{r}+
(c+\partial_{\xi^{1}}k)\cos\theta+d\sin\theta,\\
\mathcal{D}&=&2(k-1)\frac{\sin\theta\cos\theta}{r^{2}}
-(c+\partial_{\xi^{1}}k)\frac{\sin\theta}{r}
+d\frac{\cos\theta}{r}.
\end{eqnarray*}
It will be convenient to cut-off these coefficients away from the
origin.  So let $\varphi\in C^{\infty}([0,\infty))$ be a
nonnegative cut-off function with
\begin{equation*}
\varphi(r)=\begin{cases}
1 & \text{if $0<r<\frac{1}{2}\sigma$,}\\
0 & \text{if $\sigma<r$},
\end{cases}
\end{equation*}
and define
\begin{equation*}
Lu=\overline{K}u_{rr}+\overline{A}u_{r\theta}+\overline{B}u_{\theta\theta}
+\overline{C}u_{r}+\overline{D}u_{\theta}
\end{equation*}
where
\begin{equation*}
\overline{K}=\varphi^{2}\mathcal{K},\text{ }\text{ }
\overline{A}=\varphi\mathcal{A},\text{ }\text{ }
\overline{B}=\mathcal{B},\text{ }\text{
}\overline{C}=\varphi\mathcal{C},\text{ }\text{
}\overline{D}=\varphi\mathcal{D}.
\end{equation*}\par
  We will study the boundary value problem:
\begin{equation}
Lu=f\text{ }\text{ }\text{ in }\text{ }\text{ }\Omega,\text{
}\text{ }\text{ }u(r,0)=u(r,\delta)=0,\text{ }\text{ }\text{
}\partial_{r}^{s}u(0,\theta)=0,\text{ }\text{ }\text{ }0\leq s\leq
s_{0},
\end{equation}
for some large integer $s_{0}$. The motivation for considering
this problem stems from our method for constructing solutions to
the nonlinear problem (2.1) (see section $\S 5$).  Namely, we
shall construct solutions in the elliptic and hyperbolic regions
separately, and then show that they can be patched together.  This
requires certain compatibility conditions at the origin, and the
boundary conditions of (3.1) guarantee that they will be
satisfied.  Of course, a necessary condition for solving (3.1) is
that $f$ must also vanish to a corresponding high order at the
origin.  It is therefore convenient to introduce the following
weighted Sobolev spaces which control the amount of vanishing.
Define norms
\begin{equation*}
\parallel u\parallel_{(m,l,\gamma)}^{2}=\int_{\Omega}\sum_{0\leq s\leq m,\text{
} 0\leq t\leq l \atop s+t\leq\mathrm{max}(m,l)}\lambda^{-s}
r^{-\gamma+2s}(\partial_{r}^{s}\partial_{\theta}^{t}u)^{2}
\end{equation*}
where $\lambda,\gamma>0$ are large parameters, and let
$H^{(m,l,\gamma)}(\Omega)$ be the closure of
$\overline{C}^{\infty}(\Omega)$ with respect to this norm, where
$\overline{C}^{\infty}(\Omega)$ is the space of smooth functions
which vanish in a neighborhood of $r=0$.  We will always denote
the traditional Sobolev spaces having square integrable
derivatives up to and including order $m$ by $H^{m}(\Omega)$, with
norm $\parallel\cdot\parallel_{m}$.  The following simple lemma
describes the boundary behavior exhibited by elements of the
weighted spaces.\medskip

\textbf{Lemma 3.1.}  \textit{Suppose that $u\in
H^{(m,l,\gamma)}(\Omega)\cap C^{\max(m,l)-1}(\Omega)$, then for
any $0<r_{0}<\sigma$ we have}
\begin{equation*}
\int_{r=r_{0}}(\partial_{r}^{s}\partial_{\theta}^{t}u)^{2}\leq
r_{0}^{\gamma-2(s+1)}C\parallel
u\parallel_{(m,l,\gamma)}^{2},\text{ }\text{ }s\leq m-1,\text{
}\text{ }t\leq l-1,\text{ }\text{ }s+t\leq\max(m,l)-1,
\end{equation*}
\textit{where the constant $C$ depends only on
$\sigma-r_{0}$.}\medskip

\textit{Proof.}  When $s\leq m-1$, $t\leq l-1$,
$s+t\leq\max(m,l)-1$ we have
$r^{-\gamma/2+s+1}\partial_{r}^{s}\partial_{\theta}^{t}u\in
H^{1}(\Omega)\cap C^{0}(\Omega)$.  The desired result now follows
from the standard trace theorem for Sobolev spaces.
Q.E.D.\medskip

  In analogy with the theory of strictly elliptic equations in a
sector domain such as $\Omega$, the regularity of a solution to
(3.1) will depend on the size of the angle forming the domain.
More precisely, smaller angles yield higher regularity.  According
to Lemma 2.1 we are free to choose the angle $\delta$ arbitrarily
small, with the only price being paid with the blow-up of estimate
(2.4).  This blow-up, however, can be controlled in the context of
equation (3.1) by taking $\varepsilon=\varepsilon(\delta)$ to be
sufficiently small, since whenever $\xi(x)$ or its derivatives
appear in the coefficients of the operator $L$, they are always
multiplied by $\varepsilon$.  These considerations lead to
existence and regularity for (3.1), and the first step needed to
establish such a result is the following basic estimate.  Define
\begin{equation*}
a_{\lambda,\gamma}(r,\theta)=\frac{\lambda\theta^{2}-1}{r^{\gamma}},
\end{equation*}
and as in Lemma 2.1 let $w\in C^{\infty}$ throughout this
section.\medskip

\textbf{Lemma 3.2.}  \textit{Suppose that $|w|_{C^{4}}<1$ and let
$u\in H^{(2,1,\gamma+2)}(\Omega)\cap C^{2}(\Omega)$ with
$u(r,0)=u(r,\delta)=0$.  If $\delta=\delta(\lambda)$ and
$\varepsilon=\varepsilon(\delta)$ are sufficiently small, then}
\begin{equation*}
\int_{\Omega}a_{\lambda,\gamma-2}uLu\geq C\int_{\Omega}\lambda
r^{-\gamma}u^{2}+r^{-\gamma+2}(\varphi\sin\theta
u_{r}+r^{-1}\cos\theta u_{\theta})^{2}
\end{equation*}
\textit{for some constant $C>0$ independent of $\lambda$,
$\delta$, $\varepsilon$, and $w$.}\medskip

\textit{Proof.}  Let $0<r_{0}<\sigma$ and set
$\Omega_{r_{0}}=\Omega\cap\{(r,\theta)\mid r_{0}<r<\sigma\}$.
Then for any $a\in C^{\infty}(\Omega_{r_{0}})$, integration by
parts yields
\begin{eqnarray*}
\int_{\Omega_{r_{0}}}auLu
&=&\int_{\Omega_{r_{0}}}-a(\overline{K}u_{r}^{2}+\overline{A}u_{r}u_{\theta}
+\overline{B}u_{\theta}^{2})\\
&
&+\frac{1}{2}\int_{\Omega_{r_{0}}}
[(a\overline{K})_{rr}+(a\overline{A})_{r\theta}+(a\overline{B})_{\theta\theta}
-(a\overline{C})_{r}-(a\overline{D})_{\theta}]u^{2}\\
&
&+\int_{\partial\Omega_{r_{0}}}a(\overline{K}uu_{r}\nu_{1}+\overline{A}uu_{r}\nu_{2}
+\overline{B}uu_{\theta}\nu_{2})\\
& &+\frac{1}{2}\int_{\partial\Omega_{r_{0}}}
[-(a\overline{K})_{r}\nu_{1}-(a\overline{A})_{\theta}\nu_{1}
-(a\overline{B})_{\theta}\nu_{2}+a\overline{C}\nu_{1}+a\overline{D}\nu_{2}]u^{2},
\end{eqnarray*}
where $(\nu_{1},\nu_{2})$ denotes the unit outer normal to
$\partial\Omega_{r_{0}}$.  By choosing $a=a_{\lambda,\gamma-2}$
and observing that
\begin{equation*}
|(a_{\lambda,\gamma-2}\overline{K})_{rr}|+|(a_{\lambda,\gamma-2}\overline{A})_{r\theta}|
+|(a_{\lambda,\gamma-2}\overline{C})_{r}|+
|(a_{\lambda,\gamma-2}\overline{D})_{\theta}|=O(r^{-\gamma}),
\end{equation*}
\begin{equation*}
(a_{\lambda,\gamma-2}\overline{B})_{\theta\theta}=2\lambda
r^{-\gamma}\cos^{2}\theta(1+O(\theta+\lambda^{-1})),
\end{equation*}
the desired result follows since all boundary terms vanish
according to Lemma 3.1 (after letting $r_{0}\rightarrow 0$).  Note
that $\varepsilon$ is chosen small depending on $\delta$, in order
to control the blow-up (implied by the estimate (2.4)) found in
the coefficients of $L$.  Q.E.D.\medskip

  This lemma is the main tool used to establish the basic
existence result of the next theorem.  Let
$\widehat{C}^{\infty}(\Omega)$ denote the space of
$\overline{C}^{\infty}(\Omega)$ functions $v$ satisfying
$v(r,0)=v(r,\delta)=0$.  Given $f\in H^{(m,1,\gamma)}(\Omega)$, we
will refer to a function $u\in H^{(m,1,\gamma)}(\Omega)$ as a weak
solution of (3.1) if
\begin{equation}
(u,L^{*}v)=(f,v)\text{ }\text{ }\text{ all }\text{ }\text{
}v\in\widehat{C}^{\infty}(\Omega),
\end{equation}
where $(\cdot,\cdot)$ is the $L^{2}(\Omega)$ inner product and
$L^{*}$ is the formal adjoint of $L$.\medskip

\textbf{Theorem 3.1.}  \textit{Suppose that $g\in C^{m_{*}}$,
$|w|_{C^{4}}<1$, and $f\in H^{(m,1,\gamma)}(\Omega)$.  If $m\leq
m_{*}-4$ and $\delta=\delta(m)$,
$\varepsilon=\varepsilon(m,\delta)$ are sufficiently small, then
there exists a weak solution $u\in H^{(m,1,\gamma)}(\Omega)$ of
(3.1).}\medskip

\textit{Proof.} Given $v\in\widehat{C}^{\infty}(\Omega)$, let
$\zeta\in H^{(m,\infty,\gamma+2)}(\Omega)\cap C^{\infty}(\Omega)$
be the unique solution of the ODE:
\begin{eqnarray}
&
&\sum_{s=0}^{m}\lambda^{-s}(-1)^{s}\partial_{r}^{s}(a_{\lambda,\gamma-2(s-1)}
\partial_{r}^{s}\zeta)=v,\\
& &\zeta(r,0)=\zeta(r,\delta)=0,\text{ }\text{ }\text{
}\partial_{r}^{s}\zeta(\sigma,\theta)=0,\text{ }\text{ }\text{
}0\leq s\leq m-1,\nonumber\\
&
&\int_{r=r_{0}}(\partial_{r}^{s}\partial_{\theta}^{l}\zeta)^{2}\leq
r_{0}^{\gamma-2s}C,\text{ }\text{ }\text{ }0\leq s\leq 2m-1,\text{
}\text{ }\text{ }0\leq l<\infty.\nonumber
\end{eqnarray}
Here $r_{0}>0$ is assumed to be sufficiently small, and $C>0$ is a
constant depending on $m$, $\lambda$, $\gamma$, and $v$.  The
proof that such a solution exists may be found in appendix A,
section $\S 6$.\par
  Our first goal is to establish an estimate of the form
\begin{equation}
(L\zeta,r^{4}\sum_{s=0}^{m}\lambda^{-s}(-1)^{s}\partial_{r}^{s}(
a_{\lambda,\gamma-2(s-1)}\partial_{r}^{s}\zeta))\geq
C\parallel\zeta\parallel_{(m,1,\gamma)}^{2}.
\end{equation}
The boundary conditions of (3.3) allow us to integrate by parts in
a manner similar to the proof of Lemma 3.2 to find,
\begin{eqnarray}
& &(L\zeta,r^{4}\sum_{s=0}^{m}\lambda^{-s}(-1)^{s}\partial_{r}^{s}
(a_{\lambda,\gamma-2(s-1)}\partial_{r}^{s}\zeta))\\
&=&
\sum_{s=0}^{m}\lambda^{-s} ([\partial_{r}^{s},L]\zeta
+L\partial_{r}^{s}\zeta,r^{4}a_{\lambda,\gamma-2(s-1)}\partial_{r}^{s}\zeta)\nonumber\\
& &+\sum_{s=0}^{m}\lambda^{-s}\left(\sum_{l<s}
\left(\begin{array}{c}
s \\
l
\end{array}\right)\partial_{r}^{s-l}r^{4}\partial_{r}^{l}L\zeta,
a_{\lambda,\gamma-2(s-1)}\partial_{r}^{s}\zeta\right).\nonumber
\end{eqnarray}
Note that since all the coefficients of $L$ vanish at $r=\sigma$,
except $\overline{B}$, no boundary terms at $r=\sigma$ appear in
the above formula.  Furthermore since
$r^{4}a_{\lambda,\gamma-2(s-1)}=a_{\lambda,\gamma-2(s+1)}$ and
$\partial_{r}^{s}\zeta\in H^{(2,1,\gamma-2(s-1))}(\Omega)$, $0\leq
s\leq m-2$ (for $m-1\leq s\leq m$ the boundary behavior of
$\partial_{r}^{s}\zeta$ given by (3.3) is also adequate), Lemma
3.2 implies that
\begin{eqnarray}
&
&\sum_{s=0}^{m}(L\partial_{r}^{s}\zeta,\lambda^{-s}r^{4}a_{\lambda,\gamma-2(s-1)}
\partial_{r}^{s}\zeta)\\
&\geq& C\int_{\Omega}\sum_{s=0}^{m}\lambda^{-s}[
r^{-\gamma+2(s+1)}(\varphi\sin\theta(\partial_{r}^{s}\zeta)_{r}+r^{-1}\cos\theta
(\partial_{r}^{s}\zeta)_{\theta})^{2}+\lambda
r^{-\gamma+2s}(\partial_{r}^{s}\zeta)^{2}]\nonumber\\
&\geq& C\int_{\Omega}\left[\sum_{s=0}^{m-1}\lambda^{-s}
r^{-\gamma+2s}
(\partial_{r}^{s}\zeta_{\theta})^{2}+\sum_{s=0}^{m}\lambda^{1-s}r^{-\gamma+2s}
(\partial_{r}^{s}\zeta)^{2}\right]\nonumber\\
&\geq& C\parallel\zeta\parallel_{(m,1,\gamma)}^{2}\nonumber
\end{eqnarray}
if $\delta=\delta(\lambda)$ and $\varepsilon=\varepsilon(\delta)$
are sufficiently small.  Next we calculate
\begin{eqnarray*}
[\partial_{r}^{s},L]\zeta&=&\sum_{l<s}\left(\begin{array}{c}
s \\
l
\end{array}\right)[\partial_{r}^{s-l}\overline{K}(\partial_{r}^{l}\zeta)_{rr}
+\partial_{r}^{s-l}\overline{A}(\partial_{r}^{l}\zeta)_{r\theta}\\
&
&+\partial_{r}^{s-l}\overline{B}(\partial_{r}^{l}\zeta)_{\theta\theta}+
\partial_{r}^{s-l}\overline{C}(\partial_{r}^{l}\zeta)_{r}+
\partial_{r}^{s-l}\overline{D}(\partial_{r}^{l}\zeta)_{\theta}],
\end{eqnarray*}
and observe that integrating by parts, again with the help of the
boundary conditions in (3.3), produces
\begin{equation}
|(\sum_{l<s}\left(\begin{array}{c}
s \\
l
\end{array}\right)\partial_{r}^{s-l}\overline{K}(\partial_{r}^{l}\zeta)_{rr},
r^{4}a_{\lambda,\gamma-2(s-1)}\partial_{r}^{s}\zeta)|\leq
C_{s}\int_{\Omega} \sum_{l\leq s}r^{-\gamma+2l}
(\partial_{r}^{l}\zeta)^{2},
\end{equation}
\begin{equation*}
|(\sum_{l<s}\left(\begin{array}{c}
s \\
l
\end{array}\right)\partial_{r}^{s-l}\overline{A}(\partial_{r}^{l}\zeta)_{r\theta},
r^{4}a_{\lambda,\gamma-2(s-1)}\partial_{r}^{s}\zeta)|\leq
C_{s}\int_{\Omega}[r^{-\gamma+2s}(\partial_{r}^{s}\zeta)^{2}+
\sum_{l<s}r^{-\gamma+2l} (\partial_{r}^{l}\zeta_{\theta})^{2}],
\end{equation*}
\begin{equation*}
|(\sum_{l<s}\left(\begin{array}{c}
s \\
l
\end{array}\right)\partial_{r}^{s-l}\overline{B}(\partial_{r}^{l}\zeta)_{\theta\theta},
r^{4}a_{\lambda,\gamma-2(s-1)}\partial_{r}^{s}\zeta)|\leq
C_{s}\int_{\Omega}[r^{-\gamma+2s}(\partial_{r}^{s}\zeta)^{2}+
\sum_{l<s}r^{-\gamma+2l} (\partial_{r}^{l}\zeta_{\theta})^{2}],
\end{equation*}
\begin{equation*}
|(\sum_{l<s}\left(\begin{array}{c}
s \\
l
\end{array}\right)\partial_{r}^{s-l}\overline{C}(\partial_{r}^{l}\zeta)_{r},
r^{4}a_{\lambda,\gamma-2(s-1)}\partial_{r}^{s}\zeta)|\leq
C_{s}\int_{\Omega} \sum_{l\leq s}r^{-\gamma+2l}
(\partial_{r}^{l}\zeta)^{2},
\end{equation*}
\begin{equation*}
|(\sum_{l<s}\left(\begin{array}{c}
s \\
l
\end{array}\right)\partial_{r}^{s-l}\overline{D}(\partial_{r}^{l}\zeta)_{\theta},
r^{4}a_{\lambda,\gamma-2(s-1)}\partial_{r}^{s}\zeta)|\leq
C_{s}\int_{\Omega}[r^{-\gamma+2s}(\partial_{r}^{s}\zeta)^{2}+
\sum_{l<s}r^{-\gamma+2l} (\partial_{r}^{l}\zeta_{\theta})^{2}].
\end{equation*}
Also observe that
\begin{eqnarray*}
& &\sum_{l<s}\left(\begin{array}{c}
s \\
l
\end{array}\right)\partial_{r}^{s-l}r^{4}\partial_{r}^{l}L\zeta\\
&=& \sum_{s-4\leq l<s}\left(\begin{array}{c}
s \\
l
\end{array}\right)\frac{4!}{(4-s+l)!}r^{4-s+l}[\sum_{t\leq l}\left(\begin{array}{c}
l \\
t
\end{array}\right)(\partial_{r}^{l-t}\overline{K}(\partial_{r}^{t}\zeta)_{rr}\\
&
+&\!\!\partial_{r}^{l-t}\overline{A}(\partial_{r}^{t}\zeta)_{r\theta}+
\partial_{r}^{l-t}\overline{B}(\partial_{r}^{t}\zeta)_{\theta\theta}
+\partial_{r}^{l-t}\overline{C}(\partial_{r}^{t}\zeta)_{r}
+\partial_{r}^{l-t}\overline{D}(\partial_{r}^{t}\zeta)_{\theta})],
\end{eqnarray*}
so similar calculations yield
\begin{eqnarray}
& &|(\sum_{l<s}\left(\begin{array}{c}
s \\
l
\end{array}\right)\partial_{r}^{s-l}r^{4}\partial_{r}^{l}L\zeta,
a_{\lambda,\gamma-2(s-1)}\partial_{r}^{s}\zeta)|\\
&\leq& C_{s}\int_{\Omega}\left[\sum_{l\leq
s}r^{-\gamma+2l}(\partial_{r}^{l}\zeta)^{2}+
\sum_{l<s}r^{-\gamma+2l}
(\partial_{r}^{l}\zeta_{\theta})^{2}\right].\nonumber
\end{eqnarray}
Then the combination of (3.5), (3.6), (3.7), and (3.8) produces
(3.4) for $\lambda=\lambda(m)$ sufficiently large and
$\delta=\delta(\lambda)$, $\varepsilon=\varepsilon(\delta)$
sufficiently small.\par
   We need one last estimate before proving existence.  Since
$r^{4}\eta\in H^{(m,0,\gamma+2)}(\Omega)$ whenever $\eta\in
H^{(m,1,\gamma)}(\Omega)$, we have
\begin{eqnarray}
\parallel r^{4}v\parallel_{(-m,-1,\gamma)}&:=&\sup_{\eta\in H^{(m,1,\gamma)}
(\Omega)}\frac{|(\eta,r^{4}v)|}{\parallel\eta\parallel_{(m,1,\gamma)}}\\
&=&\sup_{\eta\in H^{(m,1,\gamma)}
(\Omega)}\frac{|(r^{4}(\lambda\theta^{2}-1)\eta,\zeta)_{(m,0,\gamma+2)}|}{\parallel
\eta\parallel_{(m,1,\gamma)}}\nonumber\\
&\leq&C\sup_{\eta\in H^{(m,1,\gamma)}
(\Omega)}\frac{\parallel\eta\parallel_{(m,0,\gamma)}\parallel\zeta\parallel_{(m,0,\gamma)}}{\parallel
\eta\parallel_{(m,1,\gamma)}}\nonumber\\
&\leq&C\parallel\zeta\parallel_{(m,1,\gamma)}.\nonumber
\end{eqnarray}
Here $(\cdot,\cdot)_{(m,0,\gamma+2)}$ denotes the inner product on
$H^{(m,0,\gamma+2)}(\Omega)$, and the norm
$\parallel\cdot\parallel_{(-m,-1,\gamma)}$ comes from the dual
space $H^{(-m,-1,\gamma)}(\Omega)$ of $H^{(m,1,\gamma)}(\Omega)$,
which may be obtained as the completion of $L^{2}(\Omega)$ in this
negative norm.\par
   Now apply (3.4) to obtain
\begin{eqnarray*}
\parallel \zeta\parallel_{(m,1,\gamma)}\parallel L^{*}r^{4}v
\parallel_{(-m,-1,\gamma)}&\geq&
(\zeta,L^{*}r^{4}v)\\
&=&(L\zeta,r^{4}v)\\
&=&(L\zeta,r^{4}\sum_{s=0}^{m}\lambda^{-s}(-1)^{s}\partial_{r}^{s}
(a_{\lambda,\gamma-2(s-1)}\partial_{r}^{s}\zeta))\\
&\geq& C\parallel \zeta\parallel_{(m,1,\gamma)}^{2},
\end{eqnarray*}
which when combined with (3.9) yields
\begin{equation}
\parallel r^{4}v\parallel_{(-m,-1,\gamma)}\leq C\parallel
L^{*}r^{4}v\parallel_{(-m,-1,\gamma)}.
\end{equation}
Consider the linear functional $F:X\rightarrow\mathbb{R}$ where
$X=L^{*}(r^{4}\widehat{C}^{\infty}(\Omega))$, given by
\begin{equation*}
F(L^{*}r^{4}v)=(f,r^{4}v).
\end{equation*}
According to (3.10) we have that $F$ is bounded on the subspace
$X$ of $H^{(-m,-1,\gamma)}(\Omega)$ since
\begin{eqnarray*}
|F(L^{*}r^{4}v)|=|(f,r^{4}v)|&\leq& \parallel
f\parallel_{(m,1,\gamma)}\parallel
r^{4}v\parallel_{(-m,-1,\gamma)}\\
&\leq& C\parallel f\parallel_{(m,1,\gamma)}\parallel
L^{*}r^{4}v\parallel_{(-m,-1,\gamma)}.
\end{eqnarray*}
Note that we use $f\in H^{(m,1,\gamma)}(\Omega)$ here (ie. $f$ has
to have this particular boundary behavior at $r=0$). Thus we can
apply the Hahn-Banach theorem to get a bounded extension of $F$
(still denoted $F$) defined on all of
$H^{(-m,-1,\gamma)}(\Omega)$. It follows that there exists a
unique $u\in H^{(m,1,\gamma)}(\Omega)$ such that
\begin{equation*}
F(\eta)=(u,\eta)\text{ }\text{ }\text{ all }\text{ }\text{
}\eta\in H^{(-m,-1,\gamma)}(\Omega).
\end{equation*}
Now restrict $\eta$ back to $X$ to obtain
\begin{equation*}
(u,L^{*}r^{4}v)=(f,r^{4}v)\text{ }\text{ }\text{ all }\text{
}\text{ }v\in\widehat{C}^{\infty}(\Omega).
\end{equation*}
Since every $\overline{v}\in\widehat{C}^{\infty}(\Omega)$ can be
written as $\overline{v}=r^{4}v$ for some
$v\in\widehat{C}^{\infty}(\Omega)$, $u$ is a weak solution of
(3.1).  Q.E.D.\medskip

  If $f\in H^{(m,1,\gamma)}(\Omega)\cap
C^{m_{*}-4}(\Omega)$, then the strict ellipticity of $L$ on the
interior of $\Omega$ shows that the solution given by Theorem 3.1
satisfies $u\in H^{(m,1,\gamma)}(\Omega)\cap C^{m_{*}-3}(\Omega)$,
since the coefficients of $L$ will be in $C^{m_{*}-4}$ when $g\in
C^{m_{*}}$.  In particular $Lu=f$ pointwise on the interior of
$\Omega$, so that
\begin{equation}
-u_{\theta\theta}=\overline{B}^{-1}[\overline{K}u_{rr}+\overline{A}u_{r\theta}
+\overline{C}u_{r}+\overline{D}u_{\theta}-f]\text{ }\text{ }\text{
in }\text{ }\text{ }\Omega.
\end{equation}
Since $\overline{B}^{-1}=O(r^{2})$ as $r\rightarrow 0$ and $u\in
H^{(m,1,\gamma)}(\Omega)$ it follows that $u_{\theta\theta}\in
H^{(m-2,0,\gamma)}(\Omega)$.  Let
$H_{\gamma}^{m}(\Omega):=H^{(m,m,\gamma)}(\Omega)$ with norm
$\parallel\cdot\parallel_{(m,\gamma)}$.  Then if $f\in
H_{\gamma}^{m}(\Omega)$ we may continue to differentiate the
expression for $u_{\theta\theta}$ to eventually obtain (by
induction) that $u\in H_{\gamma}^{m}(\Omega)$.\par
  In order to determine the boundary values for $u$, integrate by
parts in expression (3.2) to obtain
\begin{eqnarray*}
0&=&\int_{\partial\Omega}(\overline{K}vu_{r}\nu_{1}-(\overline{K}
v)_{r}u\nu_{1}+\overline{A}vu_{\theta}\nu_{1}-(\overline{A}v)_{r}u\nu_{2}\\
&
&+\overline{B}vu_{\theta}\nu_{2}-(\overline{B}v)_{\theta}u\nu_{2}
+\overline{C}vu\nu_{1}+\overline{D}vu\nu_{2}),
\end{eqnarray*}
for all $v\in\widehat{C}^{\infty}(\Omega)$.  This implies that
$u(r,0)=u(r,\delta)=0$, as $\overline{B}>0$.  It also shows that
we cannot arbitrarily prescribe boundary values for $u$ at
$r=\sigma$, since all the coefficients of $L$ (except
$\overline{B}$) vanish at $r=\sigma$.  However it is clear that
the boundary values at $r=\sigma$ are given explicitly in terms of
$f(\sigma,\theta)$ according to (3.11); although we will not have
need of this fact.  Moreover, the boundary behavior at $r=0$ is
completely determined by the fact that $u\in
H^{m}_{\gamma}(\Omega)$, so that if $m\geq s_{0}+1$ and
$\gamma>2m$ then $u$ will vanish to the desired order $s_{0}$ at
$r=0$ by Lemma 3.1.  We summarize all that we have found with the
following theorem, and also give an a priori estimate needed for
the Nash-Moser iteration of section $\S 5$.\medskip

\textbf{Theorem 3.2.}  \textit{Suppose that $g\in C^{m_{*}}$ and
$f\in \overline{C}^{\infty}(\Omega)$. If $s_{0}+1\leq m\leq
m_{*}-4$, $\gamma> 2m$, $|w|_{C^{6}}<1$, and $\delta=\delta(m)$,
$\varepsilon=\varepsilon(m,\delta)$ are sufficiently small then
there exists a unique solution $u\in H^{m}_{\gamma}(\Omega)\cap
C^{m_{*}-3}(\Omega)$ of (3.1). Furthermore, there exists a
constant $C_{m}$ independent of $\delta$ and $\varepsilon$ such
that}
\begin{equation*}
\parallel u\parallel_{(m,\gamma)}\leq C_{m}(\parallel
f\parallel_{m+2+\gamma}+\parallel w\parallel_{m+6}\parallel
f\parallel_{5+\gamma}),
\end{equation*}
\textit{for each $m\leq m_{*}-6$.}\medskip

\textit{Proof.}  The first half of this theorem follows from the
discussion directly above, and thus it only remains to establish
the estimate.  From Lemma 3.2 we have that
\begin{equation*}
\int_{\Omega}r^{-\gamma+2}(\lambda u^{2}+u_{\theta}^{2})\leq C
\int_{\Omega}r^{-\gamma+6}(f^{2}+r^{-2}u_{r}^{2}).
\end{equation*}
Now differentiate equation (3.1) to obtain,
\begin{equation*}
Lu_{r}=f_{r}-\overline{K}_{r}u_{rr}-\overline{A}_{r}u_{r\theta}
-\overline{B}_{r}u_{\theta\theta}-\overline{C}_{r}u_{r}
-\overline{D}_{r}u_{\theta}.
\end{equation*}
Solving for $u_{\theta\theta}$ as in (3.11) then yields
\begin{eqnarray*}
L_{1}u_{r}&:=&\overline{K}(u_{r})_{rr}+\overline{A}(u_{r\theta})_{r}+
\overline{B}(u_{r})_{\theta\theta}\\
& &+(\overline{C}+\overline{K}_{r}-\overline{B}^{-1}
\overline{B}_{r}\overline{K})(u_{r})_{r}+(\overline{D}+\overline{A}_{r}
-\overline{B}^{-1} \overline{B}_{r}\overline{A})(u_{r})_{\theta}
\\
&=&f_{r}-(\overline{C}_{r}-\overline{B}^{-1}\overline{B}_{r}
\overline{C})u_{r}-(\overline{D}_{r}-\overline{B}^{-1}
\overline{B}_{r}\overline{D})u_{\theta}.
\end{eqnarray*}
Applying the proof of Lemma 3.2 to the above equation gives
\begin{equation*}
C\int_{\Omega}[\lambda r^{-\gamma+4}u_{r}^{2}+r^{-\gamma+6}
(\varphi\sin\theta u_{rr}+r^{-1}\cos\theta u_{r\theta})^{2}]\leq
(a_{\lambda,\gamma-6}u_{r},L_{1}u_{r}),
\end{equation*}
from which we find that
\begin{equation*}
\int_{\Omega}r^{-\gamma+4}(\lambda u_{r}^{2}+u_{r\theta}^{2})\leq
C\int_{\Omega}r^{-\gamma+6}(f^{2}+r^{2}f_{r}^{2}+u_{rr}^{2})
\end{equation*}
if $\lambda$ is sufficiently large.  Eventually, with the help of
$|w|_{C^{6}}<1$ and (3.11), we obtain
\begin{equation*}
\parallel u\parallel_{(3,3,\gamma-2)}\leq C\parallel
f\parallel_{(3,3,\gamma-6)}.
\end{equation*}
By repeatedly differentiating with respect to $r$, we can continue
this procedure and apply the Gagliardo-Nirenberg inequalities
(Lemma 5.2 below) in the usual way to obtain
\begin{eqnarray}
& &\int_{\Omega}\left(\lambda\sum_{s\leq
m}r^{-\gamma+2s+2}(\partial_{r}^{s}u)^{2}+\sum_{s<
m}r^{-\gamma+2s+2}(\partial_{r}^{s}u_{\theta})^{2}\right)\\
&\leq& C_{m}(\parallel
f\parallel^{2}_{(m,0,\gamma-6)}+\Lambda_{m}^{2}\parallel
u\parallel^{2}_{(3,3,\gamma-2)})\nonumber\\
&\leq&\!\!C_{m}(\parallel
f\parallel_{(m,\gamma)}^{2}+\Lambda_{m}^{2}\parallel
f\parallel_{(3,\gamma)}^{2}),\nonumber
\end{eqnarray}
where
\begin{equation*}
\Lambda_{m}=\parallel r^{2}\overline{K}\parallel_{m}+
\parallel r^{2}\overline{A}\parallel_{m}+\parallel r^{2}\overline{B}\parallel_{m}
+\parallel r^{2}\overline{C}\parallel_{m}+
\parallel r^{2}\overline{D}\parallel_{m}+1.
\end{equation*}\par
  To see how this works, we will show the calculation for one
representative term that appears on the right-hand side of the
equation after differentiating; the remaining terms may be treated
similarly.  Differentiate (3.1) $s$-times to obtain
\begin{eqnarray*}
L_{s}\partial_{r}^{s}u&:=&\overline{K}(\partial_{r}^{s}u)_{rr}+\overline{A}(\partial_{r}^{s}u)_{r\theta}
+\overline{B}(\partial_{r}^{s}u)_{\theta\theta}\\
& &+(\overline{C}+s\overline{K}_{r}-s\overline{B}^{-1}
\overline{B}_{r}\overline{K})(\partial_{r}^{s}u)_{r}+(\overline{D}+s\overline{A}_{r}
-s\overline{B}^{-1}
\overline{B}_{r}\overline{A})(\partial_{r}^{s}u)_{\theta}
\\
&=&\partial_{r}^{s}f-\partial_{r}^{s-1}(\overline{C}_{r}u_{r})+\cdots.
\end{eqnarray*}
Since $\partial_{r}^{s}u\in H^{(2,1,\gamma-2s)}(\Omega)$ when
$s\leq m_{*}-6$, the basic estimate yields
\begin{equation*}
\int_{\Omega}r^{-\gamma+2s+2}[\lambda(\partial_{r}^{s}u)^{2}
+(\partial_{r}^{s}u_{\theta})^{2}]\leq
C\int_{\Omega}r^{-\gamma+2s+6}[(\partial_{r}^{s}f)^{2}
+(\partial_{r}^{s-1}(\overline{C}_{r}u_{r}))^{2}]+\cdots.
\end{equation*}
Furthermore observe that
\begin{equation*}
\partial_{r}^{s-1}(\overline{C}_{r}u_{r})=\sum_{l\leq s-1}\left(\begin{array}{c}
s-1 \\
l
\end{array}\right)\partial_{r}^{s-1-l}\overline{C}_{r}\partial_{r}^{l}u_{r},
\end{equation*}
and
\begin{equation*}
\partial_{r}^{s-1-l}\overline{C}_{r}=\partial_{r}^{s-1-l}(r^{-2}r^{2}\overline{C}_{r})=\sum_{t_{1}\leq
s-1-l}C_{t_{1}}r^{-2-(s-1-l-t_{1})}\partial_{r}^{t_{1}}(r^{2}\overline{C}_{r}),
\end{equation*}
\begin{equation*}
r^{-\gamma/2+s+3+[-2-(s-1-l-t_{1})]}\partial_{r}^{l}u_{r}=\sum_{t_{2}\leq
l}C_{t_{2}}\partial_{r}^{t_{2}}(r^{-\gamma/2+t_{1}+t_{2}+2}u_{r}),
\end{equation*}
for some constants $C_{t_{1}}$ and $C_{t_{2}}$.  Therefore we may
apply Lemma 5.2 $(i)$, the Sobolev embedding theorem, and
$|w|_{C^{6}}<1$ to find
\begin{eqnarray*}
& &\parallel
r^{-\gamma/2+s+3}\partial_{r}^{s-1}(\overline{C}_{r}u_{r})\parallel\\
&\leq& \sum_{2\leq t_{1}+t_{2}\leq s-1}C_{t_{1}t_{2}}\parallel
\partial_{r}^{t_{1}}(r^{2}\overline{C}_{r})\partial_{r}^{t_{2}}(r^{-\gamma/2+
t_{1}+t_{2}+2}u_{r})\parallel\\
& &+\sum_{t_{1}+t_{2}<2}C_{t_{1}t_{2}}\parallel
\partial_{r}^{t_{1}}(r^{2}\overline{C}_{r})\partial_{r}^{t_{2}}(r^{-\gamma/2+
t_{1}+t_{2}+2}u_{r})\parallel\\
&\leq&C_{s}(|r^{2}\overline{C}_{r}|_{\infty}\parallel
u_{r}\parallel_{(s-1,0,\gamma-4)}+\parallel
r^{2}\overline{C}_{r}\parallel_{s-1}|r^{-\gamma/2+4}u_{r}|_{\infty})\\
& &+C_{s}\parallel u\parallel_{(2,2,\gamma-2)}\\
&\leq& C_{s}(\parallel
u\parallel_{(s,0,\gamma-2)}+\Lambda_{s}\parallel
u\parallel_{(3,3,\gamma-2)}),
\end{eqnarray*}
where $|\cdot|_{\infty}$ denotes the $L^{\infty}(\Omega)$ norm.
Note that the first term of the last line above, may be absorbed
into the left-hand side of (3.12) for large $\lambda(s)$.\par
  The remaining derivatives of $u$ involving higher orders of $\partial_{\theta}$
may be estimated by differentiating (3.11) and using the above
estimates. Since the coefficients of $L$ depend on the derivatives
of $w$ up to and including order 3 and the derivatives of $\xi$
(the coordinates of Lemma 2.1) up to and including order 2, with
the help of (2.4) we obtain
\begin{equation*}
\parallel u\parallel_{(m,\gamma-2)}\leq C_{m}(\parallel
f\parallel_{(m,\gamma)}+\parallel w\parallel_{m+6}\parallel
f\parallel_{(3,\gamma)}).
\end{equation*}
Lastly since $f$ vanishes to all orders at $r=0$, a little
calculation shows that
\begin{equation*}
\int_{\Omega}r^{-\gamma+2s}(\partial_{r}^{s}\partial_{\theta}^{t}f)^{2}\leq
C_{s}\int_{\Omega}(\partial_{r}^{s+\gamma/2+1}\partial_{\theta}^{t}f)^{2},
\end{equation*}
from which the desired result follows.  Q.E.D.

\begin{center}
\textbf{4. Linear Theory in the Hyperbolic Regions}
\end{center}\setcounter{equation}{0}
\setcounter{section}{4}

  Here we shall study the question of existence and regularity for
the Cauchy problem associated with the operator $L(w)$ of Lemma
2.1, in the hyperbolic regions. In the previous section concerning
the elliptic regions, after cutting-off some of the coefficients
away from the origin, we were able to invert $L(w)$ in the
appropriate function spaces.  However in the hyperbolic case, we
will not necessarily be able to make such an inversion, and as a
result we must consider a regularized version of $L(w)$ as we
explain below. For convenience let $(x,y)$ denote the coordinates
$(\xi^{1},\xi^{2})$ of Lemma 2.1, so that a portion of the given
hyperbolic region $\Omega^{-}_{\varrho}\cap B_{\sigma}(0)$,
$\varrho=1,2$, may be written as
\begin{equation*}
\Omega=\{(x,y)\mid h(x)<y<\sigma\}
\end{equation*}
for some Lipschitz function $h(x)$ satisfying $h(0)=0$ and
$|h(x)-|x||_{C^{1}}=O(\sigma)$.  Then in these coordinates $L(w)$
is given by
\begin{equation*}
Lu=(\overline{K}u_{x})_{x}+u_{yy}+\overline{C}u_{x}+\overline{D}u_{y}
\end{equation*}
with
\begin{equation*}
\overline{K}=\overline{k}K,\text{ }\text{ }\text{ }\text{ }\text{
}\overline{C}=\overline{c}K+(a^{22})^{-2}\partial_{x^{1}}\xi^{1}
\partial_{x^{1}}\Phi(w),\text{ }\text{ }\text{ }\text{ }\text{
}\overline{D}=\varepsilon^{2}\overline{d},
\end{equation*}
in the notation of Lemma 2.1.  Consider the Cauchy problem
\begin{equation}
Lu=f\text{ }\text{ in }\text{ }\Omega,\text{ }\text{
}u|_{\partial\Omega_{1}}=\phi,\text{ }\text{
}u_{y}|_{\partial\Omega_{1}}=\psi,\text{ }\text{
}\partial^{\alpha}u(0,0)=0,\text{ }\text{ }0\leq|\alpha|\leq
\alpha_{0},
\end{equation}
for some given data $\phi$, $\psi$ and a large integer
$\alpha_{0}$, where $\partial\Omega_{1}$ denotes the ``bottom"
portion of the boundary given by $y=h(x)$.  The ``top" portion of
the boundary will be denoted $\partial\Omega_{2}$, and is given by
$y=\sigma$. As in problem (3.1) the solution is required to vanish
to high order at the origin, in order to satisfy certain
compatibility conditions which arise when constructing solutions
to the nonlinear problem (2.1) (see section $\S 5$).  Of course a
necessary condition to have such behavior is that $\phi$, $\psi$,
and $f$ must also vanish to a corresponding high order at the
origin.\par
  Equation (4.1) is degenerate hyperbolic, and as such,
solvability of the Cauchy problem depends on the so called Levi
conditions.  These are relations between the coefficients
$\overline{K}$ and $\overline{C}$, which if satisfies would
guarantee existence for the Cauchy problem (when $y=h(x)$ is smooth
and noncharacteristic).  A simple example of such a relation is the
condition (see [14])
\begin{equation*}
\overline{C}\leq M\sqrt{|\overline{K}|}\text{ }\text{ }\text{ in
}\text{ }\text{ }\Omega,
\end{equation*}
for some constant $M>0$.  Unfortunately the quantity
$\partial_{x^{1}}\Phi(w)$ present in $\overline{C}$ prevents the
validity of this inequality.  However in the Nash iteration of the
next section, $\partial_{x^{1}}\Phi(w)$ will be uniformly small.
Therefore it is natural to consider the regularized Cauchy problem
\begin{equation}
L_{\theta}u=f\text{ }\text{ in }\text{ }\Omega,\text{ }\text{
}u|_{\partial\Omega_{1}}=\phi,\text{ }\text{
}u_{y}|_{\partial\Omega_{1}}=\psi,\text{ }\text{
}\partial^{\alpha}u(0,0)=0,\text{ }\text{
}0\leq|\alpha|\leq\alpha_{0},
\end{equation}
where $L_{\theta}$ differs from $L$ in that $\overline{K}$ is
replaced by $\overline{K}_{\theta}:=\overline{K}-\theta$ with
$\theta=|\Phi(w)|_{C^{1}}$.  There then exists a constant $M>0$
such that
\begin{equation}
\overline{C}\leq M|\overline{K}_{\theta}|\text{ }\text{ }\text{ in
}\text{ }\text{ }\Omega.
\end{equation}
However we cannot simply apply the results of [14] to obtain the
desired solution of (4.2), since the Cauchy surface $y=h(x)$ is not
smooth.  We will therefore prove existence by hand in what follows.
Note that with regards to the Nash iteration $\theta$ is of
quadratic error (see section $\S 5$), so that the small perturbation
(4.2) will not affect convergence of this procedure.\par
  It will be convenient to first establish an existence result for
(4.2) with homogeneous Cauchy data and with $f$ vanishing to high
order on all of $\partial\Omega_{1}$.  To this end we define
$H^{(m,l)}(\Omega)$ ($H_{0}^{(m,l)}(\Omega)$) to be the closure of
all $C^{\infty}(\Omega)$ functions (which vanish to all orders at
$\partial\Omega_{1}$) in the norm
\begin{equation*}
\parallel u\parallel^{2}_{(m,l)}=\int_{\Omega}\sum_{0\leq s\leq
m\atop 0\leq t\leq
l}\lambda^{-s}(\partial_{x}^{s}\partial_{y}^{t}u)^{2},
\end{equation*}
where $\lambda>0$ is a large parameter.  The $L^{2}(\Omega)$ inner
product will as usual be denoted by $(\cdot,\cdot)$, and the
formal adjoint of $L_{\theta}$ by $L_{\theta}^{*}$.  Also as in
Lemma 2.1 $w\in C^{\infty}$ throughout this section.\medskip

\textbf{Theorem 4.1.}  \textit{Suppose that $g\in C^{m_{*}}$,
$|w|_{C^{4}}<1$, and $f\in H_{0}^{(m,0)}(\Omega)$. If $m\leq
m_{*}-6$ and $\varepsilon$ is sufficiently small, then there
exists a weak solution $u_{\theta}\in H_{0}^{(m,1)}(\Omega)$ of
(4.2) with $\phi,\psi=0$, for each $\theta>0$.  That is}
\begin{equation}
(u_{\theta},L_{\theta}^{*}v)=(f,v)\text{ }\text{ }\textit{ all
}\text{ }\text{ }v\in C^{\infty}(\Omega)
\end{equation}
\textit{with
$v|_{\partial\Omega_{2}}=v_{y}|_{\partial\Omega_{2}}=0$.}\medskip

\textit{Proof.}  Set
$b(x,y)=\overline{K}_{\theta}^{-1}(x,y)e^{-\lambda y}$ and let
$\zeta$ be the unique solution of
\begin{equation}
\sum_{s=0}^{m}(-1)^{s+1}\lambda^{-s}\partial_{x}^{s}(b\partial_{x}^{s}\zeta_{y})=v\text{
}\text{ }\text{ in }\text{ }\text{ }\Omega,\text{ }\text{
}\zeta|_{\partial\Omega_{1}}=\partial_{x}^{s}\zeta_{y}|_{\partial\Omega_{1}}=0,\text{
}\text{ }0\leq s\leq m-1,
\end{equation}
where $v$ is as stated in the theorem.  Note that for each
$y\in(0,\sigma)$, this equation may be interpreted as an ODE in
$\zeta_{y}$, and therefore the theory of such equations guarantees
the existence of a unique solution to (4.5).  Furthermore if the
metric $g\in C^{m_{*}}$ as in Lemma 2.1, then the coefficients of
(4.5) are in $C^{m_{*}-m-4}$.  Thus as long as $m\leq m_{*}-6$, we
have $\partial_{x}^{s}\zeta\in C^{2}(\Omega)$ for $0\leq s\leq
2m-1$.\par
  We first note that the solution $\zeta$ of (4.5) satisfies extra
boundary conditions, namely
\begin{equation}
\partial_{x}^{s}\partial_{y}^{t}\zeta|_{\partial\Omega_{1}}=0,\text{
}\text{ }\text{ }\text{ }s+t\leq m,\text{ }\text{ }\text{ }\text{
}0\leq t\leq 2.
\end{equation}
To see this let $\partial_{T}\zeta$ be differentiation along the
right portion of $\partial\Omega_{1}$ (that is, the curve
$y=h(x)$, $x>0$), which we denote by $\partial\Omega_{1}^{+}$.
Then since
$\zeta|_{\partial\Omega_{1}^{+}}=\zeta_{y}|_{\partial\Omega_{1}^{+}}=0$
we have
\begin{equation*}
0=\partial_{T}\zeta=\frac{1}{\sqrt{1+h'^{\!\text{ }2}}}
(\zeta_{x}+h'\zeta_{y})|_{\partial\Omega_{1}^{+}}
=\frac{\zeta_{x}}{\sqrt{1+h'^{\!\text{
}2}}}|_{\partial\Omega_{1}^{+}},
\end{equation*}
where $h'=dh/dx$.  It follows that
\begin{equation*}
0=\partial_{T}\zeta_{x}=\frac{1}{\sqrt{1+h'^{\!\text{
}2}}}(\zeta_{xx}+h'\zeta_{xy})|_{\partial\Omega_{1}^{+}}
=\frac{\zeta_{xx}}{\sqrt{1+h'^{\!\text{
}2}}}|_{\partial\Omega_{1}^{+}},
\end{equation*}
\begin{equation*}
0=\partial_{T}\zeta_{y}=\frac{1}{\sqrt{1+h'^{\!\text{
}2}}}(\zeta_{xy}+h'\zeta_{yy})|_{\partial\Omega_{1}^{+}}
=\frac{h'\zeta_{yy}}{\sqrt{1+h'^{\!\text{
}2}}}|_{\partial\Omega_{1}^{+}}.
\end{equation*}
Since $|h(x)-|x||_{C^{1}}=O(\sigma)$ this shows (4.6) for $m=2$,
and continuing this procedure establishes the full result on
$\partial\Omega_{1}^{+}$.  Also the same arguments hold to yield
(4.6) on $\partial\Omega_{1}^{-}$, the left portion of
$\partial\Omega_{1}$ (that is, the curve $y=h(x)$, $x<0$).\par
  We will now establish the basic estimate on which the existence
proof is based.  More precisely we will show that
\begin{equation}
(L_{\theta}\zeta,\sum_{s=0}^{m}(-1)^{s+1}\lambda^{-s}\partial_{x}^{s}(b\partial_{x}^{s}\zeta_{y}))
\geq C\parallel \zeta\parallel_{(m,1)}^{2},
\end{equation}
for some constant $C>0$.  Observe that the following calculations
hold for $0\leq s\leq m$ according to the boundary conditions
(4.6):
\begin{eqnarray*}
((\overline{K}_{\theta}\zeta_{x})_{x},(-1)^{s+1}\partial_{x}^{s}(b\partial_{x}^{s}\zeta_{y}))&=&
\int_{\Omega}\partial_{x}(b\partial_{x}^{s}\zeta_{y})\partial_{x}^{s}(\overline{K}_{\theta}
\zeta_{x})\\
&=&-\int_{\Omega}\left[\frac{1}{2}(b\overline{K}_{\theta})_{y}
(\partial_{x}^{s+1}\zeta)^{2}-b_{x}\overline{K}_{\theta}\partial_{x}^{s}\zeta_{y}\partial_{x}^{s+1}\zeta\right]\\
& &-\int_{\Omega}b\partial_{x}^{s}\zeta_{y}\partial_{x}\left[
\sum_{l=1}^{s}\left(\begin{array}{c}
s \\
l
\end{array}\right)\partial_{x}^{l}\overline{K}_{\theta}
\partial_{x}^{s+1-l}\zeta\right]\\
&
&+\int_{\partial\Omega}\frac{1}{2}b\overline{K}_{\theta}(\partial_{x}^{s+1}\zeta)^{2}
\nu_{2},
\end{eqnarray*}
\begin{eqnarray*}
(\zeta_{yy},(-1)^{s+1}\partial_{x}^{s}(b\partial_{x}^{s}\zeta_{y}))&=&
-\int_{\Omega}b\partial_{x}^{s}\zeta_{yy}\partial_{x}^{s}\zeta_{y}
+\int_{\partial\Omega}b\partial_{x}^{s-1}\zeta_{yy}\partial_{x}^{s}\zeta_{y}\nu_{1}\\
&=&\int_{\Omega}\frac{1}{2}b_{y}(\partial_{x}^{s}\zeta_{y})^{2}
+\int_{\partial\Omega}
\left[b\partial_{x}^{s-1}\zeta_{yy}\partial_{x}^{s}\zeta_{y}\nu_{1}
-\frac{1}{2}b(\partial_{x}^{s}\zeta_{y})^{2}\nu_{2}\right],
\end{eqnarray*}
\begin{eqnarray*}
(\overline{C}\zeta_{x},(-1)^{s+1}\partial_{x}^{s}(b\partial_{x}^{s}\zeta_{y}))&=&
-\int_{\Omega}b\partial_{x}^{s}\zeta_{y}\partial_{x}^{s}(\overline{C}\zeta_{x})\\
&=&-\int_{\Omega}\left[b\overline{C}\partial_{x}^{s}\zeta_{y}\partial_{x}^{s+1}\zeta
+b\partial_{x}^{s}\zeta_{y}\sum_{l=1}^{s}\left(\begin{array}{c}
s \\
l
\end{array}\right)\partial_{x}^{l}\overline{C}
\partial_{x}^{s+1-l}\zeta\right],
\end{eqnarray*}
\begin{eqnarray*}
(\overline{D}\zeta_{y},(-1)^{s+1}\partial_{x}^{s}(b\partial_{x}^{s}\zeta_{y}))&=&
-\int_{\Omega}b\partial_{x}^{s}\zeta_{y}\partial_{x}^{s}(\overline{D}\zeta_{y})\\
&=&-\int_{\Omega}\left[b\overline{D}(\partial_{x}^{s}\zeta_{y})^{2}
+b\partial_{x}^{s}\zeta_{y}\sum_{l=1}^{s}\left(\begin{array}{c}
s \\
l
\end{array}\right)\partial_{x}^{l}\overline{D}
\partial_{x}^{s-l}\zeta_{y}\right],
\end{eqnarray*}
where $(\nu_{1},\nu_{2})$ denotes the unit outer normal to
$\partial\Omega$.  Next observe that if $\lambda>0$ is
sufficiently large then
\begin{equation*}
b_{y}=-\frac{1}{\overline{K}_{\theta}}\left(\lambda+\frac{\partial_{y}\overline{K}_{\theta}}
{\overline{K}_{\theta}}\right)e^{-\lambda y}>0,\text{ }\text{
}\text{ }-(b\overline{K}_{\theta})_{y}=\lambda e^{-\lambda y}>0,
\end{equation*}
and in light of the Levi condition (4.3) as well as
$\overline{D}=O(\varepsilon)$ we have
\begin{eqnarray*}
&
&-(b_{y}-2b\overline{D})(b\overline{K}_{\theta})_{y}-(b_{x}\overline{K}_{\theta}-b\overline{C})^{2}\\
&\geq& -\frac{e^{-2\lambda y}}{\overline{K}_{\theta}^{2}}\left[
\left(\lambda+\frac{\partial_{y}\overline{K}_{\theta}}
{\overline{K}_{\theta}}+O(\varepsilon)\right)\lambda\overline{K}_{\theta}
+O(\overline{K}_{x}^{2}+\overline{K}_{\theta}^{2})\right]\\
&\geq&-\frac{e^{-2\lambda
y}}{\overline{K}_{\theta}^{2}}\left[\frac{\lambda^{2}}{2}\overline{K}_{\theta}
+\lambda\overline{K}_{y}+O(\overline{K}_{x}^{2})\right]\\
&>& e^{-2\lambda y},
\end{eqnarray*}
as $\lambda\overline{K}_{y}+O(\overline{K}_{x}^{2})\leq 0$ near
$\partial\Omega_{1}$ for large $\lambda$.  It follows that for
large $\lambda$ depending on $m$,
\begin{eqnarray*}
& &(L_{\theta}\zeta,\sum_{s=0}^{m}(-1)^{s+1}\lambda^{-s}\partial_{x}^{s}(b\partial_{x}^{s}\zeta_{y}))\\
&\geq& C\parallel
\zeta\parallel_{(m,1)}^{2}+\frac{\lambda^{-m}}{2}\int_{\partial\Omega}(b\overline{K}_{\theta}
(\partial_{x}^{m+1}\zeta)^{2}\nu_{2}+2b\partial_{x}^{m-1}\zeta_{yy}\partial_{x}^{m}\zeta_{y}\nu_{1}
-b(\partial_{x}^{m}\zeta_{y})^{2}\nu_{2}),
\end{eqnarray*}
where we have used a Poincar\'{e} type inequality to estimate
$\parallel\zeta\parallel_{L^{2}(\Omega)}$.  Note that the boundary
integral has the correct sign on $\partial\Omega_{2}$.  We claim
that it also has the correct sign on $\partial\Omega_{1}$.  To see
this we use the boundary condition (4.6).  That is, by (4.6)
$\partial_{x}^{m-1}\zeta_{y}|_{\partial\Omega_{1}}=0$ so if
$\partial_{T}$ is differentiation along
$\partial\Omega_{1}-\{(0,0)\}$ then
\begin{equation*}
0=\partial_{T}\partial_{x}^{m-1}\zeta_{y}=-\nu_{2}\partial_{x}^{m}\zeta_{y}+\nu_{1}
\partial_{x}^{m-1}\zeta_{yy},
\end{equation*}
which yields
\begin{equation*}
b\partial_{x}^{m-1}\zeta_{yy}\partial_{x}^{m}\zeta_{y}\nu_{1}=b(\partial_{x}^{m}\zeta_{y})^{2}
\nu_{2}\geq 0.
\end{equation*}
Moreover since $\partial_{x}^{m}\zeta|_{\partial\Omega_{1}}=0$ we
have
\begin{equation*}
0=\partial_{T}\partial_{x}^{m}\zeta|_{\partial\Omega_{1}}=
-\nu_{2}\partial_{x}^{m+1}\zeta+\nu_{1}\partial_{x}^{m}\zeta_{y},
\end{equation*}
which yields
\begin{equation*}
b\overline{K}_{\theta}(\partial_{x}^{m+1}\zeta)^{2}\nu_{2}
=b\overline{K}_{\theta}(\partial_{x}^{m}\zeta_{y})^{2}\frac{\nu_{1}^{2}}{\nu_{2}}.
\end{equation*}
It follows that the boundary integral on $\partial\Omega_{1}$ is
nonnegative, and hence (4.7) holds.\par
  We will need one last estimate before proving existence.  Namely
\begin{eqnarray}
\parallel v\parallel_{(-m,0)}&:=&\sup_{\eta\in
H_{0}^{(m,0)}(\Omega)}\frac{|(\eta,v)|}{\parallel
\eta\parallel_{(m,0)}}\\
&=&\sup_{\eta\in
H_{0}^{(m,0)}(\Omega)}\frac{|(\eta,
\sum_{s=0}^{m}(-1)^{s+1}\lambda^{-s}\partial_{x}^{s}(b\partial_{x}^{s}\zeta_{y}))|}{\parallel
\eta\parallel_{(m,0)}}\nonumber\\
&\leq&\theta^{-1}C\parallel \zeta\parallel_{(m,1)},\nonumber
\end{eqnarray}
which follows after integration by parts.  Here
$\parallel\cdot\parallel_{(-m,0)}$ is the norm on the dual space
$H_{0}^{(-m,0)}(\Omega)$ of $H_{0}^{(m,0)}(\Omega)$.  This dual
space may be obtained as the completion of $L^{2}(\Omega)$ in the
norm $\parallel\cdot\parallel_{(-m,0)}$.\par
  To prove existence apply (4.7) to obtain
\begin{eqnarray*}
\parallel \zeta\parallel_{(m,1)}\parallel L_{\theta}^{*}v
\parallel_{(-m,-1)}&\geq&
(\zeta,L_{\theta}^{*}v)=(L_{\theta}\zeta,v)\\
&=&(L_{\theta}\zeta,\sum_{s=0}^{m}(-1)^{s+1}\lambda^{-s}\partial_{x}^{s}(b\partial_{x}^{s}\zeta_{y}))\geq
C\parallel \zeta\parallel_{(m,1)}^{2},
\end{eqnarray*}
which together with (4.8) implies that
\begin{equation}
\parallel
v\parallel_{(-m,0)}\leq\theta^{-1}C
\parallel L_{\theta}^{*}v\parallel_{(-m,-1)}\text{
}\text{ }\text{ all }\text{ }\text{
}v\in\widehat{C}^{\infty}(\Omega),
\end{equation}
where $\widehat{C}^{\infty}(\Omega)$ consists of all
$C^{\infty}(\Omega)$ functions with
$v|_{\partial\Omega_{2}}=v_{y}|_{\partial\Omega_{2}}=0$.  Consider
the linear functional $F:X\rightarrow\mathbb{R}$, where
$X=L_{\theta}^{*}\widehat{C}^{\infty}(\Omega)$, given by
\begin{equation*}
F(L_{\theta}^{*}v)=(f,v).
\end{equation*}
According to (4.9) we have that $F$ is bounded on the subspace $X$
of $H_{0}^{(-m,-1)}(\Omega)$ since
\begin{equation*}
|F(L_{\theta}^{*}v)|\leq\parallel f\parallel_{(m,0)}\parallel
v\parallel_{(-m,0)}\leq\theta^{-1}C\parallel f\parallel_{(m,0)}
\parallel L_{\theta}^{*}v\parallel_{(-m,-1)}.
\end{equation*}
Note that the generalized Schwarz inequality (the first inequality
in the above sequence) holds because $f\in H_{0}^{(m,0)}(\Omega)$,
that is, $f$ vanishes appropriately on $\partial\Omega_{1}$.  Thus
we can apply the Hahn-Banach theorem to obtain a bounded extension
of $F$ defined on all of $H_{0}^{(-m,-1)}(\Omega)$.  It follows
that there exists $u_{\theta}\in H^{(m,1)}_{0}(\Omega)$ such that
\begin{equation*}
F(\eta)=(u_{\theta},\eta)\text{ }\text{ }\text{ all }\text{
}\text{ }\eta\in H_{0}^{(-m,-1)}(\Omega).
\end{equation*}
Now restrict $\eta$ back to $X$ to obtain
\begin{equation*}
(u_{\theta},L_{\theta}^{*}v)=(f,v)\text{ }\text{ }\text{ all
}\text{ }\text{ }v\in\widehat{C}^{\infty}(\Omega).
\end{equation*}
Q.E.D.\medskip

  In order to obtain higher regularity for the solution given by
Theorem 4.1, we will utilize the following standard lemma
concerning the difference quotient
\begin{equation*}
u^{q}(x,y):=\frac{u(x,y+q)-u(x,y)}{q}.
\end{equation*}

\textbf{Lemma 4.1.} $i)$ \textit{ Let $u\in H^{(0,1)}(\Omega)$ and
$\Omega'\subset\subset\Omega$ (that is, $\Omega'$ is compactly
contained in $\Omega$). Then}
\begin{equation*}
\parallel u^{q}\parallel_{L^{2}(\Omega')}\leq\parallel
u_{y}\parallel_{L^{2}(\Omega)}
\end{equation*}
\textit{for all
$0<|q|<\frac{1}{2}\mathrm{dist}(\Omega',\partial\Omega)$.}\smallskip

$ii)$  \textit{ If $u\in L^{2}(\Omega)$ and $\parallel
u^{q}\parallel_{L^{2}(\Omega')}\leq C$ for all
$0<|q|<\frac{1}{2}\mathrm{dist}(\Omega',\partial\Omega)$, then
$u\in H^{(0,1)}(\Omega')$}.\medskip

  The Sobolev space of square integrable derivatives up to and
including order $m$ will be denoted by $H^{m}(\Omega)$ with norm
$\parallel\cdot\parallel_{m}$, and the completion of
$C^{\infty}(\Omega)$ functions which vanish to all orders at
$\partial\Omega_{1}$ in the norm $\parallel\cdot\parallel_{m}$
shall be denoted by $H_{0}^{m}(\Omega)$.\medskip

\textbf{Corollary 4.1.}  \textit{Under the hypotheses of Theorem
4.1, if $f\in H_{0}^{m}(\Omega)$ there exists a unique solution
$u_{\theta}\in H_{0}^{m}(\Omega)$ of (4.2) with $\phi,\psi=0$, for
each $\theta>0$.}\medskip

\textit{Proof.}  Let $u_{\theta}\in H_{0}^{(m,1)}(\Omega)$ be the
weak solution given by Theorem 4.1, so that (4.4) holds.  If
$m\leq 1$ then this corollary follows directly from Theorem 4.1,
so assume that $m\geq 2$. We may integrate by parts to obtain
\begin{eqnarray*}
-(\partial_{y}u_{\theta}+\overline{D}u_{\theta},v_{y})&=&
(f-\partial_{x}(\overline{K}_{\theta}\partial_{x}u_{\theta})-\overline{C}\partial_{x}
u_{\theta}+\overline{D}_{y}u_{\theta},v)\\
&
&+\int_{\partial\Omega}(\overline{C}vu_{\theta}\nu_{1}-v_{y}u_{\theta}
\nu_{2}-\overline{K}_{\theta}v_{x}u_{\theta}\nu_{1}+\overline{K}_{\theta}v\partial_{x}
u_{\theta}\nu_{1})
\end{eqnarray*}
for all $v\in\widehat{C}^{\infty}(\Omega)$, where
$(\nu_{1},\nu_{2})$ is the unit outer normal to $\partial\Omega$.
Note that since $u_{\theta},\partial_{x}u_{\theta}\in
H^{1}(\Omega)$ both $u_{\theta}|_{\partial\Omega}$ and
$\partial_{x}u_{\theta}|_{\partial\Omega}$ are meaningful in
$L^{2}(\partial\Omega)$, and in particular as
$u_{\theta},\partial_{x}u_{\theta}\in H_{0}^{1}(\Omega)$ we have
$u_{\theta}|_{\partial\Omega_{1}}=\partial_{x}u_{\theta}|_{\partial\Omega_{1}}=0$
in the $L^{2}(\partial\Omega_{1})$ sense.  Thus we may write
\begin{equation*}
(\overline{u}_{\theta},v_{y})=(\overline{f},v)\text{ }\text{
}\text{ all }\text{ }\text{ }v\in\widehat{C}^{\infty}(\Omega),
\end{equation*}
where
\begin{equation*}
\overline{u}_{\theta}=-\partial_{y}u_{\theta}-\overline{D}u_{\theta},\text{
}\text{ }\text{ }\text{
}\overline{f}=f-\partial_{x}(\overline{K}_{\theta}\partial_{x}u_{\theta})
-\overline{C}\partial_{x}u_{\theta}+\overline{D}_{y}u_{\theta}.
\end{equation*}
Furthermore
\begin{equation*}
(\overline{u}_{\theta}^{q},v_{y})=(\overline{f}^{q},v)\text{
}\text{ }\text{ all }\text{ }\text{ }v\in C_{c}^{\infty}(\Omega),
\end{equation*}
so that choosing a sequence $v_{i}\in C^{\infty}_{c}(\Omega)$ with
$v_{i}\rightarrow-\eta u_{\theta}^{q}$ in $H^{(0,1)}(\Omega)$ for
some nonnegative $\eta\in C^{\infty}_{c}(\Omega)$, implies that
\begin{eqnarray*}
\parallel\sqrt{\eta}\overline{u}_{\theta}^{q}\parallel^{2}&\leq&
|(\overline{f}^{q},\eta
u_{\theta}^{q})|+|(\overline{u}_{\theta}^{q},\eta_{y}u_{\theta}^{q})|
+|(\overline{u}_{\theta}^{q},\eta(\overline{D}u_{\theta})^{q})|\\
&\leq&
\parallel\sqrt{\eta}\overline{f}^{q}\parallel\parallel\sqrt{\eta}u_{\theta}^{q}\parallel
+\parallel\sqrt{\eta}\overline{u}_{\theta}^{q}\parallel\parallel\frac{\eta_{y}}{\sqrt{\eta}}u_{\theta}^{q}\parallel
+\parallel\sqrt{\eta}\overline{u}_{\theta}^{q}\parallel\parallel\sqrt{\eta}(\overline{D}u_{\theta})^{q}\parallel.
\end{eqnarray*}
Then since $u_{\theta},\overline{f}\in H^{(0,1)}(\Omega)$ and
$|\nabla\eta|^{2}\leq C\eta$, Lemma 4.1 $(i)$ yields
$\parallel\sqrt{\eta}\overline{u}_{\theta}^{q}\parallel\leq C$ for
some constant $C$ independent of $q$, if $|q|$ is sufficiently
small.  Now Lemma 4.1 $(ii)$ shows that $\overline{u}_{\theta}\in
H^{(0,1)}_{loc}(\Omega)$, as $\eta$ was arbitrary.  Hence
$\partial_{y}^{2}u_{\theta}\in L^{2}_{loc}(\Omega)$.  It follows
that the equation $L_{\theta}u_{\theta}=f$ holds in
$L^{2}_{loc}(\Omega)$, and since we can solve for
$\partial_{y}^{2}u_{\theta}$, we can boot-strap in the usual way
to obtain $u_{\theta}\in H^{m}_{loc}(\Omega)$.\par
  Next observe that the above restrictions on $\eta$ may be
relaxed if $q>0$, that is in this case $\eta$ is only required to
vanish in a neighborhood of $\partial\Omega_{2}$.  Then the same
procedure yields $\eta u_{\theta}\in H^{m}(\Omega)$ for all such
$\eta$.  Furthermore since $L_{\theta}$ is strictly hyperbolic, we
can use the regularity theory for such operators to obtain
estimates for $u_{\theta}$ near $\partial\Omega_{2}$.  It follows
that $u_{\theta}\in H^{m}(\Omega)$.\par
  To show that $u_{\theta}\in H_{0}^{m}(\Omega)$, it is enough to
observe that
\begin{equation}
\partial_{x}^{s}\partial_{y}^{t}u_{\theta}|_{\partial\Omega_{1}}=0,\text{
}\text{ }\text{ }s+t\leq m-1,
\end{equation}
where the equality is interpreted in the
$L^{2}(\partial\Omega_{1})$ sense when $s+t=m-1$.  This follows
from the fact that $f$ satisfies (4.10), in the following way.
First note that as in the arguments used to establish (4.6),
$u_{\theta}|_{\partial\Omega_{1}}=\partial_{x}u_{\theta}|_{\partial\Omega_{1}}=0$
implies that $\partial_{y}u_{\theta}|_{\partial\Omega_{1}}=0$.
Furthermore using the notation of those arguments we have
\begin{equation*}
0=\partial_{T}(\partial_{x}u_{\theta})=\frac{1}{\sqrt{1+h'^{\!\text{
}2}}}(\partial_{x}^{2}u_{\theta}
+h'\partial_{x}\partial_{y}u_{\theta})|_{\partial\Omega_{1}^{+}},
\end{equation*}
\begin{equation*}
0=\partial_{T}(\partial_{y}u_{\theta})=\frac{1}{\sqrt{1+h'^{\!\text{
}2}}}(\partial_{x}\partial_{y}u_{\theta}
+h'\partial_{y}^{2}u_{\theta})|_{\partial\Omega_{1}^{+}}.
\end{equation*}
However from equation (4.2) we find that
\begin{equation*}
(-\theta\partial_{x}^{2}u_{\theta}+\partial_{y}^{2}u_{\theta})|_{\partial\Omega_{1}^{+}}=0,
\end{equation*}
hence
\begin{equation*}
\partial_{x}^{2}u_{\theta}|_{\partial\Omega_{1}^{+}}=\partial_{y}^{2}u_{\theta}|_{\partial\Omega_{1}^{+}}=
\partial_{x}\partial_{y}u_{\theta}|_{\partial\Omega_{1}^{+}}=0.
\end{equation*}
The same arguments also apply to $\partial\Omega_{1}^{-}$, so by
differentiating equation (4.2) we can continue this procedure to
obtain (4.10).\par
  Lastly we note that since
$u_{\theta}|_{\partial\Omega_{1}}=|\nabla
u_{\theta}|_{\partial\Omega_{1}}=0$, (4.7) with $m=1$ yields
\begin{equation*}
(L_{\theta}u_{\theta},\sum_{s=0}^{1}(-1)^{s+1}\lambda^{-s}\partial_{x}^{s}(b\partial_{x}^{s}\partial_{y}u_{\theta}))
\geq C\parallel u_{\theta}\parallel_{(1,1)}^{2},
\end{equation*}
from which uniqueness follows.  Q.E.D.\medskip

  This corollary yields the existence of a regular solution to
(4.2) when $\phi=\psi=0$ and $f$ vanishes to high order on
$\partial\Omega_{1}$.  However, we are interested in solving (4.2)
in the general case when $\phi$, $\psi$, and $f$ vanish to high
order at the origin but are otherwise arbitrary.  The next lemma
will enable us to obtain the general case from Corollary 4.1. Here
and below $\overline{H}_{0}^{m}(\Omega)$ will denote the
completion of $\overline{C}^{\infty}(\Omega)$ in the norm
$\parallel\cdot\parallel_{m}$, and
$\overline{H}_{0}^{m}(\partial\Omega_{1})$ will be defined
similarly with respect to the
$\parallel\cdot\parallel_{m,\partial\Omega_{1}}$ norm. Recall that
$\overline{C}^{\infty}(\Omega)$ consists of $C^{\infty}(\Omega)$
functions which vanish in a neighborhood of the origin.\medskip

\textbf{Lemma 4.2.}  \textit{Suppose that $g\in C^{m_{*}}$,
$|w|_{C^{6}}<1$, and
$\phi\in\overline{H}_{0}^{m+1}(\partial\Omega_{1})$,
$\psi\in\overline{H}_{0}^{m}(\partial\Omega_{1})$,
$f\in\overline{H}_{0}^{m}(\Omega)$ with $m\leq m_{*}-6$.  Then
there exists a function
$\eta_{\theta}\in\overline{H}_{0}^{m+2}(\Omega)$ such that
$\eta_{\theta}|_{\partial\Omega_{1}}=\phi$,
$\partial_{y}\eta_{\theta}|_{\partial\Omega_{1}}=\psi$, and
$\partial_{y}^{t}(f-L_{\theta}\eta_{\theta})|_{\partial\Omega_{1}}=0$,
$0\leq t\leq m-1$, with}
\begin{eqnarray*}
\parallel \eta_{\theta}\parallel_{m+2}&\leq& C_{m}(\parallel
f\parallel_{m}+\parallel\phi\parallel_{m+1,\partial\Omega_{1}}
+\parallel\psi\parallel_{m,\partial\Omega_{1}}\\
& &+\parallel w\parallel_{m+6}(\parallel
f\parallel_{2}+\parallel\phi\parallel_{2,\partial\Omega_{1}}
+\parallel\psi\parallel_{2,\partial\Omega_{1}})).
\end{eqnarray*}

\textit{Proof.}  We may assume that a unique solution
$u_{\theta}\in\overline{H}_{0}^{m+2}(\Omega)$ of (4.2) exists, since
here we shall only use its boundary values which can be explicitly
determined in terms of $\phi$, $\psi$, and $f$, as
$\partial\Omega_{1}$ is noncharacteristic for $L_{\theta}$. Then
because $\Omega$ is a Lipschitz domain, the linear restriction map
$H^{m+2}(\Omega)\rightarrow H^{m+1}(\partial\Omega_{1})$ is bounded
and onto (see [18]).  By quotienting with the kernel and applying
the closed graph theorem, we obtain a bounded inverse
$H^{m+1}(\partial\Omega_{1})\rightarrow
H^{m+2}(\Omega)/H^{m+2}_{0}(\Omega)$ with respect to the quotient
norm.  We may then use this map to obtain an extension
$\eta_{\theta}$ of $u_{\theta}$ from $\partial\Omega_{1}$ to
$\Omega$ with
$\partial^{\alpha}\eta_{\theta}|_{\partial\Omega_{1}}=\partial^{\alpha}u_{\theta}|_{\partial\Omega_{1}}$
for all $|\alpha|\leq m+1$, and
\begin{equation*}
\parallel \eta_{\theta}\parallel_{m+2}\leq
C_{m}\sum_{|\alpha|\leq
m+1}\parallel\partial^{\alpha}u_{\theta}\parallel_{0,\partial\Omega_{1}}.
\end{equation*}
Applying the Gagliardo-Nirenberg inequalities (Lemma 5.2) to the
expression for $\partial^{\alpha}u_{\theta}|_{\partial\Omega_{1}}$
in terms of $\phi$, $\psi$, and $f$ yields the desired
estimate.\par
  We remark that an equivalent and perhaps more concrete way to obtain the
extension is as the unique weak solution
$\eta_{\theta}\in\overline{H}_{0}^{m+2}(\Omega)$ of the boundary
value problem:
\begin{equation*}
\sum_{s=0}^{m+2}(-1)^{s}\Delta^{s}\eta_{\theta}=0\text{ }\text{
}\text{ in }\text{ }\text{ }\Omega,\text{ }\text{ }\text{ }\text{
}\partial_{y}^{s}\eta_{\theta}|_{\partial\Omega_{1}}
=\partial_{y}^{s}u_{\theta}|_{\partial\Omega_{1}},\text{ }\text{
}\text{ }\text{ }0\leq s\leq m+1,
\end{equation*}
\begin{equation*}
\left(\sum_{l=s}^{m+1-s}(-1)^{l}\partial_{y}\Delta^{l}\eta_{\theta}\right)
_{\partial\Omega_{2}}=0,\text{ }\text{ }\text{ }\text{ }0\leq
s\leq\left[\frac{m+1}{2}\right],
\end{equation*}
\begin{equation*}
\left(\sum_{l=s+1}^{m+1-s}(-1)^{l}\Delta^{l}\eta_{\theta}\right)
_{\partial\Omega_{2}}=0,\text{ }\text{ }\text{ }\text{ }0\leq
s\leq\left[\frac{m+(-1)^{m}}{2}\right].
\end{equation*}
Q.E.D.\medskip

  Now in order to solve (4.2) with
$\phi\in\overline{H}_{0}^{m+1}(\partial\Omega_{1})$,
$\psi\in\overline{H}_{0}^{m}(\partial\Omega_{1})$, and
$f\in\overline{H}_{0}^{m}(\Omega)$, we note that if
$\eta_{\theta}\in\overline{H}_{0}^{m+2}(\Omega)$ is as in Lemma
4.2 and $\overline{u}_{\theta}\in H_{0}^{m}(\Omega)$ is given by
Corollary 4.1 with
$\overline{f}_{\theta}=f-L_{\theta}\eta_{\theta}\in
H_{0}^{m}(\Omega)$, then
$u_{\theta}=\overline{u}_{\theta}+\eta_{\theta}\in\overline{H}_{0}^{m}(\Omega)$
satisfies (4.2).  Observe that in order for $u_{\theta}$ to have
the desired vanishing at the origin we require
$m\geq\alpha_{0}+2$.\par
  Our next task shall be to estimate
$u_{\theta}$ independent of $\theta$, in a manner suitable for
application to the Nash iteration of the next section.  A
significant difference between the following theorem and its
analogue for the elliptic regions (Theorem 3.2), is that the loss
of derivatives here depends on the degree to which the Gaussian
curvature vanishes at the origin.\medskip

\textbf{Theorem 4.2.}  \textit{Suppose that $g\in C^{m_{*}}$,
$\phi,\psi\in\overline{C}^{\infty}(\partial\Omega_{1})$,
$f\in\overline{C}^{\infty}(\Omega)$, and $|w|_{C^{2N+4}}<1$ where
$N\leq m_{*}-2$ is the largest integer such that
$\partial^{\alpha}K(0,0)=0$ for all $|\alpha|\leq N$. If
$\alpha_{0}+2\leq m\leq m_{*}-6$ and $\varepsilon=\varepsilon(m)$
is sufficiently small, then there exists a unique solution
$u_{\theta}\in\overline{H}_{0}^{m}(\Omega)$ of (4.2) for each
$\theta>0$.  Furthermore there exists a constant $C_{m}$
independent of $\varepsilon$ and $\theta$ such that}
\begin{eqnarray*}
\parallel u_{\theta}\parallel_{m}&\leq& C_{m}(\parallel
f\parallel_{m+N}+\parallel\phi\parallel_{m+N+1,\partial\Omega_{1}}
+\parallel\psi\parallel_{m+N,\partial\Omega_{1}}\\
& &+\parallel w\parallel_{m+N+6}(\parallel
f\parallel_{N+2}+\parallel\phi\parallel_{N+3,\partial\Omega_{1}}
+\parallel\psi\parallel_{N+2,\partial\Omega_{1}}))
\end{eqnarray*}
\textit{for each $m\leq m_{*}-N-8$.}\medskip

\textit{Proof.}  The existence of a solution
$u_{\theta}\in\overline{H}_{0}^{m_{*}-6}(\Omega)$ for each
$\theta>0$ follows directly from the discussion preceding the
statement of this theorem.  In order to make estimates it will be
advantageous to have a zeroth order term for $L_{\theta}$.
Therefore we set $v_{\theta}=e^{-\frac{1}{2}y^{2}}u_{\theta}$ and
observe that
\begin{equation*}
\overline{L}_{\theta}v_{\theta}:=\partial_{x}(\overline{K}_{\theta}\partial_{x}v_{\theta})
+\partial_{y}^{2}v_{\theta}+\overline{C}\partial_{x}v_{\theta}+(2y+\overline{D})
\partial_{y}v_{\theta}+(1+y^{2}+y\overline{D})v_{\theta}=e^{-\frac{1}{2}y^{2}}f:=\overline{f}.
\end{equation*}
With the aim of treating the $x$-derivatives first, we
differentiate the above equation to find
\begin{eqnarray*}
\overline{L}_{\theta}^{(m)}\partial_{x}^{m}v_{\theta}&=&\partial_{x}^{m}\overline{f}
-\sum_{s=3}^{m}\left(\begin{array}{c}
m \\
s
\end{array}\right)\partial_{x}^{s}\overline{K}_{\theta}\partial_{x}^{m-s}
(\partial_{x}^{2}v_{\theta})\\
& &-\sum_{s=2}^{m}\left(\begin{array}{c}
m \\
s
\end{array}\right)\partial_{x}^{s}(\overline{C}+\partial_{x}\overline{K}_{\theta})
\partial_{x}^{m-s}(\partial_{x}v_{\theta})\\
& &-\sum_{s=1}^{m}\left(\begin{array}{c}
m \\
s
\end{array}\right)[\partial_{x}^{s}\overline{D}\partial_{x}^{m-s}(\partial_{y}v_{\theta})
+y\partial_{x}^{s}\overline{D}\partial_{x}^{m-s}v_{\theta}]\\
&:=&\partial_{x}^{m}\overline{f}+\overline{f}^{(m)}(v_{\theta}),
\end{eqnarray*}
where
\begin{eqnarray*}
\overline{L}_{\theta}^{(m)}v&:=&(\overline{K}_{\theta}v_{x})_{x}+v_{yy}
+(\overline{C}+m\partial_{x}\overline{K}_{\theta})v_{x}
+(2y+\overline{D})v_{y}\\
&
&+\left(1+y^{2}+y\overline{D}+m\partial_{x}\overline{C}+\frac{m(m+1)}{2}\partial_{x}^{2}
\overline{K}_{\theta}\right)v.
\end{eqnarray*}\par
  We first assume that $m\leq N+1$.  In this case let
$\eta_{\theta}$ be given by Lemma 4.2 such that
$\eta_{\theta}|_{\partial\Omega_{1}}=v_{\theta}|_{\partial\Omega_{1}}$,
$\partial_{y}\eta_{\theta}|_{\partial\Omega_{1}}=\partial_{y}v_{\theta}|_{\partial\Omega_{1}}$
and
\begin{equation}
\partial_{y}^{t}(\overline{f}-\overline{L}_{\theta}\eta_{\theta})|_{\partial\Omega_{1}}=0,\text{
}\text{ }\text{ }\text{ }0\leq t\leq m+N.
\end{equation}
Note that this implies that
$\eta_{\theta}\in\overline{H}_{0}^{m+N+3}(\Omega)$,
$\partial^{\alpha}\eta_{\theta}|_{\partial\Omega_{1}}=\partial^{\alpha}v_{\theta}|_{\partial\Omega_{1}}$
for all $|\alpha|\leq m+N+2$, and we have the estimate
\begin{eqnarray}
\parallel\eta_{\theta}\parallel_{m+N+3}&\leq& C(\parallel
f\parallel_{m+N+1}+\parallel\phi\parallel_{m+N+2,\partial\Omega_{1}}
+\parallel\psi\parallel_{m+N+1,\partial\Omega_{1}}\\
& &+\parallel w\parallel_{m+N+7}(\parallel
f\parallel_{2}+\parallel\phi\parallel_{2,\partial\Omega_{1}}
+\parallel\psi\parallel_{2,\partial\Omega_{1}})).\nonumber
\end{eqnarray}
Furthermore the function
$\overline{v}_{\theta}:=v_{\theta}-\eta_{\theta}$ satisfies
\begin{equation*}
\overline{L}_{\theta}^{(m)}\partial_{x}^{m}\overline{v}_{\theta}=
\partial_{x}^{m}(\overline{f}-\overline{L}_{\theta}\eta_{\theta})
+\overline{f}^{(m)}(\overline{v}_{\theta}).
\end{equation*}
As in the proof of Theorem 4.1, we set
$b=\overline{K}_{\theta}^{-1}e^{-\lambda y}$ and integrate by
parts:
\begin{eqnarray}
& &\!\!\!(-b\partial_{y}\partial_{x}^{m}\overline{v}_{\theta},
\overline{L}_{\theta}^{(m)}\partial_{x}^{m}\overline{v}_{\theta})\\
&=&\!\!\!\int_{\Omega}\left[-\frac{1}{2}(b\overline{K}_{\theta})_{y}(\partial_{x}^{m+1}\overline{v}_{\theta})^{2}
+(b_{x}\overline{K}_{\theta}-mb\partial_{x}\overline{K}_{\theta}-b\overline{C})\partial_{x}^{m+1}
\overline{v}_{\theta}
\partial_{y}\partial_{x}^{m}\overline{v}_{\theta}\right]\nonumber\\
& &\!\!\!+\int_{\Omega}
\left(\frac{1}{2}b_{y}-b(2y+\overline{D})\right)(\partial_{y}\partial_{x}^{m}\overline{v}_{\theta})^{2}\nonumber\\
& &\!\!\!+\int_{\Omega}
\frac{1}{2}\left[b(1+y^{2}+y\overline{D}+m\partial_{x}\overline{C}+\frac{m(m+1)}{2}\partial_{x}^{2}\overline{K}
_{\theta})\right]_{y}(\partial_{x}^{m}\overline{v}_{\theta})^{2}\nonumber\\
&
&\!\!\!+\int_{\partial\Omega}\left[\frac{1}{2}b\overline{K}_{\theta}(\partial_{x}^{m+1}\overline{v}_{\theta})^{2}\nu_{2}
-b\overline{K}_{\theta}\partial_{x}^{m+1}\overline{v}_{\theta}\partial_{y}\partial_{x}^{m}
\overline{v}_{\theta}\nu_{1}\right]\nonumber\\
&
&\!\!\!-\int_{\partial\Omega}\left[\frac{1}{2}b(\partial_{y}\partial_{x}^{m}\overline{v}_{\theta})^{2}\nu_{2}
+\frac{1}{2}b\left(1+y^{2}+y\overline{D}
+m\partial_{x}\overline{C}+\frac{m(m+1)}{2}\partial_{x}^{2}\overline{K}_{\theta}\right)
(\partial_{x}^{m}\overline{v}_{\theta})^{2}\nu_{2}\right].\nonumber
\end{eqnarray}
The boundary integral along $\partial\Omega_{2}$ is nonnegative,
and according to the choice of $\eta_{\theta}$ it vanishes along
$\partial\Omega_{1}$.  Moreover, the same calculations as in the
proof of Theorem 4.1 apply to the interior integral to yield
\begin{eqnarray*}
&
&\lambda(\parallel\partial_{x}^{m+1}\overline{v}_{\theta}\parallel+
\parallel\sqrt{|b|}\partial_{y}\partial_{x}^{m}\overline{v}_{\theta}\parallel
+\parallel\sqrt{|b|}\partial_{x}^{m}\overline{v}_{\theta}\parallel)\\
&\leq&
C(\parallel\sqrt{|b|}\partial_{x}^{m}(\overline{f}-\overline{L}_{\theta}\eta_{\theta})
\parallel
+\parallel\sqrt{|b|}\overline{f}^{(m)}(\overline{v}_{\theta})\parallel).
\end{eqnarray*}\par
  We proceed to estimate each term on the right-hand side
separately.  First note that since $m\leq N+1$ and
$|w|_{C^{N+4}}<1$, we have
\begin{equation*}
\parallel\sqrt{|b|}\overline{f}^{(m)}(\overline{v}_{\theta})\parallel^{2}
\leq C_{m}\int_{\Omega}\sum_{s=0}^{m-1}e^{-\lambda
y}|\overline{K}_{\theta}|^{-1}[(\partial_{x}^{s}
\overline{v}_{\theta})^{2}+(\partial_{y}\partial_{x}^{s}
\overline{v}_{\theta})^{2}].
\end{equation*}
Next observe that since $\overline{K}$ vanishes (at most) to order
$N$ at the origin, there exists a constant $C_{0}>0$ such that
$|\overline{K}|\geq C_{0}^{-1}(y-h(x))^{N+1}$ in $\Omega$. Then in
light of (4.11), a little calculation shows that
\begin{eqnarray*}
\parallel\sqrt{|b|}\partial_{x}^{m}(\overline{f}-\overline{L}_{\theta}\eta_{\theta})\parallel^{2}
&\leq&
C_{0}\int_{\Omega}\frac{[\partial_{x}^{m}(\overline{f}-\overline{L}_{\theta}\eta_{\theta})]^{2}}{(y-h(x))^{N+1}}\\
&\leq&
C_{1}\int_{\Omega}[\partial_{y}^{N+1}\partial_{x}^{m}(\overline{f}-\overline{L}_{\theta}\eta_{\theta})]^{2}.
\end{eqnarray*}
It follows that applying (4.12) and summing from 0 to $m$ produces
\begin{eqnarray}
&
&\sum_{s=0}^{m+1}\parallel\partial_{x}^{s}\overline{v}_{\theta}\parallel
+\sum_{s=0}^{m}\parallel\sqrt{|b|}\partial_{y}\partial_{x}^{s}\overline{v}_{\theta}\parallel
+\sum_{s=0}^{m}\parallel\sqrt{|b|}\partial_{x}^{s}\overline{v}_{\theta}\parallel\\
&\leq& C_{m}(\parallel f\parallel_{m+N+1}+\parallel
\phi\parallel_{m+N+2,\partial\Omega_{1}}+\parallel
\psi\parallel_{m+N+1,\partial\Omega_{1}}\nonumber\\
& &+\parallel w\parallel_{m+N+7}(\parallel
f\parallel_{2}+\parallel\phi\parallel_{2,\partial\Omega_{1}}
+\parallel\psi\parallel_{2,\partial\Omega_{1}})),\nonumber
\end{eqnarray}
if $\lambda$ is sufficiently large.  From this inequality we may
obtain an estimate for the $x$-derivatives of $v_{\theta}$ (and
hence for $u_{\theta}$), with the help of (4.12).  In addition, by
solving for $\partial_{y}^{2}u_{\theta}$ in equation (4.2) all
remaining derivatives up to and including order $m$ may also be
estimated.\par
  We now assume that $m\geq N+2$.  In order to isolate terms involving high order
derivatives of $w$ we break $\overline{f}^{(m)}$ into two parts.
Let
$v_{\theta}^{(s)}:=\partial_{x}^{s}v_{\theta}-\eta_{\theta}^{(s)}$
with $\eta_{\theta}^{(s)}$ to be given below, then
$\overline{f}^{(m)}=\overline{f}_{1}^{(m)}+\overline{f}_{2}^{(m)}$
where
\begin{eqnarray*}
\overline{f}_{1}^{(m)}(v_{\theta})&=&
-\sum_{s=N+2}^{m}\left(\begin{array}{c}
m \\
s
\end{array}\right)\partial_{x}^{s}\overline{K}_{\theta}\partial_{x}^{m-s+2}v_{\theta}
-\sum_{s=N+1}^{m}\left(\begin{array}{c}
m \\
s
\end{array}\right)\partial_{x}^{s}(\overline{C}+\partial_{x}\overline{K}_{\theta})
\partial_{x}^{m-s+1}v_{\theta}\\
& &-\sum_{s=N+1}^{m}\left(\begin{array}{c}
m \\
s
\end{array}\right)\partial_{x}^{s}\overline{D}\partial_{y}\partial_{x}^{m-s}v_{\theta}
-y\sum_{s=N}^{m}\left(\begin{array}{c}
m \\
s
\end{array}\right)\partial_{x}^{s}\overline{D}\partial_{x}^{m-s}v_{\theta}\\
& &-\sum_{s=3}^{N+1}\left(\begin{array}{c}
m \\
s
\end{array}\right)\partial_{x}^{s}\overline{K}_{\theta}\eta_{\theta}^{(m-s+2)}
-\sum_{s=2}^{N}\left(\begin{array}{c}
m \\
s
\end{array}\right)\partial_{x}^{s}(\overline{C}+\partial_{x}\overline{K}_{\theta})
\eta_{\theta}^{(m-s+1)}\\
& &-\sum_{s=1}^{N}\left(\begin{array}{c}
m \\
s
\end{array}\right)\partial_{x}^{s}\overline{D}\partial_{y}\eta_{\theta}^{(m-s)}
-y\sum_{s=1}^{N-1}\left(\begin{array}{c}
m \\
s
\end{array}\right)\partial_{x}^{s}\overline{D}\eta_{\theta}^{(m-s)}
\end{eqnarray*}
and
\begin{eqnarray*}
\overline{f}_{2}^{(m)}(v_{\theta})&=&
-\sum_{s=3}^{N+1}\left(\begin{array}{c}
m \\
s
\end{array}\right)\partial_{x}^{s}\overline{K}_{\theta}v_{\theta}^{(m-s+2)}
-\sum_{s=2}^{N}\left(\begin{array}{c}
m \\
s
\end{array}\right)\partial_{x}^{s}(\overline{C}+\partial_{x}\overline{K}_{\theta})
v_{\theta}^{(m-s+1)}\\
& &-\sum_{s=1}^{N}\left(\begin{array}{c}
m \\
s
\end{array}\right)\partial_{x}^{s}\overline{D}\partial_{y}v_{\theta}^{(m-s)}
-y\sum_{s=1}^{N-1}\left(\begin{array}{c}
m \\
s
\end{array}\right)\partial_{x}^{s}\overline{D}v_{\theta}^{(m-s)}.
\end{eqnarray*}
The functions
$\eta_{\theta}^{(s)}\in\overline{H}_{0}^{N+3}(\Omega)$, $0\leq
s\leq m$, are defined recursively in the following way.  For
$0\leq s\leq N+1$ we set
$\eta_{\theta}^{(s)}=\partial_{x}^{s}\eta_{\theta}$, and if
$N+2\leq s\leq m$ we apply Lemma 4.2 to obtain
$\eta_{\theta}^{(s)}$ such that
$\eta_{\theta}^{(s)}|_{\partial\Omega_{1}}=\partial_{x}^{s}v_{\theta}|_{\partial\Omega_{1}}$,
$\partial_{y}\eta_{\theta}^{(s)}|_{\partial\Omega_{1}}
=\partial_{y}\partial_{x}^{s}v_{\theta}|_{\partial\Omega_{1}}$
with
\begin{equation}
\partial_{y}^{t}(\partial_{x}^{s}\overline{f}+\overline{f}^{(s)}_{1}(v_{\theta})
-\overline{L}_{\theta}^{(s)}\eta_{\theta}^{(s)})|_{\partial\Omega_{1}}=0,\text{
}\text{ }\text{ }\text{ }0\leq t\leq N.
\end{equation}
Note that since
$\overline{f}^{(s)}_{1}\in\overline{H}_{0}^{\min(m_{*}-s-6,N+2)}(\Omega)$
and
$\partial_{x}^{s}v_{\theta}|_{\partial\Omega_{1}}\in\overline{H}_{0}^{m_{*}-s-7}(\partial\Omega_{1})$,
$\partial_{y}\partial_{x}^{s}v_{\theta}|_{\partial\Omega_{1}}\in\overline{H}_{0}^{m_{*}-s-8}(\partial\Omega_{1})$,
we must have $N+1\leq m_{*}-s-8$ for the construction of
$\eta_{\theta}^{(s)}$ to be valid.  In this case the following
estimate holds
\begin{eqnarray}
\parallel\eta_{\theta}^{(s)}\parallel_{N+3}\!\!\!&\leq&\!\!\!
C(\parallel\partial_{x}^{s}\overline{f}+\overline{f}^{(s)}_{1}(v_{\theta})\parallel_{N+1}
+\parallel\partial_{x}^{s}v_{\theta}\parallel_{N+2,\partial\Omega_{1}}
+\parallel\partial_{y}\partial_{x}^{s}v_{\theta}\parallel_{N+1,\partial\Omega_{1}}\\
& &\!\!\!+\parallel
w\parallel_{N+7}(\parallel\partial_{x}^{s}\overline{f}+\overline{f}^{(s)}_{1}(v_{\theta})\parallel_{2}
+\parallel\partial_{x}^{s}v_{\theta}\parallel_{2,\partial\Omega_{1}}
+\parallel\partial_{y}\partial_{x}^{s}v_{\theta}\parallel_{2,\partial\Omega_{1}}))\nonumber\\
&\leq&\!\!\!
C_{s}(\parallel\overline{f}^{(s)}_{1}(v_{\theta})\parallel_{N+1}+\parallel
f\parallel_{s+N+1}+\parallel\phi\parallel_{s+N+2,\partial\Omega_{1}}
+\parallel\psi\parallel_{s+N+1,\partial\Omega_{1}}\nonumber\\
& &\!\!\!+\parallel w\parallel_{s+N+7}(\parallel
f\parallel_{2}+\parallel\phi\parallel_{2,\partial\Omega_{1}}
+\parallel\psi\parallel_{2,\partial\Omega_{1}})),\nonumber
\end{eqnarray}
where we have used $|w|_{C^{N+7}}<1$ and the proof of Lemma 4.2 to
estimate
$\partial_{y}^{t}\partial_{x}^{s}v_{\theta}|_{\partial\Omega_{1}}$,
$t=0,1$.\par
  Observe that $v_{\theta}^{(s)}$, $N+2\leq s\leq m$, satisfies
\begin{equation*}
\overline{L}_{\theta}^{(s)}v_{\theta}^{(s)}
=(\partial_{x}^{s}\overline{f}+\overline{f}^{(s)}_{1}(v_{\theta})-\overline{L}_{\theta}^{(s)}\eta_{\theta}^{(s)})
+\overline{f}^{(s)}_{2}(v_{\theta}).
\end{equation*}
Therefore (4.13) applies to yield
\begin{eqnarray}
& &\lambda(\parallel\partial_{x}v_{\theta}^{(s)}\parallel+
\parallel\sqrt{|b|}\partial_{y}v_{\theta}^{(s)}\parallel
+\parallel\sqrt{|b|}v_{\theta}^{(s)}\parallel)\\
&\leq&
C(\parallel\sqrt{|b|}(\partial_{x}^{s}\overline{f}+\overline{f}_{1}^{(s)}(v_{\theta})
-\overline{L}^{(s)}_{\theta}\eta_{\theta}^{(s)})\parallel
+\parallel\sqrt{|b|}\overline{f}_{2}^{(s)}(v_{\theta})\parallel).\nonumber
\end{eqnarray}
We now proceed to estimate each term on the right-hand side of
(4.17) separately.  First note that since $|w|_{C^{N+4}}<1$, we
have
\begin{eqnarray}
\parallel\sqrt{|b|}\overline{f}_{2}^{(s)}(v_{\theta})\parallel^{2}
&\leq& C_{s}\int_{\Omega}\sum_{l=0}^{s-1}e^{-\lambda
y}|\overline{K}_{\theta}|^{-1}[(
v_{\theta}^{(l)})^{2}+(\partial_{y}v_{\theta}^{(l)})^{2}]\\
&=&C_{s}\int_{\Omega}\sum_{l=N+2}^{s-1}e^{-\lambda
y}|\overline{K}_{\theta}|^{-1}[(
v_{\theta}^{(l)})^{2}+(\partial_{y}v_{\theta}^{(l)})^{2}]\nonumber\\
& &+C_{s}\int_{\Omega}\sum_{l=0}^{N+1}e^{-\lambda
y}|\overline{K}_{\theta}|^{-1}[(\partial_{x}^{l}
\overline{v}_{\theta})^{2}+(\partial_{y}\partial_{x}^{l}
\overline{v}_{\theta})^{2}].\nonumber
\end{eqnarray}
Next observe again that since $|\overline{K}|\geq
C_{0}^{-1}(y-h(x))^{N+1}$ in $\Omega$, (4.15) implies that
\begin{eqnarray}
\parallel\sqrt{|b|}(\partial_{x}^{s}\overline{f}+\overline{f}_{1}^{(s)}(v_{\theta})
-\overline{L}^{(s)}_{\theta}\eta_{\theta}^{(s)})\parallel^{2}
&\leq&
C_{0}\int_{\Omega}\frac{[\partial_{x}^{s}\overline{f}+\overline{f}_{1}^{(s)}(v_{\theta})
-\overline{L}^{(s)}_{\theta}\eta_{\theta}^{(s)}]^{2}}{(y-h(x))^{N+1}}\\
&\leq&
C_{1}\int_{\Omega}[\partial_{y}^{N+1}(\partial_{x}^{s}\overline{f}+\overline{f}_{1}^{(s)}(v_{\theta})
-\overline{L}^{(s)}_{\theta}\eta_{\theta}^{(s)})]^{2}.\nonumber
\end{eqnarray}
Furthermore applying the Gagliardo-Nirenberg inequalities (Lemma
5.2), (4.12), (4.14), (4.16), and using $|w|_{C^{2N+4}}<1$
produces
\begin{eqnarray}
\parallel\overline{f}_{1}^{(s)}(v_{\theta})\parallel_{N+1}
&\leq& C_{s}((|\overline{K}_{\theta}|_{C^{N+2}}
+|\overline{C}|_{C^{N+1}}+|\overline{D}|_{C^{N+1}})
\parallel v_{\theta}\parallel_{s+1}\\
& &+(\parallel\overline{K}_{\theta}\parallel_{s+N+3}
+\parallel\overline{C}\parallel_{s+N+2}+
\parallel\overline{D}\parallel_{s+N+2})|v_{\theta}|_{C^{0}})\nonumber\\
&
&+C_{s}\left(\sum_{l=s-N}^{N+1}\parallel\eta_{\theta}\parallel_{l+N+2}
+\sum_{l=N+2}^{s-1}\parallel\eta_{\theta}^{(l)}\parallel_{N+2}\right)\nonumber\\
&\leq& C_{s}(\parallel f\parallel_{s+N+1}+\varepsilon\parallel
v_{\theta}\parallel_{s+1}+\parallel\phi\parallel_{s+N+2,\partial\Omega_{1}}
+\parallel\psi\parallel_{s+N+1,\partial\Omega_{1}}\nonumber\\
& &+\parallel w\parallel_{s+N+7}(\parallel
f\parallel_{N+2}+\parallel\phi\parallel_{N+3,\partial\Omega_{1}}
+\parallel\psi\parallel_{N+2,\partial\Omega_{1}})).\nonumber
\end{eqnarray}
It follows that we may combine (4.17)-(4.20) and utilize (4.12),
(4.14), as well as (4.16) to obtain
\begin{eqnarray*}
& &\sum_{s=N+2}^{m}\parallel
\partial_{x}v_{\theta}^{(s)}\parallel
+\sum_{s=N+2}^{m}\parallel\sqrt{|b|}\partial_{y}v_{\theta}^{(s)}\parallel
+\sum_{s=N+2}^{m}\parallel\sqrt{|b|}v_{\theta}^{(s)}\parallel\\
&\leq& C_{m}(\parallel f\parallel_{m+N+1}+\varepsilon\parallel
u_{\theta}\parallel_{m+1}+\parallel
\phi\parallel_{m+N+2,\partial\Omega_{1}}+\parallel
\psi\parallel_{m+N+1,\partial\Omega_{1}}\nonumber\\
& &+\parallel w\parallel_{m+N+7}(\parallel
f\parallel_{N+2}+\parallel\phi\parallel_{N+3,\partial\Omega_{1}}
+\parallel\psi\parallel_{N+2,\partial\Omega_{1}})).\nonumber
\end{eqnarray*}
Since
$v_{\theta}^{(s)}=\partial_{x}^{s}v_{\theta}-\eta_{\theta}^{(s)}$,
with the help of (4.16) the above inequality yields an estimate
for the $x$-derivatives of $u_{\theta}$.  The remaining
derivatives of $u_{\theta}$ may be estimated in the usual way, by
solving for $\partial_{y}^{2}u_{\theta}$ from equation (4.2).
Lastly taking $\varepsilon=\varepsilon(m)$ sufficiently small
yields the desired result.  Q.E.D.\medskip

\textbf{Remark 4.1.} \textit{In this section we have focused on
the Cauchy problem in the domains $\Omega_{\varrho}^{-}\cap
B_{\sigma}(0)$, $\varrho=1,2$. However analogous existence and
regularity results, requiring only slight modifications of the
arguments above, hold for the Cauchy problem in the domains
$\Omega_{\varrho}^{-}\cap B_{\sigma}(0)$, $\varrho=3,4$, when the
Cauchy data are prescribed on either the ``upper" or ``lower"
parts of these domains (that is, on one of the two differentiable
components of the curve $x=h(y)$).}

\begin{center}
\textbf{5. The Nash-Moser Iteration}
\end{center}\setcounter{equation}{0}
\setcounter{section}{5}

  In this section we will carry out Nash-Moser type iteration
procedures to obtain solutions of (2.1) in each of the elliptic,
hyperbolic, and mixed type regions separately.  The solutions will
then be patched together by choosing appropriate boundary values,
to yield a solution on a full neighborhood of the origin.  As a
consequence of Lemma 2.1, we can assume (by a judicious choice of
coordinates) that each elliptic region is given by
\begin{equation*}
\Omega_{\kappa}^{+}=\{(\xi^{1},\xi^{2})\mid
0<\xi^{2}<(\tan\delta)\xi^{1},\!\text{ } |\xi|<\sigma\},\text{
}\text{ }\text{ }\text{ }1\leq\kappa\leq\kappa_{0},
\end{equation*}
each hyperbolic region is given by
\begin{equation*}
\Omega_{\varrho}^{-}=\{(\xi^{1},\xi^{2})\mid
h(\xi^{1})<\xi^{2}<\sigma\},\text{ }\text{ }\text{ }\text{
}\varrho=1,2,
\end{equation*}
or
\begin{equation*}
\Omega_{\varrho}^{-}= \{(\xi^{1},\xi^{2})\mid
h(\xi^{2})<\xi^{1}<\sigma\},\text{ }\text{ }\text{ }\text{
}\varrho=3,4.
\end{equation*}
If $\overline{\partial}\Omega_{\kappa}^{+}$ denotes the portion of
the boundary consisting of the curves $\xi^{2}=0$ and
$\xi^{2}=(\tan\delta)\xi^{1}$,
$\overline{\partial}\Omega_{\varrho}^{-}$, $\varrho=1,2$, denotes
the portion of the boundary given by $\xi^{2}=h(\xi^{1})$, and
$\overline{\partial}\Omega_{\varrho}^{-}$, $\varrho=3,4$, denotes
either the upper or lower part of the boundary curve
$\xi^{1}=h(\xi^{2})$, then we aim to solve
\begin{equation}
\Phi(w_{\kappa}^{+})=0\text{ }\text{ }\text{ in }\text{ }\text{
}\Omega_{\kappa}^{+},\text{ }\text{ }\text{ }\text{
}w_{\kappa}^{+}|_{\overline{\partial}\Omega_{\kappa}^{+}}=0,
\end{equation}
\begin{equation*}
\partial^{\alpha}w_{\kappa}^{+}(0,0)=0,\text{ }\text{ }\text{
}\text{ }|\alpha|\leq\alpha_{0},
\end{equation*}
for each $\kappa=1,\ldots,\kappa_{0}$, and
\begin{equation}
\Phi(w_{\varrho}^{-})=0\text{ }\text{ }\text{ in }\text{ }\text{
}\Omega_{\varrho}^{-},\text{ }\text{ }\text{ }\text{
}w_{\varrho}^{-}|_{\overline{\partial}\Omega_{\varrho}^{-}}=\phi_{\varrho}^{-},\text{
}\text{ }\text{ }\text{
}\partial_{\nu}w_{\varrho}^{-}|_{\overline{\partial}\Omega_{\varrho}^{-}}=\psi_{\varrho}^{-},
\end{equation}
\begin{equation*}
\partial^{\alpha}w_{\varrho}^{-}(0,0)=0,\text{ }\text{ }\text{
}\text{ }|\alpha|\leq\alpha_{0},
\end{equation*}
for each $\varrho=1,\ldots,\varrho_{0}$, where $\partial_{\nu}$
denotes the outward normal derivative, $\alpha_{0}$ is a large
integer, and $\phi_{\varrho}^{-}$, $\psi_{\varrho}^{-}$ will be
specified below.\par
    Both of the problems (5.1) and (5.2) will require slight
modifications of the standard Nash-Moser procedure.  This arises
from the fact that instead of solving the linearized equation at
each iteration, the theories developed in sections $\S 3$ and $\S
4$ require us to solve modified versions of the linearized
equation.  However the error incurred at each step by these
modifications is of quadratic type and therefore does not prevent
the procedure from converging to a solution.  Below we shall carry
out the proof in full detail for the hyperbolic regions, problem
(5.2).  Moreover, since the corresponding iteration for problem
(5.1) differs only slightly from that of (5.2), we shall merely
indicate the necessary changes required in this case.\bigskip

\textsc{Hyperbolic Regions}\bigskip

  In the hyperbolic regions, the linearization $\mathcal{L}(w)$
and the operator $L_{\theta}(w)$ that we invert in section $\S 4$
differ by perturbation terms coming from Lemma 2.1 as well as a
regularizing term involving $\theta=|\Phi(w)|_{C^{1}}$.  More
precisely according to Lemma 2.1
\begin{eqnarray}
\mathcal{L}(w)u&=&\varepsilon
a^{22}(w)L_{\theta}(w)u+\varepsilon\theta
a^{22}(w)\partial_{\xi^{1}}^{2}u\\
&
&+\varepsilon(a^{22}(w))^{-1}\Phi(w)[\partial_{x^{1}}^{2}u-\partial_{x^{1}}\log
(a^{22}(w)\sqrt{|g|})\partial_{x^{1}}u],\nonumber
\end{eqnarray}
where $\xi^{i}(x^{1},x^{2})$ are the coordinates constructed in
Lemma 2.1. Furthermore, as with all Nash-Moser iterative schemes
we will need smoothing operators.  Since the theory of section $\S
4$ is based on the Sobolev spaces $\overline{H}_{0}^{m}$, the
smoothing operators we employ should respect these spaces.  For
convenience, in the remainder of this section the hyperbolic
region in question $\Omega_{\varrho}^{-}$ will be denoted by
$\Omega$, and the $\overline{H}_{0}^{m}(\Omega)$ norm will be
denoted by $\parallel\cdot\parallel_{m}$.\medskip

\textbf{Lemma 5.1.}  \textit{Given $\mu\geq 1$ there exists a
linear smoothing operator
$S_{\mu}:L^{2}(\Omega)\rightarrow\overline{H}_{0}^{\infty}(\Omega)$
such that for all $l,m\in\mathbb{Z}_{\geq 0}$ and $u\in
\overline{H}_{0}^{l}(\Omega)$,}\smallskip

$i)$ $\parallel S_{\mu}u\parallel_{m}\leq C_{l,m}\parallel
u\parallel_{l}$, $m\leq l$,\smallskip

$ii)$ $\parallel S_{\mu}u\parallel_{m}\leq
C_{l,m}\mu^{m-l}\parallel u\parallel_{l}$, $l\leq m$,\smallskip

$iii)$ $\parallel u- S_{\mu}u\parallel_{m}\leq
C_{l,m}\mu^{m-l}\parallel u\parallel_{l}$, $m\leq l$.\smallskip

\noindent\textit{Furthermore there exists a linear smoothing
operator $S^{'}_{\mu}:L^{2}(\Omega)\rightarrow H^{\infty}(\Omega)$
such that $(i)$, $(ii)$, and $(iii)$ hold whenever $u\in
H^{l}(\Omega)$.}\medskip

\textit{Proof.}  See appendix B.  Q.E.D.\medskip

  The next lemma contains the so called Gagliardo-Nirenberg inequalities,
which will be used frequently throughout this section.\medskip

\textbf{Lemma 5.2.} \textit{Let $u,v\in
C^{k}(\overline{\Omega})$.}

$i)$\textit{  If $\alpha$ and $\beta$ are multi-indices such that
$|\alpha|+|\beta|=m$, then there exists a constant $C_{1}$
depending on $m$ such that}
\begin{equation*}
\parallel\partial^{\alpha}u\partial^{\beta}v\parallel_{L^{2}(\Omega)}\leq
C_{1}(|u|_{L^{\infty}(\Omega)}\parallel
v\parallel_{H^{m}(\Omega)}+\parallel
u\parallel_{H^{m}(\Omega)}|v|_{L^{\infty}(\Omega)}).
\end{equation*}

$ii)$\textit{  If $\alpha_{1},\ldots,\alpha_{l}$ are multi-indices
such that $|\alpha_{1}|+\cdots+|\alpha_{l}|=m$, then there exists
a constant $C_{2}$ depending on $l$ and $m$ such that}
\begin{equation*}
\parallel\partial^{\alpha_{1}}u_{1}\cdots\partial^{\alpha_{l}}u_{l}
\parallel_{L^{2}(\Omega)}\leq
C_{2}\sum_{j=1}^{l}(|u_{1}|_{L^{\infty}(\Omega)}\cdots
\widehat{|u_{j}|}_{L^{\infty}(\Omega)}
\cdots|u_{l}|_{L^{\infty}(\Omega)})\parallel
u_{j}\parallel_{H^{m}(\Omega)},
\end{equation*}
\textit{where $\widehat{|u_{j}|}_{L^{\infty}(\Omega)}$ indicates
the absence of $|u_{j}|_{L^{\infty}(\Omega)}$.}

$iii)$\textit{  Let $\mathcal{D}\subset\mathbb{R}^{l}$ be compact
and contain the origin, and let $G\in C^{\infty}(\mathcal{D})$. If
$u\in H^{m}(\Omega,\mathcal{D})\cap
L^{\infty}(\Omega,\mathcal{D})$, then there exists a constant
$C_{3}$ depending on $m$ such that}
\begin{equation*}
\parallel G\circ u\parallel_{H^{m}(\Omega)}\leq C_{3}|u|_{L^{\infty}(\Omega)}(|G(0)|+\parallel
u\parallel_{H^{m}(\Omega)}).
\end{equation*}

\textit{Proof.}  These estimates are standard consequences of the
interpolation inequalities, and may be found in, for instance, [19].
Q.E.D.\medskip

  We now set up the underlying iterative procedure.  Suppose that
$\phi_{\varrho}^{-}\in\overline{H}_{0}^{m_{*}-m_{0}+1}(\overline{\partial}\Omega)$
and
$\psi_{\varrho}^{-}\in\overline{H}_{0}^{m_{*}-m_{0}}(\overline{\partial}\Omega)$
for some $m_{0}\geq 0$.  Then according to the proof of Lemma 4.2
there exists $w_{0}\in\overline{H}_{0}^{m_{*}-m_{0}+2}(\Omega)$
such that
\begin{equation}
w_{0}|_{\overline{\partial}\Omega}=\phi_{\varrho}^{-},\text{
}\text{ }\text{ }\text{
}\partial_{\nu}w_{0}|_{\overline{\partial}\Omega}=\psi_{\varrho}^{-}.
\end{equation}
Now suppose that in addition to $w_{0}$, functions
$w_{1},w_{2},\ldots,w_{n}$ have been defined on $\Omega$, and put
$v_{i}=S_{i}w_{i}$, $0\leq i\leq n$, where $S_{i}=S_{\mu^{i}}$.
Then define $w_{n+1}=w_{n}+u_{n}$ where $u_{n}$ is the unique
solution of
\begin{equation}
L_{\theta_{n}}(v_{n})u_{n}=f_{n}\text{ }\text{ }\text{ in }\text{
}\text{ }\Omega,\text{ }\text{ }\text{
}u_{n}|_{\overline{\partial}\Omega}=\partial_{\nu}u_{n}|_{\overline{\partial}\Omega}=0,
\end{equation}
given by Theorem 4.2, where $\theta_{n}=|\Phi(v_{n})|_{C^{1}}$ and
$f_{n}$ will be specified below.\par
  Let $Q_{n}(w_{n},u_{n})$ denote the quadratic error in the
Taylor expansion of $\Phi$ at $w_{n}$.  Then according to (5.3) we
have
\begin{eqnarray}
\Phi(w_{n+1})&=&\Phi(w_{n})+\mathcal{L}(w_{n})u_{n}+Q_{n}(w_{n},u_{n})\\
&=&\Phi(w_{n})+\varepsilon
S^{'}_{n}a^{22}(v_{n})L_{\theta_{n}}(v_{n})u_{n}+e_{n},\nonumber
\end{eqnarray}
where
\begin{eqnarray*}
e_{n}&=&(\mathcal{L}(w_{n})-\mathcal{L}(v_{n}))u_{n}+\varepsilon(I-S_{n}^{'})a^{22}(v_{n})L_{\theta_{n}}(v_{n})u_{n}
+\varepsilon\theta_{n}a^{22}(v_{n})\partial_{\xi_{n}^{1}}^{2}u_{n}\\
&
&+Q_{n}(w_{n},u_{n})+\varepsilon(a^{22}(v_{n}))^{-1}\Phi(v_{n})[\partial_{x^{1}}^{2}u_{n}-\partial_{x^{1}}(\log
a^{22}(v_{n})\sqrt{|g|})\partial_{x^{1}}u_{n}],
\end{eqnarray*}
and $\xi_{n}^{i}$ are the coordinates of Lemma 2.1 with respect to
$v_{n}$.\par
  We now define $f_{n}$.  In order to solve (5.5) with the theory
of section $\S 4$, we require
$f_{n}\in\overline{C}^{\infty}(\Omega)$.  Furthermore, we need the
right-hand side of (5.6) to tend to zero sufficiently fast, to
make up for the error incurred at each step by solving (5.5)
instead of solving the unmodified linearized equation. Therefore
we set $E_{0}=0$, $E_{n}=\sum_{i=0}^{n-1}e_{i}$, and define
\begin{equation}
f_{0}=-[\varepsilon S_{0}^{'}a^{22}(v_{0})]^{-1}S_{0}\Phi(w_{0}),
\end{equation}
\begin{equation*}
f_{n}=[\varepsilon
S_{n}^{'}a^{22}(v_{n})]^{-1}(S_{n-1}E_{n-1}-S_{n}E_{n}+(S_{n-1}-S_{n})\Phi(w_{0})).
\end{equation*}
It follows that
\begin{eqnarray}
\Phi(w_{n+1})&=&\Phi(w_{0})+\sum_{i=0}^{n}\varepsilon
S_{i}^{'}a^{22}(v_{i})f_{i}+E_{n}+e_{n}\\
&=&(I-S_{n})\Phi(w_{0})+(I-S_{n})E_{n}+e_{n}.\nonumber
\end{eqnarray}
The following theorem contains the Moser estimate for solutions of
(5.5), upon which the whole iteration scheme is based.\medskip

\textbf{Theorem 5.1.}  \textit{Suppose that $g\in C^{m_{*}}$ and
$N$ is as in Theorem 4.2.  If $m\leq m_{*}-N-8$,
$|v_{n}|_{C^{2N+4}}<1$, and $\varepsilon=\varepsilon(m)$ is
sufficiently small then there exists a unique solution
$u_{n}\in\overline{H}^{m}_{0}(\Omega)$ of (5.5) which satisfies
the estimate}
\begin{equation}
\parallel u_{n}\parallel_{m}\leq C_{m}(\parallel
f_{n}\parallel_{m+N}+\parallel v_{n}\parallel_{m+N+6}
\parallel f_{n}\parallel_{N+2}),
\end{equation}
\textit{for some constant $C_{m}$ independent of $\varepsilon$ and
$\theta_{n}$.}\medskip

\textit{Proof.}  This will follow from Theorem 4.2 with
$\phi=\psi\equiv 0$.  The only difference is that the Sobolev
norms appearing in Theorem 4.2 are with respect to the coordinates
$\xi_{n}^{i}(x^{1},x^{2})$ of Lemma 2.1.  In order to obtain the
current estimate from that of Theorem 4.2 we may utilize (2.4).
Note that $\delta$ (of Lemma 2.1) is not chosen arbitrarily small
in the hyperbolic regions, and so it does not appear in the above
estimate.  Q.E.D.\medskip

  In what follows, we will show that the right-hand side of (5.8)
tends to zero sufficiently fast to guarantee convergence of
$\{w_{n}\}_{n=0}^{\infty}$ to a solution of (5.2).  Let $\rho$ be
a positive number that will be chosen as large as possible, and
set $\mu=\varepsilon^{-\frac{1}{2\rho}}$, $\mu_{n}=\mu^{n}$.
Furthermore, note that
$\Phi(w_{0})\in\overline{H}_{0}^{m_{*}-m_{0}}(\Omega)$. The
convergence of $\{w_{n}\}_{n=0}^{\infty}$ will follow from the
following eight statements, valid for $0\leq m\leq m_{*}-N-8$
unless specified otherwise.  These statements shall be proven by
induction on $n$, for some constants $C_{1},\ldots,C_{5}$
independent of $n$ and $\varepsilon$, but dependent on
$m$.\bigskip

   I$_{n}$:   $\parallel u_{n-1}\parallel_{m}\leq\varepsilon\mu_{n-1}^{m+N+2-\rho}$,\bigskip

   II$_{n}$:   $\parallel w_{n}\parallel_{m}\leq\begin{cases}
C_{1}\varepsilon & \text{if $m+N+2-\rho\leq-1/2$},\\
C_{1}\varepsilon\mu_{n}^{m+N+2-\rho} & \text{if $m+N+2-\rho\geq
1/2$,}
\end{cases}$\bigskip

   III$_{n}$:   $\parallel w_{n}\parallel_{2N+6}\leq C_{1}\varepsilon,\text{ }
   \text{ }\parallel v_{n}\parallel_{2N+6}\leq
   C_{3}\varepsilon$,\bigskip

   IV$_{n}$:   $\parallel w_{n}-v_{n}\parallel_{m}\leq
   C_{2}\varepsilon\mu_{n}^{m+N+2-\rho}$,\bigskip

   V$_{n}$:   $\parallel v_{n}\parallel_{m}\leq\begin{cases}
C_{3}\varepsilon & \text{if $m+N+2-\rho\leq-1/2$},\\
C_{3}\varepsilon\mu_{n}^{m+N+2-\rho} & \text{if $m+N+2-\rho\geq
1/2$,}
\end{cases}$ $m<\infty$,\bigskip

   VI$_{n}$:   $\parallel e_{n-1}\parallel_{m}\leq
   C_{4}\varepsilon^{3}\mu_{n-1}^{m-\rho}$,\text{ } $m\leq\min(m_{*}-N-10,m_{*}-m_{0})$,\bigskip

   VII$_{n}$:   $\parallel f_{n}\parallel_{m}
   \leq C_{5}\varepsilon^{2}(1+\mu^{\rho-m})\mu_{n}^{m-\rho}$, $m<\infty$,\bigskip

   VIII$_{n}$:
   $\parallel\Phi(w_{n})\parallel_{m}+\parallel\Phi(v_{n})\parallel_{m}
   \leq\varepsilon\mu_{n}^{m+N+4-\rho}$, $m\leq\min(m_{*}-N-10,m_{*}-m_{0})$.\bigskip

  Assume that the above eight statements hold for all nonnegative integers less than or equal
to $n$.  The next four propositions will show that they also hold
for $n+1$.  The case $n=0$ will be proven shortly
thereafter.\medskip

\textbf{Proposition 5.1.}  \textit{If $3N+8< \rho< m_{*}-6$,
$0\leq m\leq m_{*}-N-8$, and $\varepsilon$ is sufficiently small,
then} I$_{n+1}$, II$_{n+1}$, III$_{n+1}$, IV$_{n+1}$, \textit{and}
V$_{n+1}$ \textit{hold.}
\medskip

\textit{Proof.}  I$_{n+1}$:  First note that by III$_{n}$,
\begin{equation*}
|v_{n}|_{C^{2N+4}}\leq C\parallel v_{n}
\parallel_{2N+6}\leq CC_{3}\varepsilon<1
\end{equation*}
for small $\varepsilon$.  Therefore if $m\leq m_{*}-N-8$ we may
apply Theorem 5.1 to obtain a solution
$u_{n}\in\overline{H}_{0}^{m}(\Omega)$ of (5.5) which satisfies
the estimate (5.9).  When $m+2N+8-\rho\geq 1/2$ this may be
combined with V$_{n}$, VII$_{n}$, and $\rho\geq 2N+8$ to obtain
\begin{eqnarray*}
\parallel u_{n}\parallel_{m}&\leq& C_{m}(\parallel
f_{n}\parallel_{m+N}+\parallel v_{n}\parallel_{m+N+6}\parallel
f_{n}\parallel_{N+2})\\
&\leq&
C_{m}(C_{5}\varepsilon^{2}(1+\mu^{\rho-m})\mu_{n}^{m+N-\rho}+C_{3}C_{5}\varepsilon^{3}
(1+\mu^{\rho-N-2})\mu_{n}^{m+2N+8-\rho}\mu_{n}^{N+2-\rho})\\
&\leq&\varepsilon\mu_{n}^{m+N+2-\rho}
\end{eqnarray*}
for small $\varepsilon$.  When $m+2N+8-\rho\leq-1/2$, the estimate
$\parallel v_{n}\parallel_{m+N+6}\leq C_{3}\varepsilon$ placed in
the above calculation gives the desired result.\par
   II$_{n+1}$:  Since $w_{n+1}=\sum_{i=0}^{n}u_{i}$, we have
\begin{equation*}
\parallel w_{n+1}\parallel_{m}\leq\sum_{i=0}^{n}\parallel u_{i}
\parallel_{m}\leq\varepsilon\sum_{i=0}^{n}\mu_{i}^{m+N+2-\rho}.
\end{equation*}
Hence, if $m+N+2-\rho\leq-1/2$ then
\begin{equation*}
\parallel
w_{n+1}\parallel_{m}\leq\varepsilon\sum_{i=0}^{\infty}(\mu^{i})^{-1/2}
\leq\varepsilon\sum_{i=0}^{\infty}(2^{i})^{-1/2}:=C_{1}\varepsilon,
\end{equation*}
and if $m+N+2-\rho\geq 1/2$ then
\begin{equation*}
\parallel w_{n+1}\parallel_{m}\leq\varepsilon\mu_{n+1}^{m+N+2-\rho}
\sum_{i=0}^{n}(\frac{\mu_{i}}{\mu_{n+1}})^{m+N+2-\rho}
\leq\varepsilon\mu_{n+1}^{m+N+2-\rho}\sum_{i=0}^{\infty}(\mu^{-i})^{1/2}
\leq C_{1}\varepsilon\mu_{n+1}^{m+N+2-\rho}.
\end{equation*}\par
   III$_{n+1}$:  By the largeness assumption on $\rho$ we have $3N+8-\rho\leq-1/2$.
Therefore II$_{n+1}$ and V$_{n+1}$ (proven below) imply that
\begin{equation*}
\parallel w_{n+1}\parallel_{2N+6}\leq C_{1}\varepsilon
\text{ }\text{ }\text{ and }\text{ }
\parallel v_{n+1}\parallel_{2N+6}\leq
C_{3}\varepsilon.
\end{equation*}\par
   IV$_{n+1}$:  Since $\rho< m_{*}-6$ we have $(m_{*}-N-8)+N+2-\rho\geq 1/2$.
Therefore Lemma 5.1 and II$_{n+1}$ yield,
\begin{eqnarray*}
\parallel w_{n+1}-v_{n+1}\parallel_{m}&=&\parallel
(I-S_{n+1})w_{n+1}\parallel_{m}\\
&\leq& C_{m}\mu^{m-(m_{*}-N-8)}_{n+1}\parallel w_{n+1}\parallel_{m_{*}-N-8}\\
&\leq& C_{m}\mu^{m-(m_{*}-N-8)}_{n+1}C_{1}\varepsilon\mu^{(m_{*}-N-8)+N+2-\rho}_{n+1}\\
&:=& C_{2}\varepsilon\mu^{m+N+2-\rho}_{n+1}.
\end{eqnarray*}\par
   V$_{n+1}$:  From Lemma 5.1 and $\rho< m_{*}-6$ we have for all
$m\geq 0$,
\begin{equation*}
\parallel \!v_{n+1}\!\parallel_{m}=\parallel\!
S_{n+1}w_{n+1}\!\parallel_{m}\leq C_{m}
\begin{cases}
\parallel \!w_{n+1}\!\parallel_{\rho-N-3}
& \text{if }m+N+2-\rho\leq-1/2,\\
\mu_{n+1}^{m+N+1-\rho}\!\parallel\! w_{n+1}\!\parallel_{\rho-N-1}
& \text{if }m+N+2-\rho\geq 1/2.
\end{cases}
\end{equation*}
V$_{n+1}$ now follows from II$_{n+1}$.  Q.E.D.\medskip

   Write $e_{n}=e^{'}_{n}+e^{''}_{n}+e^{'''}_{n}$, where
\begin{eqnarray*}
e^{'}_{n}&=&(\mathcal{L}(w_{n})-\mathcal{L}(v_{n}))u_{n},\\
e^{''}_{n}&=&\varepsilon(I-S_{n}^{'})a^{22}(v_{n})L_{\theta_{n}}(v_{n})u_{n}
+\varepsilon\theta_{n}a^{22}(v_{n})\partial_{\xi_{n}^{1}}^{2}u_{n}\\
&
&+\varepsilon(a^{22}(v_{n}))^{-1}\Phi(v_{n})[\partial_{x^{1}}^{2}u_{n}-\partial_{x^{1}}(
\log
a^{22}(v_{n})\sqrt{|g|})\partial_{x^{1}}u_{n}],\\
e^{'''}_{n}&=&Q_{n}(w_{n},u_{n}).
\end{eqnarray*}

\textbf{Proposition 5.2.}  \textit{If the hypotheses of
Proposition 5.1 hold in addition to $n>0$, $\rho\geq 2N+12$, and
$0\leq m\leq \min(m_{*}-N-10,m_{*}-m_{0})$ then} VI$_{n+1}$
\textit{holds.}\medskip

\textit{Proof.}  We will estimate $e^{'}_{n}$, $e^{''}_{n}$, and
$e^{'''}_{n}$ separately.  According to (2.3) we may write
\begin{equation*}
(\mathcal{L}(w_{n})-\mathcal{L}(v_{n}))u_{n}
=\varepsilon\sum_{i,j}A_{ij}\partial_{x^{i}x^{j}}u_{n}
+\varepsilon\sum_{i}A_{i}\partial_{x^{i}}u_{n}.
\end{equation*}
Then Lemma 5.2 $(i)$ and $(iii)$, I$_{n+1}$, and IV$_{n}$ show
that
\begin{eqnarray*}
\parallel e_{n}^{'}\parallel_{m}&\leq& \varepsilon C_{m,1}
[(\sum_{i,j}\parallel A_{ij}\parallel_{m}+\sum_{i}\parallel
A_{i}\parallel_{m})\parallel
u_{n}\parallel_{4}\\
& &+(\sum_{i,j}\parallel A_{ij}\parallel_{2}+\sum_{i}\parallel
A_{i}\parallel_{2})\parallel u_{n}\parallel_{m+2}]\\
&\leq& \varepsilon C_{m,2}(\parallel
w_{n}-v_{n}\parallel_{m+2}\parallel u_{n}\parallel_{4}+\parallel
w_{n}-v_{n}\parallel_{4}\parallel u_{n}\parallel_{m+2})\\
&\leq&
C_{m,3}\varepsilon^{3}\mu_{n}^{m+N+4-\rho}\mu_{n}^{N+6-\rho}\\
&\leq&C_{m,3}\varepsilon^{3}\mu_{n}^{m-\rho}.
\end{eqnarray*}
Note that we have used $\rho\geq 2N+10$, as well as $m\leq
m_{*}-N-10$ which allows us to apply I$_{n+1}$ and IV$_{n}$.\par
   We now estimate $e^{''}_{n}$.  By Lemma 5.2 $(i)$ and $(iii)$,
I$_{n+1}$, V$_{n}$, and VIII$_{n}$,
\begin{eqnarray*}
\parallel\varepsilon\theta_{n}a^{22}(v_{n})\partial_{\xi_{n}^{1}}^{2}u_{n}
\parallel_{m}\!\!\!\!&\leq&\!\!\!\!\varepsilon\theta_{n}C_{m,4}(\parallel
a^{22}(v_{n})\parallel_{m}\parallel\partial_{\xi_{n}^{1}}^{2}u_{n}\parallel_{2}
+\parallel
a^{22}(v_{n})\parallel_{2}\parallel\partial_{\xi_{n}^{1}}^{2}u_{n}\parallel_{m})\\
&\leq&\!\!\!\!\varepsilon\theta_{n}C_{m,5}[(1+\parallel
v_{n}\parallel_{m+2})\parallel u_{n}\parallel_{4}+(1+\parallel
v_{n}\parallel_{4})\parallel u_{n}\parallel_{m+2}]\\
&\leq&\!\!\!\!\varepsilon^{2}\mu_{n}^{3-\rho}C_{m,6}[(1\!+C_{3}\varepsilon\mu_{n}^{m+N+4-\rho})\varepsilon
\mu_{n}^{N+4-\rho}\!+(1\!+C_{3}\varepsilon)\varepsilon\mu_{n}^{m+N+4-\rho}]\\
&\leq&\!\!\!\!C_{m,7}\varepsilon^{3}\mu_{n}^{m-\rho}
\end{eqnarray*}
if $\mu$ is large and $m+N+4-\rho\geq 1/2$.  If $m+N+4-\rho\leq
-1/2$ then we may use the estimate $\parallel
v_{n}\parallel_{m+2}\leq C_{3}\varepsilon$ to obtain the same
outcome.  Another application of Lemma 5.2 gives
\begin{eqnarray*}
\parallel\varepsilon(a^{22}(v_{n}))^{-1}\Phi(v_{n})\partial_{x^{1}}^{2}u_{n}\parallel_{m}&\leq&
\varepsilon C_{m,8}[\parallel\Phi(v_{n})\parallel_{m}(1+\parallel
v_{n}\parallel_{4})\parallel u_{n}\parallel_{4}\\
& &+\parallel\Phi(v_{n})\parallel_{2}(1+\parallel
v_{n}\parallel_{m+2})\parallel u_{n}\parallel_{4}\\
& &+\parallel\Phi(v_{n})\parallel_{2}(1+\parallel
v_{n}\parallel_{4})\parallel u_{n}\parallel_{m+2}]\\
&\leq&C_{m,9}\varepsilon^{3}\mu_{n}^{m-\rho}
\end{eqnarray*}
after noting that $m\leq\min(m_{*}-N-10,m_{*}-m_{0})$ is required
for VIII$_{n}$ to be valid, and similar methods yield
\begin{eqnarray*}
\parallel\varepsilon(a^{22}(v_{n}))^{-1}\Phi(v_{n})\partial_{x^{1}}(
\log a^{22}(v_{n})\sqrt{|g|})\partial_{x^{1}}u_{n}\parallel_{m}
\leq C_{m,10}\varepsilon^{3}\mu_{n}^{m-\rho}.
\end{eqnarray*}
Moreover, if $l=\rho+2\leq m_{*}-2$ and $n>0$ then we may apply
Lemma 5.1 and recall that $\mu=\varepsilon^{-\frac{1}{2\rho}}$ to
obtain
\begin{eqnarray*}
&
&\parallel\varepsilon(I-S_{n}^{'})a^{22}(v_{n})L_{\theta_{n}}(v_{n})u_{n}\parallel_{m}\\
&\leq&\!\!\!\!\varepsilon
C_{m,11}[\parallel(I-S_{n}^{'})a^{22}(v_{n})\parallel_{m}\parallel
f_{n}\parallel_{2}+\parallel(I-S_{n}^{'})a^{22}(v_{n})\parallel_{2}\parallel f_{n}\parallel_{2}]\\
&\leq&\!\!\!\!\varepsilon C_{m,12}[\mu_{n}^{m-l}(1+\parallel
v_{n}\parallel_{l+2})\varepsilon^{2}(1\!+\mu^{\rho-2})\mu_{n}^{2-\rho}+\mu_{n}^{2-l}(1+\parallel
v_{n}\parallel_{l+2})\varepsilon^{2}(1\!+\mu^{\rho-m})\mu_{n}^{m-\rho}]\\
&\leq&\!\!\!\! C_{m,13}\varepsilon^{3}\mu_{n}^{m-\rho}.
\end{eqnarray*}
Therefore
\begin{equation*}
\parallel e^{''}_{n}\parallel\leq
C_{m,14}\varepsilon^{3}\mu^{m-\rho}_{n}.
\end{equation*}\par
   We now estimate $e^{'''}_{n}$.  We have
\begin{equation*}
e^{'''}_{n}=Q_{n}(w_{n},u_{n})=\int_{0}^{1}(1-t)\frac{\partial^{2}}
{\partial t^{2}}\Phi(w_{n}+tu_{n})dt.
\end{equation*}
Apply Lemma 5.2 $(i)$ and $(ii)$, as well as the Sobolev Lemma to
obtain
\begin{eqnarray*}
\parallel e^{'''}_{n}\parallel_{m}&\leq&\int_{0}^{1}
\sum_{|\alpha|,|\beta|\leq
2}\parallel\nabla_{\overline{\alpha}\overline{\beta}}
\Phi(w_{n}+tu_{n})\partial^{\alpha}u_{n}\partial^{\beta}u_{n}\parallel
_{m}dt\\
&\leq&\int_{0}^{1}C_{m,15}(\parallel\nabla^{2}\Phi(w_{n}+tu_{n})\parallel_{2}
\parallel u_{n}\parallel_{4}\parallel u_{n}\parallel_{m+2}\\
& &+\parallel\nabla^{2}\Phi(w_{n}+tu_{n})\parallel_{m}
\parallel u_{n}\parallel_{4}^{2})dt,
\end{eqnarray*}
where $\overline{\alpha}=\partial^{\alpha}(w_{n}+tu_{n})$ and
$\overline{\beta}=\partial^{\beta}(w_{n}+tu_{n})$.  The notation
$\nabla^{2}\Phi$ represents the collection of second partial
derivatives with respect to the variables $\overline{\alpha}$ and
$\overline{\beta}$.  Furthermore it is easy to see that
$|\nabla^{2}\Phi(0)|=O(\varepsilon)$.  Therefore using Lemma 5.2
$(iii)$, I$_{n+1}$, and II$_{n}$, we have
\begin{eqnarray*}
\parallel e^{'''}_{n}\parallel_{m}&\leq& C_{m,16}[(\varepsilon+\parallel
w_{n}\parallel_{4}+\parallel u_{n}\parallel_{4})\parallel u_{n}
\parallel_{4}\parallel u_{n}\parallel_{m+2}\\
& &+(\varepsilon+\parallel w_{n}
\parallel_{m+2}+\parallel u_{n}\parallel_{m+2})\parallel u_{n}
\parallel^{2}_{4}]\\
&\leq&C_{m,17}\varepsilon^{3}\mu^{m-\rho}_{n}
\end{eqnarray*}
if $\rho\geq 2N+12$.  Combining the estimates of $e^{'}_{n}$,
$e^{''}_{n}$, and $e^{'''}_{n}$ yields the desired result.
Q.E.D.\medskip
   According to the above proposition if
$\rho+1\leq m_{*}-m_{0}$ (in addition to the other required
restrictions) then $E_{n}\in H^{\rho+1}(\Omega)$ and the following
estimate holds, which will be utilized in the next proposition:
\begin{eqnarray}
\parallel E_{n}\parallel_{\rho+1}\leq\sum_{i=0}^{n-1}\parallel
e_{i}\parallel_{\rho+1}\!\!\!&\leq&\!\!\!
C_{4}\varepsilon^{3}\sum_{i=0}^{n-1}
\mu_{i}\\
&\leq&\!\!\!
C_{4}\varepsilon^{3}(\sum_{i=0}^{\infty}\mu^{-1}_{i})\mu^{n} \leq
C_{4}\varepsilon^{3}(\sum_{i=0}^{\infty}2^{-i})\mu_{n}.\nonumber
\end{eqnarray}

\textbf{Proposition 5.3.}  \textit{If the hypotheses of
Proposition 5.2 hold and $\rho+1\leq m_{*}-m_{0}$, then}
VII$_{n+1}$ \textit{holds for all $0\leq m<\infty$.}\medskip

\textit{Proof.}  By (5.7) as well as Lemma 5.2 $(i)$ and $(iii)$,
\begin{eqnarray}
\parallel f_{n+1}\parallel_{m}&\leq&\varepsilon^{-1}
C_{m}(\parallel
S_{n}E_{n}-S_{n+1}E_{n+1}+(S_{n}-S_{n+1})\Phi(w_{0})\parallel_{m}\\
& &+\parallel v_{n+1}\parallel_{m+2}\parallel
S_{n}E_{n}-S_{n+1}E_{n+1}+(S_{n}-S_{n+1})\Phi(w_{0})\parallel_{2}).\nonumber
\end{eqnarray}
Next observe that (2.2) together with (5.14) below yields
\begin{equation}
\parallel\Phi(w_{0})\parallel_{\rho+1}\leq C(\varepsilon^{3}+
\parallel w_{0}\parallel_{\rho+3})\leq C\varepsilon^{3}.
\end{equation}
Then (5.10) implies that for all $m\geq \rho+1$,
\begin{eqnarray}
& &\!\!\!\parallel
S_{n}E_{n}-S_{n+1}E_{n+1}+(S_{n}-S_{n+1})\Phi(w_{0})\parallel_{m}\\
&\leq& \!\!\!C_{m}(\mu_{n}^{m-\rho-1}\parallel
E_{n}\parallel_{\rho+1}+\mu_{n+1}^{m-\rho-1}\parallel
E_{n+1}\parallel_{\rho+1}
+(\mu_{n}^{m-\rho-1}+\mu_{n+1}^{m-\rho-1})\parallel\Phi(w_{0})\parallel_{\rho+1})\nonumber\\
&\leq&\!\!\!
C_{m}\varepsilon^{3}(1+\mu^{\rho-m})\mu_{n+1}^{m-\rho}.\nonumber
\end{eqnarray}
If $m<\rho+1$, then applying similar methods along with VI$_{n+1}$
to
\begin{eqnarray*}
& &\parallel
S_{n}E_{n}-S_{n+1}E_{n+1}+(S_{n}-S_{n+1})\Phi(w_{0})\parallel_{m}\\
&\leq& \parallel (I-S_{n})E_{n}\parallel_{m}+\parallel
(I-S_{n+1})E_{n}\parallel_{m}+\parallel
S_{n+1}e_{n}\parallel_{m}\\
& &+\parallel (I-S_{n})\Phi(w_{0})\parallel_{m}+\parallel
(I-S_{n+1})\Phi(w_{0})\parallel_{m},
\end{eqnarray*}
produces the same estimate found in (5.13).  Therefore plugging
into (5.11) produces
\begin{eqnarray*}
\parallel f_{n+1}\parallel_{m}&\leq&
C_{m}[\varepsilon^{2}(1+\mu^{\rho-m})\mu^{m-\rho}_{n+1}
+\varepsilon^{3}(1+\mu^{\rho-2})\mu^{m+N+6-2\rho}_{n+1}]\\
&\leq& C_{m}\varepsilon^{2}(1+\mu^{\rho-m})\mu^{m-\rho}_{n+1},
\end{eqnarray*}
if $m+N+4-\rho\geq 1/2$.  If $m+N+4-\rho\leq-1/2$ and $m\geq 2$,
then using $\parallel v_{n+1}\parallel_{m+2}\leq C_{3}\varepsilon$
in the estimate above gives the desired result. Moreover if $0\leq
m<2$, then in place of (5.11) we use the estimate
\begin{equation*}
\parallel f_{n+1}\parallel_{m}\leq\varepsilon^{-1}C_{m}\parallel
S_{n}E_{n}-S_{n+1}E_{n+1}+(S_{n}-S_{n+1})\Phi(w_{0})\parallel_{m}
\end{equation*}
combined with the above method to obtain the desired result.
Lastly if $m+N+4-\rho=0$, then replace $\parallel
v_{n+1}\parallel_{m+2}$ in (5.11) by $\parallel
v_{n+1}\parallel_{m+3}$ and follow the above method.
Q.E.D.\medskip

\textbf{Proposition 5.4.}  \textit{If the hypotheses of
Proposition 5.3 hold and $\rho+1=\min(m_{*}-N-10,m_{*}-m_{0})$,
then} VIII$_{n+1}$ \textit{holds for $0\leq m\leq\min(m_{*}-N-10,
m_{*}-m_{0})$.}\medskip

\textit{Proof.}  By (5.8) and VI$_{n+1}$ and $m\leq\rho+1$, we
have
\begin{eqnarray*}
\parallel\Phi(w_{n+1})\parallel_{m}
&\leq&\parallel(I-S_{n})\Phi(w_{0})\parallel_{m}+\parallel(I-S_{n})E_{n}\parallel
_{m}+\parallel e_{n}\parallel_{m}\\
&\leq&
C_{m}(\mu^{m-\rho-1}_{n}\parallel\Phi(w_{0})\parallel_{\rho+1}+
\mu^{m-\rho-1}_{n}\parallel
E_{n}\parallel_{\rho+1}+C_{4}\varepsilon^{3}\mu^{m-\rho}_{n}).
\end{eqnarray*}
Applying the estimate (5.10) along with (5.12) and
$\varepsilon\mu^{\rho-m}\leq\varepsilon^{1/2}$ produces
\begin{equation*}
\parallel\Phi(w_{n+1})\parallel_{m}\leq
C_{m}\varepsilon^{3}
\mu^{\rho-m}\mu_{n+1}^{m-\rho}\leq\frac{1}{3}\varepsilon\mu_{n+1}^{m-\rho},
\end{equation*}
if $\varepsilon$ is sufficiently small.  Lastly a similar estimate
may be obtained for $\Phi(v_{n+1})$ by writing
\begin{eqnarray*}
\parallel\Phi(v_{n+1})\parallel_{m}&\leq&\parallel\Phi(w_{n+1})\parallel_{m}
+\parallel\Phi(v_{n+1})-\Phi(w_{n+1})\parallel_{m}\\
&\leq&\frac{1}{3}\varepsilon\mu_{n+1}^{m-\rho}+\varepsilon\parallel
v_{n+1}-w_{n+1}\parallel_{m+2}\\
&\leq&(\frac{1}{3}+C_{2}\varepsilon^{2})\mu_{n+1}^{m+N+4-\rho}.
\end{eqnarray*}
Q.E.D.\medskip

   To complete the proof by induction we will now prove the case
$n=0$.  Here we will assume that the initial data are
appropriately small:
\begin{equation}
\parallel\phi_{\varrho}^{-}\parallel_{m_{*}-m_{0}+1,\overline{\partial}\Omega}+
\parallel\psi_{\varrho}^{-}\parallel_{m_{*}-m_{0},\overline{\partial}\Omega}\leq
C\varepsilon^{l},\text{ }\text{ }\text{ }\text{ }l\geq 3.
\end{equation}
Then according to (5.4), II$_{0}$, III$_{0}$, IV$_{0}$, and
V$_{0}$ are trivial as long as $\varepsilon$ is small enough.
Furthermore applying (5.12) and again taking $\varepsilon$ to be
sufficiently small yields VII$_{0}$ and VIII$_{0}$.  In addition
by the proof of Proposition 5.1 we obtain the following stronger
version of I$_{1}$,
\begin{equation*}
\parallel u_{0}\parallel_{m}\leq C_{0}\varepsilon^{2},\text{ }\text{ }\text{ }\text{ }
m\leq m_{*}-N-8.
\end{equation*}
Now the proof of Proposition 5.2 may be appropriately modified to
show that VI$_{1}$ is valid.  This completes the proof by
induction.\par
   In view of the hypotheses of Propositions 5.1-5.4, we will
choose
\begin{equation*}
\rho=\min(m_{*}-N-10,m_{*}-m_{0})-1.
\end{equation*}
Since $\rho\geq 3N+9$ we must then have
\begin{equation*}
m_{*}\geq\max(3N+16, 3N+m_{0}+10).
\end{equation*}
The following corollary yields a solution of (5.2) with
$\alpha_{0}=m_{*}-m_{0}-N-6$.\medskip

\textbf{Corollary 5.1.}  \textit{If $m_{0}\geq N+10$ then under
the above assumptions $w_{n}\rightarrow w$\textit{ in
}$\overline{H}_{0}^{m_{*}-m_{0}-N-4}(\Omega)$.  Furthermore
$\Phi(w_{n})\rightarrow 0$\textit{ in }$C^{0}(\Omega)$.}\medskip

\textit{Proof.}  When $m_{0}\geq N+10$ we have
$\rho-1=m_{*}-m_{0}-2$.  Then for $m+N+2\leq \rho-1$ and $i>j$,
I$_{n}$ implies that
\begin{equation*}
\parallel w_{i}-w_{j}\parallel_{m}\leq
\sum_{n=j}^{i-1}\parallel u_{n}\parallel_{m}\leq
\varepsilon\sum_{n=j}^{i-1}\mu_{n}^{m+N+2-\rho}\leq\varepsilon
\sum_{n=j}^{i-1}\mu^{-n}.
\end{equation*}
Hence, $\{w_{n}\}_{n=0}^{\infty}$ is Cauchy in
$\overline{H}_{0}^{m}(\Omega)$ for all $m\leq
m_{*}-m_{0}-N-4$.\par
  Lastly by the Sobolev Lemma and VIII$_{n}$,
\begin{equation*}
|\Phi(w_{n})|_{C^{0}(\Omega)}\leq
C\parallel\Phi(w_{n})\parallel_{2}\leq
\varepsilon\mu^{N+6-\rho}_{n}.
\end{equation*}
The desired conclusion follows since $\rho>N+6$. Q.E.D.\bigskip

\textsc{Elliptic Regions}\bigskip

  Here we shall set up the iteration procedure for problem (5.1).
For convenience we will denote the domain $\Omega_{\kappa}^{+}$ by
$\Omega$.  Set $w_{0}=0$ and suppose that functions
$w_{1},\ldots,w_{n}$ have been defined on $\Omega$.  If
$S_{i}=S_{\mu^{i}}$ are smoothing operators given by Lemma 5.1,
then we put $v_{i}=S_{i}w_{i}$, $0\leq i\leq n$, and define
$w_{n+1}=w_{n}+u_{n}$ where $u_{n}$ is the solution of
\begin{equation}
L(v_{n})u_{n}=f_{n}\text{ }\text{ }\text{ in }\text{ }\text{
}\Omega,\text{ }\text{ }\text{ }\text{
}u_{n}|_{\overline{\partial}\Omega}=0,
\end{equation}
given by Theorem 5.2 below, $L(v_{n})$ is the operator of Lemma
2.1, and $f_{n}$ will also be specified below.  Let
$Q_{n}(w_{n},u_{n})$ again denote the quadratic error and
$\mathcal{L}(w_{n})$ the linearization of (5.1), then according to
(5.3) (with $\theta_{n}=0$) we have
\begin{eqnarray*}
\Phi(w_{n+1})&=&\Phi(w_{n})+\mathcal{L}(w_{n})u_{n}+Q_{n}(w_{n},u_{n})\\
&=&\Phi(w_{n})+\varepsilon
S_{n}^{'}a^{22}(v_{n})L(v_{n})u_{n}+e_{n},
\end{eqnarray*}
with
\begin{eqnarray*}
e_{n}&=&(\mathcal{L}(w_{n})-\mathcal{L}(v_{n}))u_{n}+\varepsilon(I-S_{n}^{'})
a^{22}(v_{n})L(v_{n})u_{n}
+Q_{n}(w_{n},u_{n})\\
&
&+\varepsilon(a^{22}(v_{n}))^{-1}\Phi(v_{n})[\partial_{x^{1}}^{2}u_{n}-\partial_{x^{1}}(
\log a^{22}(v_{n})\sqrt{|g|})\partial_{x^{1}}u_{n}].
\end{eqnarray*}
Lastly we set $E_{0}=0$, $E_{n}=\sum_{i=0}^{n-1}e_{i}$, and define
$f_{n}$ according to (5.7).\par
  It is clear that similar arguments as those used for the
hyperbolic regions will show that $\{w_{n}\}_{n=0}^{\infty}$
converges to a solution of (5.1) if a Moser estimate (like that
found in Theorem 5.1) holds for the solution of (5.15).  In order
to establish such an estimate using the theory of section $\S 3$,
we need to extend the coefficients of $L(v_{n})$ outside of
$\Omega$ and cut them off.  For this purpose we will use the
following extension lemma.\medskip

\textbf{Lemma 5.3 [18].}  \textit{Let $X$ be a bounded convex domain
in $\mathbb{R}^{2}$, with Lipschitz smooth boundary. Then there
exists a linear operator $E_{X}:L^{2}(X)\rightarrow
L^{2}(\mathbb{R}^{2})$ such that:}

$i)$  $E_{X}(u)|_{X}=u$,

$ii)$  $E_{X}:H^{m}(X)\rightarrow H^{m}(\mathbb{R}^{2})$
\textit{continuously for each $m\in\mathbb{Z}_{\geq 0}$}.\medskip

\textbf{Theorem 5.2.}  \textit{Suppose that $g\in C^{m_{*}}$. If
$m\leq \frac{1}{3}(m_{*}-8)$, $|v_{n}|_{C^{6}}<1$, and
$\delta=\delta(m)$, $\varepsilon=\varepsilon(m,\delta)$ are
sufficiently small, then there exists a solution
$u_{n}\in\overline{H}^{m}_{0}(\Omega)$ of (5.15) which satisfies
the estimate}
\begin{equation*}
\parallel u_{n}\parallel_{m}\leq\delta^{-1}C_{m}(\parallel
f_{n}\parallel_{m+2+\gamma}+\sum_{i+j+l\leq
m+23+\gamma}(1+\parallel v_{n}\parallel_{i})\parallel
v_{n}\parallel_{j}\parallel f_{n}\parallel_{l})
\end{equation*}
\textit{for some constant $C_{m}$ independent of $\delta$ and
$\varepsilon$ and where $2m<\gamma<m_{*}-m-6$.}\medskip

\textit{Proof.}  This will follow from Theorem 3.2.  However we
must first change to the coordinates $\xi_{n}^{i}(x,y)$ of Lemma
2.1, and then change to polar coordinates so that
\begin{equation*}
\Omega=\{(r,\theta)\mid 0<r<\sigma,0<\theta<\delta\}.
\end{equation*}
Set
\begin{equation*}
\Omega^{1}=\{(r,\theta)\mid 0<r<\sigma+1,0<\theta<\delta\}
\end{equation*}
and let $\varphi=\varphi(r)$ be a smooth nonnegative cut-off
function with $\varphi(r)\equiv 1$ for $0<r<\sigma$, and
$\varphi(r)\equiv 0$ for $\sigma+1<r$.  If we cut-off the
coefficients of $L(E_{\Omega}v_{n})$ as in (3.1), we may use
Theorem 3.2 to solve
\begin{equation*}
L(E_{\Omega}v_{n})u_{n}=E_{\Omega}f_{n}\text{ }\text{ }\text{ in
}\text{ }\text{ }\Omega^{1},\text{ }\text{ }\text{ }\text{
}u_{n}|_{\overline{\partial}\Omega^{1}}=0,
\end{equation*}
with
\begin{equation*}
\parallel u_{n}\parallel_{(m,\gamma),\Omega^{1}}^{'}\leq C_{m}(\parallel
E_{\Omega}f_{n}\parallel^{'}_{m+2+\gamma,\Omega^{1}}+\parallel
E_{\Omega}v_{n}\parallel^{'}_{m+6,\Omega^{1}}
\parallel E_{\Omega}f_{n}\parallel^{'}_{5+\gamma,\Omega^{1}})
\end{equation*}
for $m\leq m_{*}-6$ where $\gamma> 2m$ and
$\parallel\cdot\parallel^{'}$ indicates that the norm is with
respect to these polar coordinates. By Lemma 5.3
\begin{equation*}
\parallel E_{\Omega}f_{n}\parallel^{'}_{m+2+\gamma,\Omega^{1}}
\leq C_{m}\parallel f_{n}\parallel^{'}_{m+2+\gamma,\Omega},\text{
}\text{ }\text{ }\text{ }\parallel
E_{\Omega}v_{n}\parallel^{'}_{m+6,\Omega^{1}} \leq C_{m}\parallel
v_{n}\parallel^{'}_{m+6,\Omega}.
\end{equation*}
Therefore with the help of (2.4) and Lemma 5.2 it follows that
\begin{equation*}
\parallel u_{n}\parallel_{m,\Omega}\leq\delta^{-1}C_{m}(\parallel
f_{n}\parallel_{m+2+\gamma,\Omega}+\sum_{i+j+l\leq
m+23+\gamma}(1+\parallel v_{n}\parallel_{i,\Omega})\parallel
v_{n}\parallel_{j,\Omega}
\parallel f_{n}\parallel_{l,\Omega}).
\end{equation*}
The result is now obtained by noting that
$\max(m+2+\gamma,5+\gamma)\leq m_{*}-2$ is required to apply
(2.4). Q.E.D.\medskip

  We may now apply arguments similar to those in the hyperbolic
regions to obtain a solution of (5.1).  More precisely, the proofs
of Propositions 5.1-5.4 yield the following restrictions on
$\rho$, $\gamma$, and $m_{*}$ in the elliptic regions:
\begin{equation*}
\rho\geq 2\gamma+54,\text{ }\text{ }\text{ }\text{
}\rho+1=\frac{1}{3}(m_{*}-14).
\end{equation*}
Choosing the largest possible value for $\gamma$ and noting that
the hypothesis of Theorem 5.2 requires $2m<\gamma$, implies that
we must have $m\leq\frac{1}{12}(m_{*}-185)$.  The following
corollary produces a solution of (5.1) for
$\alpha_{0}=\frac{1}{12}m_{*}-18$.\medskip

\textbf{Corollary 5.2.}  \textit{If $m_{*}\geq 192$ then
$w_{n}\rightarrow w$ in
$\overline{H}_{0}^{\frac{1}{12}m_{*}-16}(\Omega)$ with}
\begin{equation*}
\parallel w\parallel_{\frac{1}{12}m_{*}-16}\leq
C\varepsilon^{3}.
\end{equation*}
\textit{Furthermore $\Phi(w_{n})\rightarrow 0$ in
$C^{0}(\Omega)$.}\medskip

\textit{Proof.}  The same arguments used for Corollary 5.1 apply.
Moreover, we use the analogue of II$_{n}$ to obtain the estimate
for $w$. Q.E.D.
\bigskip

\textsc{Proof of Theorem 1.1}\bigskip

  Here we shall construct a solution of (2.1) in a full
neighborhood of the origin.  First consider the case in which
there are exactly two elliptic regions, each bordering two
hyperbolic regions. Then on each elliptic region
$\Omega_{\kappa}^{+}$ let
$w_{\kappa}^{+}\in\overline{H}_{0}^{\frac{1}{12}m_{*}-16}(\Omega_{\kappa}^{+})$
be the solution of (5.1) given by Corollary 5.2.  On each boundary
of the hyperbolic regions
$\overline{\partial}\Omega_{\varrho}^{-}$ set
$\phi_{\varrho}^{-}=0$,
$\psi_{\varrho}^{-}=-\partial_{\nu}w_{\kappa(\varrho)}^{+}|_{\overline{\partial}
\Omega_{\kappa(\varrho)}^{+}}\in\overline{H}_{0}^{\frac{1}{12}m_{*}-18}
(\overline{\partial}\Omega_{\kappa(\varrho)})$ where
$\Omega_{\kappa(\varrho)}^{+}$ is the bordering elliptic region,
and note that (5.14) is valid with $m_{0}=\frac{11}{12}m_{*}+18$.
Then Corollary 5.1 yields a solution
$w_{\varrho}^{-}\in\overline{H}^{m_{*}-m_{0}-N-4}_{0}(\Omega_{\varrho}^{-})$
of (5.2).  Under the hypotheses of Corollaries 5.1 and 5.2 we
require $m_{*}\geq \max(192,3N+m_{0}+10)$ or rather $m_{*}\geq
36(N+10)$.\par
  Suppose that $\Omega_{\kappa}^{+}$ borders on
$\Omega_{\varrho}^{-}$.  Then since the common boundary curve
$\Upsilon$ is noncharacteristic for (2.1) (according to our
original choice of approximate solution $z_{0}$), the functions
$w_{\kappa}^{+}$ and $w_{\varrho}^{-}$ agree along with their
derivatives up to and including order $\frac{1}{12}m_{*}-N-24$
along $\Upsilon$.  It follows that the individual solutions
$\{w_{\kappa}^{+}\}_{\kappa=1}^{\kappa_{0}}$ and
$\{w_{\varrho}^{+}\}_{\varrho=1}^{\varrho_{0}}$ combine to form a
$C^{\frac{1}{12}m_{*}-N-24}$ solution of (2.1) on some
neighborhood of the origin.\par
  Now consider the general case in which elliptic and hyperbolic
regions are allowed to border regions of the same type.  If an
elliptic region borders another elliptic region, they may be
combined to form a single elliptic region which contains a single
curve of degeneracy on the interior.  By appropriately
regularizing the linearized equation in this combined region to
eliminate the interior degeneracy, we may apply the theory of
section $\S 3$ to obtain Theorem 5.2, and hence a solution of
(5.1) in this combined region.  Therefore we may assume that each
elliptic region is bordered by hyperbolic regions (unless no
hyperbolic regions are present).  On the other hand, if two
hyperbolic regions share a common boundary, for instance
$\Omega_{1}^{-}$ and $\Omega_{3}^{-}$, then Cauchy data may be
prescribed appropriately on the portion of
$\partial\Omega_{1}^{-}$ which is shared with
$\partial\Omega_{3}^{-}$, so that the solution on both regions may
be glued together. Moreover, Cauchy data may be arbitrarily
prescribed on the portion of $\partial\Omega_{3}^{-}$ emanating
from the origin and which is not shared with
$\partial\Omega_{1}^{-}$. It follows that in the general case, the
solutions of the elliptic and hyperbolic regions may be patched
together in the usual way.

\begin{center}
\textbf{6.  Appendix A}
\end{center}\setcounter{equation}{0}
\setcounter{section}{6}

  The purpose of this appendix is to show existence for the ODE
occurring in the proof of Theorem 3.1:
\begin{eqnarray}
&
&\sum_{s=0}^{m}\lambda^{-s}(-1)^{s}\partial_{r}^{s}(a_{\lambda,\gamma-2(s-1)}
\partial_{r}^{s}\zeta)=v,\\
& &\zeta(r,0)=\zeta(r,\delta)=0,\text{ }\text{ }\text{
}\partial_{r}^{s}\zeta(\sigma,\theta)=0,\text{ }\text{ }\text{
}0\leq s\leq m-1,\nonumber\\
&
&\int_{r=r_{0}}(\partial_{r}^{s}\partial_{\theta}^{l}\zeta)^{2}\leq
r_{0}^{\gamma-2s}C,\text{ }\text{ }\text{ }0\leq s\leq 2m-1,\text{
}\text{ }\text{ }0\leq l<\infty,\nonumber
\end{eqnarray}
where $v\in\widehat{C}^{\infty}(\Omega)$, $\zeta\in
H^{(m,\infty,\gamma+2)}(\Omega)\cap C^{\infty}(\Omega)$, $r_{0}$
is sufficiently small, and all other definitions/notation may be
found in section $\S 3$.\par
  First note that $\eta\mapsto(\eta,(\lambda\theta^{2}-1)^{-1}v)$
is a bounded linear functional on $H^{(m,0,\gamma+2)}(\Omega)$,
and thus by the Riesz representation theorem there exists a unique
$\zeta\in H^{(m,0,\gamma+2)}(\Omega)$ such that
\begin{equation*}
(\eta,\zeta)_{(m,0,\gamma+2)}=(\eta,(\lambda\theta^{2}-1)^{-1}v)\text{
}\text{ }\text{ all }\text{ }\text{ }\eta\in
H^{(m,0,\gamma+2)}(\Omega),
\end{equation*}
where $(\cdot,\cdot)_{(m,0,\gamma+2)}$ denotes the inner product
on $H^{(m,0,\gamma+2)}(\Omega)$.  It follows that $\zeta$ is a
weak solution of (6.1), and according to the basic regularity
theory for ODEs we have $\zeta\in C^{\infty}(\Omega)$. Furthermore
the desired boundary behavior of the solution at $\theta=0,\delta$
arises from the requirement that $v(r,0)=v(r,\delta)=0$, and the
vanishing at $r=\sigma$ is a result of the trace theorem for
Sobolev spaces.\par
  Lastly we observe that since $v$ vanishes in a neighborhood of
$r=0$, the solution $\zeta$ satisfies a version of the so called
Euler differential equation in this domain.  All solutions of this
equation may be written down explicitly.  In particular, for $r$
sufficiently small $\zeta$ must be a linear combination of $2m$
functions of the form: $r^{\alpha}(\log r)^{\beta}$ where
$\alpha\in\mathbb{C}$ and $\beta\in\mathbb{Z}_{\geq 0}$.  However
according to Lemma 3.1
\begin{equation}
\int_{r=r_{0}}(\partial_{r}^{s}\partial_{\theta}^{l}\zeta)^{2}\leq
r_{0}^{\gamma-2s}C\parallel\zeta\parallel_{(m,l+1,\gamma+2)}^{2},\text{
}\text{ }\text{ }s\leq m-1,\text{ }\text{ }\text{ }0\leq l<\infty.
\end{equation}
Therefore each term in the linear combination must satisfy
\begin{equation}
r^{\alpha}(\log r)^{\beta}=O(r^{\gamma/2})\text{ }\text{ }\text{
as }\text{ }\text{ }r\rightarrow 0.
\end{equation}
The desired boundary behavior at $r=0$ now follows from (6.2) and
(6.3).

\begin{center}
\textbf{7.  Appendix B}
\end{center}\setcounter{equation}{0}
\setcounter{section}{7}

  The purpose of this section is to construct the smoothing
operators $S_{\mu}$ of Lemma 5.1.  The construction will differ
from the standard one for $S_{\mu}^{'}$ (see [16]), in that the
smoothed functions are required to vanish identically at the
origin.  This of course is only possible if the function being
smoothed already vanishes in an appropriate sense at the
origin.\par
  We first construct smoothing operators on the plane, and will later
restrict them back to the bounded domain $\Omega$.  Fix
$\widehat{\chi}\in C^{\infty}_{c}(\mathbb{R}^{2})$ such that
$\widehat{\chi}\equiv 1$ on some neighborhood of the origin, and
let
\begin{equation*}
\chi(x)=\int_{\mathbb{R}^{2}}e^{2\pi i\xi\cdot
x}\widehat{\chi}(\xi)d\xi
\end{equation*}
be its inverse Fourier Transform.  Then $\chi$ is a Schwartz
function and satisfies
\begin{equation*}
\int_{\mathbb{R}^{2}}\chi(x)dx=1,\text{ }\text{ }\text{ }\text{
}\int_{\mathbb{R}^{2}}x^{\alpha}\chi(x)dx=0,\text{ }\text{ }\text{
}\text{ }|\alpha|>0.
\end{equation*}
Furthermore let $\eta\in C^{\infty}(\mathbb{R}^{2})$ be a radial
function vanishing to all orders at the origin, and satisfying
\begin{equation*}
\eta(x)=\begin{cases}
1 & \text{if $|x|> 1$},\\
0 & \text{if $|x|<\frac{1}{2}$}.
\end{cases}
\end{equation*}
For $\mu\geq 1$ we will write $\eta_{\mu}(x)=\eta(\mu x)$,
$\chi_{\mu}(x)=\mu^{2}\chi(\mu x)$, and define smoothing operators
$\overline{S}_{\mu}:L^{2}(\mathbb{R}^{2})\rightarrow\overline{H}_{0}^{\infty}(\mathbb{R}^{2})$
by
\begin{equation*}
(\overline{S}_{\mu}u)(x)=\eta_{\mu}(x)(\chi_{\mu}*u)(x)=
\mu^{2}\eta(\mu x)\int_{\mathbb{R}^{2}}\chi(\mu(x-y))u(y)dy.
\end{equation*}
Here the space $\overline{H}_{0}^{l}(\mathbb{R}^{2})$ is the
completion of $\overline{C}_{c}^{\infty}(\mathbb{R}^{2})$ in the
Sobolev norm $\parallel\cdot\parallel_{l}$, where
$\overline{C}_{c}^{\infty}(\mathbb{R}^{2})$ denotes all
$C_{c}^{\infty}(\mathbb{R}^{2})$ functions vanishing in a
neighborhood of the origin.\par
  We now proceed to show statements $(i)$, $(ii)$, and $(iii)$ of
Lemma 5.1 with respect to $\overline{S}_{\mu}$.  Note that it is
sufficient to prove these for
$u\in\overline{C}_{c}^{\infty}(\mathbb{R}^{2})$ as this space of
functions is dense in $\overline{H}^{l}_{0}(\mathbb{R}^{2})$. We
begin with $(ii)$.  Let
\begin{equation*}
u(y)=\sum_{|\alpha|<|\rho|}\frac{1}{|\alpha|!}\partial^{\alpha}u(x)
(y-x)^{\alpha}+\frac{1}{(|\rho|-1)!}\sum_{|\alpha|=|\rho|}
(y-x)^{\alpha}\int_{0}^{1}(1-t)^{|\rho|-1}\partial^{\alpha}u(x+t(y-x))dt
\end{equation*}
be a Taylor expansion of $u$ with integral remainder.  Then
according to the properties of $\chi$,
\begin{eqnarray}
& &\!\!\!(\chi_{\mu}*u)(x)\\
&=&\!\!\!u(x)+
\frac{\mu^{2}}{(|\rho|-1)!}\sum_{|\alpha|=|\rho|}\int_{\mathbb{R}^{2}}\int_{0}^{1}\!
\chi(\mu(y\!-\!x))(1\!-\!t)^{|\rho|-1}(y\!-\!x)^{\alpha}\partial^{\alpha}u(x\!+t(y\!-\!x))dtdy.\nonumber
\end{eqnarray}
Suppose that $l\leq |\sigma|\leq m$, and notice that
\begin{equation*}
\parallel\partial^{\sigma}\overline{S}_{\mu}u\parallel\leq\sum_{\beta+\gamma=\sigma}
\parallel\partial^{\beta}\eta_{\mu}\partial^{\gamma}(\chi_{\mu}*u)\parallel.
\end{equation*}
If $|\rho|=l-|\gamma|> 0$ we may apply (7.1) to find
\begin{eqnarray*}
\parallel\partial^{\beta}\eta_{\mu}\partial^{\gamma}(\chi_{\mu}*u)\parallel^{2}&=&
\mu^{2|\beta|}\int_{\mathbb{R}^{2}}(\partial^{\beta}\eta)^{2}(\mu
x)
(\chi_{\mu}*\partial^{\gamma}u)^{2}(x)dx\\
&\leq&
C_{1}\mu^{2|\beta|}\int_{\mathbb{R}^{2}}(\partial^{\beta}\eta)^{2}(\mu
x)
(\partial^{\gamma}u)^{2}(x)dx+C_{2}\mu^{2(|\beta|+|\gamma|-l)}\parallel
u\parallel_{l}^{2}.
\end{eqnarray*}
Under the current assumptions $|\beta|\neq 0$ which implies that
$\mathrm{supp}\text{
}\partial^{\beta}\eta_{\mu}\subset\{|x|<\mu^{-1}\}$, so applying
the Taylor expansion of $\partial^{\gamma}u$ at $x=0$ with
$|\rho|=l-|\gamma|$ and recalling that $u$ vanishes to all orders
at the origin, yields
\begin{equation}
\int_{\mathbb{R}^{2}}(\partial^{\beta}\eta)^{2}(\mu x)
(\partial^{\gamma}u)^{2}(x)dx\leq C_{3}
\int_{B_{\mu^{-1}}(0)}(\partial^{\gamma}u)^{2}(x)dx\leq
C_{4}\mu^{2(|\gamma|-l)}\parallel u\parallel_{l}^{2}.
\end{equation}
Moreover the case $|\gamma|\geq l$ may be treated by Young's
inequality:
\begin{eqnarray*}
\parallel\partial^{\beta}\eta_{\mu}\partial^{\gamma}(\chi_{\mu}*u)\parallel^{2}&\leq&
C_{5}\mu^{2(|\beta|+|\gamma|-|\tau|)}\int_{\mathbb{R}^{2}}(\partial^{\gamma-\tau}
\chi_{\mu}*\partial^{\tau}u)^{2}(x)dx\\
&\leq&C_{6}\mu^{2(m-l)}\parallel u\parallel^{2}_{l},
\end{eqnarray*}
where $|\tau|=l$.  Therefore $(ii)$ follows once $(i)$ is
established, and $(i)$ is established by similar arguments which
will be omitted here.\par
  We now show $(iii)$.  Let
$|\alpha|=m\leq l$ and observe that
\begin{eqnarray*}
\parallel\partial^{\alpha}(u-\overline{S}_{\mu}u)\parallel&\leq&
\parallel\partial^{\alpha}[(1-\eta_{\mu})u]\parallel+
\parallel\partial^{\alpha}[\eta_{\mu}(u-\chi_{\mu}*u)]\parallel\\
&\leq&\sum_{\beta+\gamma=\alpha}(\parallel\partial^{\beta}(1-\eta_{\mu})\partial^{\gamma}u
\parallel+\parallel\partial^{\beta}\eta_{\mu}\partial^{\gamma}(u-\chi_{\mu}*u)\parallel).
\end{eqnarray*}
According to the standard construction [17],
\begin{eqnarray*}
\parallel\partial^{\beta}\eta_{\mu}\partial^{\gamma}(u-
\chi_{\mu}*u)\parallel&\leq& C_{7}\mu^{|\beta|}\parallel
u-\chi_{\mu}*u\parallel_{|\gamma|}\\
&\leq& C_{8}\mu^{|\beta|}\mu^{|\gamma|-l}\parallel
u\parallel_{l}\\
&=& C_{8}\mu^{m-l}\parallel u\parallel_{l}.
\end{eqnarray*}
Furthermore since $\mathrm{supp }\text{
}\partial^{\beta}(1-\eta_{\mu})\subset\{|x|<\mu^{-1}\}$ we may
apply the same methods used to establish (7.2) to obtain
\begin{equation*}
\parallel\partial^{\beta}(1-\eta_{\mu})\partial^{\gamma}u
\parallel\leq C_{9}\mu^{m-l}\parallel u\parallel_{l}.
\end{equation*}
It follows that $(iii)$ holds.\par
  The desired smoothing operators
on $\Omega$ may be obtained from $\overline{S}_{\mu}$ in the
following way. If $\Omega$ is a bounded convex Lipschitz domain,
then Lemma 5.3 yields an extension operator
$E:H^{m}(\Omega)\rightarrow H^{m}(\mathbb{R}^{2})$.  We then
define smoothing operators
$S_{\mu}:L^{2}(\Omega)\rightarrow\overline{H}_{0}^{\infty}(\Omega)$
by $S_{\mu}u=(\overline{S}_{\mu}Eu)|_{\Omega}$.  As $E$ is
bounded, it is clear that Lemma 5.1 will also hold for $S_{\mu}$.

\begin{center}
\textbf{References}
\end{center}

\noindent[1]\hspace{.06in} Q. Han, \textit{On the isometric
embedding of surfaces with Gauss curvature
changing}\par\hspace{.06in}\textit{sign cleanly}, Comm. Pure Appl.
Math., $\mathbf{58}$ (2005), 285-295.\bigskip

\noindent[2]\hspace{.06in} Q. Han, \textit{Local isometric
embedding of surfaces with Gauss curvature
changing}\par\hspace{.06in}\textit{sign stably across a curve},
Calc. Var. \& P.D.E., $\mathbf{25}$ (2006), 79-103.\bigskip

\noindent[3]\hspace{.06in} Q. Han, \textit{Smooth local isometric
embedding of surfaces with Gauss
curvature}\par\hspace{.06in}\textit{changing sign cleanly},
preprint.\bigskip

\noindent[4]\hspace{.06in} Q. Han, \& J.-X. Hong,
\textit{Isometric Embedding of Riemannian Manifolds in
Eu-}\par\hspace{.06in}\textit{clidean Spaces}, Mathematical
Surveys and Monographs, Vol. 130, AMS,
Provi-\par\hspace{.06in}dence, RI, 2006.\bigskip

\noindent[5]\hspace{.06in} Q. Han, J.-X. Hong, \& C.-S. Lin,
\textit{Local isometric embedding of surfaces
with}\par\hspace{.06in}\textit{nonpositive Gaussian curvature,} J.
Differential Geom., $\mathbf{63}$ (2003), 475-520.\bigskip

\noindent[6]\hspace{.06in} Q. Han, \& M. Khuri, \textit{Smooth
solutions of a class of mixed-type differential
equa-}\par\hspace{.06in}\textit{tions}, preprint.\bigskip

\noindent[7]\hspace{.06in} P. Hartman, \& A. Wintner,
\textit{Gaussian curvature and local embedding}, Amer.
J.\par\hspace{.06in}Math., $\mathbf{73}$ (1951), 876-884.\bigskip

\noindent[8]\hspace{.06in} M. Khuri, \textit{The local isometric
embedding in $\mathbb{R}^{3}$ of two-dimensional
Riemannian}\par\hspace{.06in}\textit{manifolds with Gaussian
curvature changing sign to finite order on a
curve},\par\hspace{.06in}J. Differential Geom., $\mathbf{76}$
(2007), 249-291.\bigskip

\noindent[9]\hspace{.06in} M. Khuri, \textit{Local solvability of
degenerate Monge-Amp\`{e}re equations and
applica-}\par\hspace{.06in}\textit{tions to geometry}, Electron. J.
Diff. Eqns., $\mathbf{2007}$ (2007), No. 65, 1-37.\bigskip

\noindent[10] M. Khuri, \textit{Counterexamples to the local
solvability of Monge-Amp\`{e}re
equations}\par\hspace{.06in}\textit{in the plane,} Comm. PDE,
$\mathbf{32}$ (2007), 665-674.\bigskip

\noindent[11]  C.-S. Lin, \textit{The local isometric embedding in
$\mathbb{R}^{3}$ of} 2-\textit{dimensional
Riemannian}\par\hspace{.06in}\textit{manifolds with nonnegative
curvature}, J. Differential Geom., $\mathbf{21}$ (1985), no.
2,\par\hspace{.02in} 213-230.\bigskip

\noindent[12]  C.-S. Lin, \textit{The local isometric embedding in
$\mathbb{R}^{3}$ of two-dimensional
Riemannian}\par\hspace{.06in}\textit{manifolds with Gaussian
curvature changing sign cleanly}, Comm. Pure
Appl.\par\hspace{.06in}Math., $\mathbf{39}$ (1986), no. 6,
867-887.\bigskip

\noindent[13]  N. Nadirashvili, \& Y. Yuan, \textit{Improving
Pogorelov's isometric embedding
couter-}\par\hspace{.06in}\textit{example,} Calc. Var. \& P.D.E.,
$\mathbf{21}$ (2008), no. 3, 319-323.\bigskip

\noindent[14]  O. Oleinik, \textit{On the Cauchy problem for weakly
hyperbolic equations}, Comm.\par\hspace{.01in} Pure Appl. Math.,
$\mathbf{23}$ (1970), 569-586.\bigskip

\noindent[15]  A. Pogorelov, \textit{An example of a two-dimensional
Riemannian metric not admitting}\par\textit{ a local realization in
$E_{3}$}, Dokl. Akad. Nauk. USSR, $\mathbf{198}$ (1971),
42-43.\bigskip

\noindent[16]  L. Schlaefli, \textit{Nota alla memoria del Sig.
Beltrami, Sugli spazii di curvatura}\par\textit{ constante,}  Ann.
di mat., 2e s\'{e}rie, $\mathbf{5}$ (1873), 170-193.\bigskip

\noindent[17]  J. Schwartz, \textit{Nonlinear Functional Analysis},
New York University, New York,\par\text{ }1964.\bigskip

\noindent[18]  E. Stein, \textit{Singular Integrals and
Differentiability Properties of Functions},\par\text{ }Princeton
University Press, Princeton, 1970.\bigskip

\noindent[19]  M. Taylor, \textit{Partial Differential Equations}
III, Springer-Verlag, New York,\par\text{ }1996.\bigskip

\noindent[20]  J. Weingarten, \textit{\"{U}ber die theorie der
Aubeinander abwickelbarren Oberfl\"{a}chen},\par\text{ }Berlin,
1884.

\bigskip\bigskip\footnotesize

\noindent\textsc{Department of Mathematics, University of Notre
Dame, Notre Dame, IN 46556}\par

\noindent\textit{E-mail address}: \verb"qhan@nd.edu"\bigskip

\noindent\textsc{Department of Mathematics, Stony Brook
University, Stony Brook, NY 11794}\par

\noindent\textit{E-mail address}: \verb"khuri@math.sunysb.edu"

\end{document}